\documentclass[11pt,letterpaper]{article}
\usepackage[T1]{fontenc}
\usepackage{lmodern}
\usepackage[margin=1in]{geometry}
\usepackage{amsmath,amssymb,amsthm,mathtools}
\usepackage{microtype}
\usepackage{indentfirst}
\usepackage[shortlabels]{enumitem}
\setlist[itemize]{leftmargin=*,itemsep=3pt,topsep=5pt}
\usepackage{needspace}
\usepackage{titlesec}
\titleformat{\subsection}
  {\normalfont\fontsize{12}{14}\selectfont\bfseries}
  {\thesubsection}{.75em}{}
\titleformat{\subsubsection}
  {\normalfont\normalsize\itshape}
  {\thesubsubsection}{.75em}{}
\titlespacing{\subsection}{0pt}{3ex plus 1ex minus .2ex}{1.2ex plus .2ex}
\titlespacing{\subsubsection}{0pt}{1.8ex plus .5ex minus .2ex}{.7ex plus .2ex}
\renewcommand{\thesubsubsection}{\thesubsection.\alph{subsubsection}}
\usepackage{graphicx,xcolor,tikz}
\usepackage{caption}
\usetikzlibrary{arrows,arrows.meta,decorations.markings}
\tikzset{>=latex,inner sep=3pt}
\usepackage[all]{xy}
\usepackage{mathrsfs,upgreek,stmaryrd}
\usepackage{cite,etoolbox}
\apptocmd{\thebibliography}{\small\raggedright\setlength{\itemsep}{4pt}}{}{}
\usepackage[colorlinks=true,linkcolor=blue!45!black,
  citecolor=blue!45!black,urlcolor=blue!45!black,
  bookmarksnumbered=true]{hyperref}
\hypersetup{pdftitle={Koszul duality and homotopy algebras for operads over quivers},
  pdfauthor={Hang Yuan}}

\numberwithin{equation}{section}
\theoremstyle{plain}
\newcounter{mainthm}
\newtheorem{main}[mainthm]{Theorem}

\newtheorem{thm}{Theorem}[section]
\newtheorem{lem}[thm]{Lemma}
\newtheorem{prop}[thm]{Proposition}
\newtheorem{cor}[thm]{Corollary}
\newtheorem{defn-thm}[thm]{Definition-Theorem}

\newtheorem{defn-lem}[thm]{Definition-Lemma}

\theoremstyle{definition}
\newtheorem{ex}[thm]{Example}

\newtheorem{defn}[thm]{Definition}
\newtheorem{defn-prop}[thm]{Definition-Proposition}
\theoremstyle{remark}
\newtheorem{rmk}[thm]{Remark}
\newtheorem{axiom*}[thm]{Axiom}

\newcommand{\address}[1]{\AfterEndEnvironment{thebibliography}{\par\bigskip{\small\itshape #1\par}}}

\newcommand{\antishriek}{%
	\mathord{\mathchoice
		{\rotatebox[origin=c]{180}{$\displaystyle !$}}
		{\rotatebox[origin=c]{180}{$\textstyle !$}}
		{\rotatebox[origin=c]{180}{$\scriptstyle !$}}
		{\rotatebox[origin=c]{180}{$\scriptscriptstyle !$}}}}

\newcommand{\Hom}{\operatorname{Hom}}
\newcommand{\id}{\operatorname{id}}

\newcommand{\ev}{\mathrm{ev}}

\newcommand{\f}{\mathfrak f}

\newcommand{\p}{\mathfrak p}
\newcommand{\q}{\mathfrak q}

\DeclareMathOperator{\Ob}{Ob}

\begin{document}
\title
{Koszul duality and homotopy algebras for operads over quivers}
\author{Hang Yuan}
\date{}
\address{Beijing Institute of Mathematical Sciences and Applications (BIMSA), Beijing, 101408, China; E-mail: yuanhang@bimsa.cn}
\maketitle
\begin{abstract}
For a quiver $Q$, we introduce and study \textit{$Q$-operads} (or \textit{operads over $Q$}), defined as a special sort of $\mathbf{fc}$-multicategories or virtual double categories, in which the directions of the arrows of \(Q\) encode the permitted input orderings.
In this higher operadic framework, we construct the bar–cobar adjunction and develop a theory of Koszul duality and homotopy algebras, together with compatibility results under changes of quiver.
We treat dg categories and their left modules, right modules, and bimodules as algebras over quiver-operads. We prove that these operads are Koszul and that their resolutions govern the corresponding \(A_\infty\) categories and modules.
We also categorify Sagave's notion of derived \(A_\infty\) algebras and introduce modules over derived \(A_\infty\) algebras and categories.
Finally, we explain that restriction along a folding of quivers constructs diagonal \(A_\infty\) bimodules, while extension along an inclusion of quivers realizes two-sided bar models of derived tensor products.
\end{abstract}
\begingroup
\scriptsize
\tableofcontents
\endgroup

\section{Introduction}

Homotopy algebras have their roots in Stasheff's work on homotopy associativity \cite{stasheff1963homotopy} and the operadic framework developed by Boardman--Vogt and May \cite{boardman1973homotopy,may2006geometry}. In parallel, Priddy introduced Koszul resolutions for quadratic associative algebras \cite{priddy1970koszul}. Koszul duality was subsequently extended to operads by Ginzburg--Kapranov \cite{ginzburg1994koszul} and further developed by Getzler--Jones \cite{getzler1994operads}, connecting these algebraic and homotopical perspectives. For a systematic treatment, see \cite{loday2012algebraic}.

Koszul duality has since developed in several directions, including representation theory \cite{beilinson1996koszul}, colored operads \cite{vanderlaan2003coloured}, and structures with several inputs and outputs, treated through dioperads \cite{gan2003koszul} and properads \cite{Vallette_Koszul_props}. Extensions beyond homogeneous quadratic relations include the quadratic--linear theory used for Batalin--Vilkovisky algebras \cite{GalvezCarrillo_Tonks_Vallette_BV} and the curved theory of Hirsh--Mill\`es, which allows constant terms \cite{hirsh2012curved,hirsh2024correction}.
The literature on homotopy algebras includes strongly homotopy Lie algebras \cite{LadaMarkl1995}, open--closed homotopy algebras and their colored Koszul theory \cite{kajiura2006homotopy,Hoefel_Livernet_OCHA_SHLP}, homotopy Batalin--Vilkovisky algebras \cite{GalvezCarrillo_Tonks_Vallette_BV}, and homotopy diagrams of algebras \cite{Markl_Homotopy_Diagrams}. Sagave's derived $A_\infty$-algebras \cite{sagave2007dg} and their operadic treatment by Livernet--Roitzheim--Whitehouse \cite{Livernet_Roitzheim_Whitehouse_derived_Ainfinity} extend this picture to a bigraded setting.

In this paper, we develop Koszul duality over quivers as a natural framework for categorifying homotopy algebras governed by nonsymmetric operads. Categories and their modules fit naturally into the quiver setting, where incidence records sources, targets, input order, and module positions. We show that the quiver-operads governing dg categories and their modules have Koszul resolutions whose algebras are precisely the usual $A_\infty$ categories and their left modules, right modules, and bimodules. We further extend Sagave's derived $A_\infty$ algebras to categories and introduce bimodules over pairs of derived $A_\infty$ algebras and categories. We also apply restriction and extension along quiver morphisms to construct diagonal bimodules and realize the two-sided bar tensor product as extension of algebras over Koszul resolutions.

A \textit{$Q$-operad} is a nonsymmetric operad colored by the arrows of a quiver $Q$, supported on profiles $(e_1,\ldots,e_n;e_0)$ for which $(e_1,\ldots,e_n)$ is a composable path with the same endpoints as $e_0$. The terminal quiver and a one-vertex quiver with several loops recover nonsymmetric operads and nonsymmetric colored operads, respectively. In the language of higher operads, these are also called \textit{vertically discrete $\mathbf{fc}$-multicategories}, or \textit{virtual double categories} \cite{leinster1999fc,leinster2004higher,cruttwell2009unified}.

%
%The left-module example opening \S~\ref{s_operads_over} illustrates the input rule with a two-arrow quiver. Compositions are encoded by planar rooted trees or their dual diagrams, with operations on vertices or faces, respectively (Figure~\ref{figure:tree_to_planar_diagram}). To form the dual, embed the tree in a disk with its external endpoints on the boundary and place a dual vertex in each complementary region; see \S~\ref{ss_rooted_planar_diagram}.

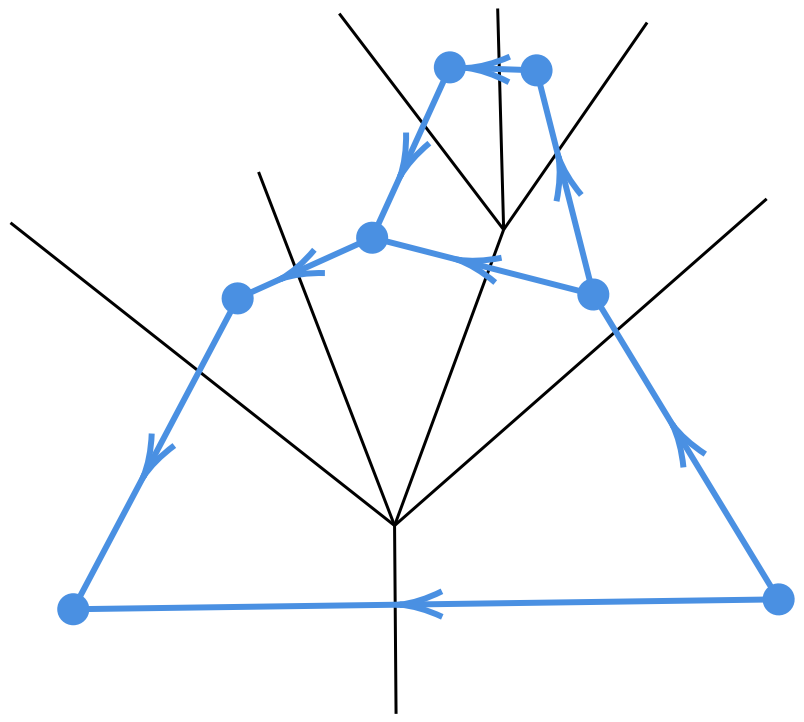
\begin{figure}[htbp]
	\centering

	\scalebox{0.75}{%
	\begin{tikzpicture}[x=0.75pt,y=0.75pt,yscale=-1,xscale=1]

		\draw [color={rgb, 255:red, 0; green, 0; blue, 0 }  ,draw opacity=1 ][fill={rgb, 255:red, 255; green, 255; blue, 255 }  ,fill opacity=1 ][line width=0.75]    (283.5,135.7) -- (412,237) ;
		\draw [color={rgb, 255:red, 0; green, 0; blue, 0 }  ,draw opacity=1 ][fill={rgb, 255:red, 255; green, 255; blue, 255 }  ,fill opacity=1 ][line width=0.75]    (366.5,118.7) -- (412,237) ;
		\draw [color={rgb, 255:red, 0; green, 0; blue, 0 }  ,draw opacity=1 ][fill={rgb, 255:red, 255; green, 255; blue, 255 }  ,fill opacity=1 ][line width=0.75]    (448.5,138) -- (412,237) ;
		\draw [color={rgb, 255:red, 0; green, 0; blue, 0 }  ,draw opacity=1 ][fill={rgb, 255:red, 255; green, 255; blue, 255 }  ,fill opacity=1 ][line width=0.75]    (536.5,127.7) -- (412,237) ;
		\draw [color={rgb, 255:red, 0; green, 0; blue, 0 }  ,draw opacity=1 ][fill={rgb, 255:red, 255; green, 255; blue, 255 }  ,fill opacity=1 ][line width=0.75]    (393.5,65.7) -- (448.5,138) ;
		\draw [color={rgb, 255:red, 0; green, 0; blue, 0 }  ,draw opacity=1 ][fill={rgb, 255:red, 255; green, 255; blue, 255 }  ,fill opacity=1 ][line width=0.75]    (496.5,68.7) -- (448.5,138) ;
		\draw [color={rgb, 255:red, 0; green, 0; blue, 0 }  ,draw opacity=1 ][fill={rgb, 255:red, 255; green, 255; blue, 255 }  ,fill opacity=1 ][line width=0.75]    (446.5,64) -- (448.5,138) ;
		\draw [color={rgb, 255:red, 0; green, 0; blue, 0 }  ,draw opacity=1 ][fill={rgb, 255:red, 255; green, 255; blue, 255 }  ,fill opacity=1 ][line width=0.75]    (412,237) -- (412.5,300) ;
		\draw [color={rgb, 255:red, 74; green, 144; blue, 226 }  ,draw opacity=1 ][line width=1.5]    (304.5,265) -- (359.5,161) ;
		\draw [shift={(359.5,161)}, rotate = 297.87] [color={rgb, 255:red, 74; green, 144; blue, 226 }  ,draw opacity=1 ][fill={rgb, 255:red, 74; green, 144; blue, 226 }  ,fill opacity=1 ][line width=1.5]      (0, 0) circle [x radius= 4.36, y radius= 4.36]   ;
		\draw [shift={(327.89,220.78)}, rotate = 297.87] [color={rgb, 255:red, 74; green, 144; blue, 226 }  ,draw opacity=1 ][line width=1.5]    (14.21,-4.28) .. controls (9.04,-1.82) and (4.3,-0.39) .. (0,0) .. controls (4.3,0.39) and (9.04,1.82) .. (14.21,4.28)   ;
		\draw [color={rgb, 255:red, 74; green, 144; blue, 226 }  ,draw opacity=1 ][line width=1.5]    (304.5,265) -- (540.5,261.7) ;
		\draw [shift={(413.7,263.47)}, rotate = 359.20] [color={rgb, 255:red, 74; green, 144; blue, 226 }  ,draw opacity=1 ][line width=1.5]    (14.21,-4.28) .. controls (9.04,-1.82) and (4.3,-0.39) .. (0,0) .. controls (4.3,0.39) and (9.04,1.82) .. (14.21,4.28)   ;
		\draw [shift={(304.5,265)}, rotate = 359.2] [color={rgb, 255:red, 74; green, 144; blue, 226 }  ,draw opacity=1 ][fill={rgb, 255:red, 74; green, 144; blue, 226 }  ,fill opacity=1 ][line width=1.5]      (0, 0) circle [x radius= 4.36, y radius= 4.36]   ;
		\draw [color={rgb, 255:red, 74; green, 144; blue, 226 }  ,draw opacity=1 ][line width=1.5]    (478.5,159.7) -- (540.5,261.7) ;
		\draw [shift={(540.5,261.7)}, rotate = 58.71] [color={rgb, 255:red, 74; green, 144; blue, 226 }  ,draw opacity=1 ][fill={rgb, 255:red, 74; green, 144; blue, 226 }  ,fill opacity=1 ][line width=1.5]      (0, 0) circle [x radius= 4.36, y radius= 4.36]   ;
		\draw [shift={(504.93,203.18)}, rotate = 58.71] [color={rgb, 255:red, 74; green, 144; blue, 226 }  ,draw opacity=1 ][line width=1.5]    (14.21,-4.28) .. controls (9.04,-1.82) and (4.3,-0.39) .. (0,0) .. controls (4.3,0.39) and (9.04,1.82) .. (14.21,4.28)   ;
		\draw [color={rgb, 255:red, 74; green, 144; blue, 226 }  ,draw opacity=1 ][line width=1.5]    (404.5,140.7) -- (478.5,159.7) ;
		\draw [shift={(478.5,159.7)}, rotate = 14.4] [color={rgb, 255:red, 74; green, 144; blue, 226 }  ,draw opacity=1 ][fill={rgb, 255:red, 74; green, 144; blue, 226 }  ,fill opacity=1 ][line width=1.5]      (0, 0) circle [x radius= 4.36, y radius= 4.36]   ;
		\draw [shift={(432.98,148.01)}, rotate = 14.40] [color={rgb, 255:red, 74; green, 144; blue, 226 }  ,draw opacity=1 ][line width=1.5]    (14.21,-4.28) .. controls (9.04,-1.82) and (4.3,-0.39) .. (0,0) .. controls (4.3,0.39) and (9.04,1.82) .. (14.21,4.28)   ;
		\draw [shift={(404.5,140.7)}, rotate = 14.4] [color={rgb, 255:red, 74; green, 144; blue, 226 }  ,draw opacity=1 ][fill={rgb, 255:red, 74; green, 144; blue, 226 }  ,fill opacity=1 ][line width=1.5]      (0, 0) circle [x radius= 4.36, y radius= 4.36]   ;
		\draw [color={rgb, 255:red, 74; green, 144; blue, 226 }  ,draw opacity=1 ][line width=1.5]    (359.5,161) -- (404.5,140.7) ;
		\draw [shift={(373.98,154.47)}, rotate = 335.72] [color={rgb, 255:red, 74; green, 144; blue, 226 }  ,draw opacity=1 ][line width=1.5]    (14.21,-4.28) .. controls (9.04,-1.82) and (4.3,-0.39) .. (0,0) .. controls (4.3,0.39) and (9.04,1.82) .. (14.21,4.28)   ;
		\draw [color={rgb, 255:red, 74; green, 144; blue, 226 }  ,draw opacity=1 ][line width=1.5]    (459.5,84.7) -- (478.5,159.7) ;
		\draw [shift={(466.84,113.67)}, rotate = 75.78] [color={rgb, 255:red, 74; green, 144; blue, 226 }  ,draw opacity=1 ][line width=1.5]    (14.21,-4.28) .. controls (9.04,-1.82) and (4.3,-0.39) .. (0,0) .. controls (4.3,0.39) and (9.04,1.82) .. (14.21,4.28)   ;
		\draw [color={rgb, 255:red, 74; green, 144; blue, 226 }  ,draw opacity=1 ][line width=1.5]    (430.5,83.7) -- (459.5,84.7) ;
		\draw [shift={(459.5,84.7)}, rotate = 1.97] [color={rgb, 255:red, 74; green, 144; blue, 226 }  ,draw opacity=1 ][fill={rgb, 255:red, 74; green, 144; blue, 226 }  ,fill opacity=1 ][line width=1.5]      (0, 0) circle [x radius= 4.36, y radius= 4.36]   ;
		\draw [shift={(436.21,83.9)}, rotate = 1.97] [color={rgb, 255:red, 74; green, 144; blue, 226 }  ,draw opacity=1 ][line width=1.5]    (14.21,-4.28) .. controls (9.04,-1.82) and (4.3,-0.39) .. (0,0) .. controls (4.3,0.39) and (9.04,1.82) .. (14.21,4.28)   ;
		\draw [shift={(430.5,83.7)}, rotate = 1.97] [color={rgb, 255:red, 74; green, 144; blue, 226 }  ,draw opacity=1 ][fill={rgb, 255:red, 74; green, 144; blue, 226 }  ,fill opacity=1 ][line width=1.5]      (0, 0) circle [x radius= 4.36, y radius= 4.36]   ;
		\draw [color={rgb, 255:red, 74; green, 144; blue, 226 }  ,draw opacity=1 ][line width=1.5]    (404.5,140.7) -- (430.5,83.7) ;
		\draw [shift={(413.85,120.21)}, rotate = 294.52] [color={rgb, 255:red, 74; green, 144; blue, 226 }  ,draw opacity=1 ][line width=1.5]    (14.21,-4.28) .. controls (9.04,-1.82) and (4.3,-0.39) .. (0,0) .. controls (4.3,0.39) and (9.04,1.82) .. (14.21,4.28)   ;

	\end{tikzpicture}
	}%
	\caption{A planar tree (black) and its dual diagram (blue), with external endpoints understood to lie on a surrounding disk; compare Keller's survey \cite[\S 2.2]{keller2005infinity}.}
	\label{figure:tree_to_planar_diagram}
\end{figure}

\subsection*{General theory}
We work over a field $\Bbbk$ with nonsymmetric operads, allowing unary operations and retaining operadic identities. Differentials have cohomological degree $+1$, and tensor permutations follow the Koszul sign rule. The associative algebras, categories, and modules in our Koszul examples are non-unital. We impose strict units locally for the derived tensor products in \S~\ref{sss_derived_tensor_construction}. The Koszul statements themselves do not include strictly unital extensions.

The bar--cobar adjunction for $Q$-operads (Theorem~\ref{thm_bar_cobar_adjunction}) provides the starting point for our construction of homotopy algebras. The notion of \textit{Koszulness} extends to dg quadratic quiver-operads (Definition~\ref{defn_koszul_fc_multicategory}). Let $\mathscr M$ be a quadratic Koszul dg $Q$-operad arising from a dg quadratic datum, and write $\mathscr M^{\antishriek}$ for its Koszul dual dg cooperad. The differential on the generators is allowed to be nonzero and preserves the quadratic relations. The canonical map $\mathscr M_\infty=\Omega_Q\mathscr M^{\antishriek}\to\mathscr M$ is called its \textit{Koszul resolution}. Algebras over $\mathscr M_\infty$ are called \textit{homotopy $\mathscr M$-algebras}.

The familiar constructions of endomorphism operads and cofree coalgebras carry over to the quiver-operadic setting. For a fixed family $X=\{X(e)\}_{e\in Q}$ of cochain complexes indexed by the arrows of $Q$, write $\mathbf{End}_Q(X)$ for its endomorphism $Q$-operad and $\mathbb T^{\mathscr M^{\antishriek}}(X)$ for the cofree $\mathscr M^{\antishriek}$-coalgebra on $X$.

\begin{main}[Homotopy algebras; Theorem~\ref{thm_three_homotopy_M_algebra_definitions}]
\label{thm_introduction_three_homotopy_M_algebra_definitions}
The following data are naturally equivalent:
\begin{enumerate}[(i)]
\itemsep 1pt
\item an $\mathscr M_\infty$-algebra structure on $X$;
\item a degree-one twisting map $\theta:\mathscr M^{\antishriek}\to\mathbf{End}_Q(X)$ vanishing on the coaugmentation;
\item a square-zero differential $d=\bar d+\mathcal D$ on $\mathbb T^{\mathscr M^{\antishriek}}(X)$, where $\bar d$ is the canonical differential induced by those of $\mathscr M^{\antishriek}$ and $X$, and $\mathcal D$ is a degree-one coderivation vanishing on the cogenerators $X$.
\end{enumerate}
\end{main}
The total differential $d$ is compatible with the cooperad differential in the sense of Definition~\ref{defn:C_coderivation}. The vanishing of $\mathcal D$ on $X$ ensures that $d$ restricts to the prescribed differential on $X$. The twisting-map description requires no augmentation of $\mathbf{End}_Q(X)$.

Quiver morphisms relate these constructions for different indexing
quivers. Profilewise pullback of operads extends the usual change of
colours \cite[\S~1.6]{BergerMoerdijkColoured}. On algebra categories,
morphisms of dg quiver-operads give the following adjunction.

\begin{main}[Restriction and extension of algebras;
Propositions~\ref{prop_change_quiver_algebras}--\ref{prop_change_quiver_algebra_adjunction}]
\label{thm_introduction_algebra_adjunction}
A morphism of dg quiver-operads
$
 \Phi=(f,\phi):(Q,\mathscr O)\longrightarrow(Q',\mathscr P)$
induces a restriction functor $\Phi^*$ admitting a left adjoint
$\Phi_!$. Thus there are bijections, natural in
$X\in\operatorname{Alg}(\mathscr O)$ and
$Y\in\operatorname{Alg}(\mathscr P)$,
\[
 \operatorname{Hom}_{\operatorname{Alg}(\mathscr P)}(\Phi_!X,Y)
 \cong
 \operatorname{Hom}_{\operatorname{Alg}(\mathscr O)}(X,\Phi^*Y).
\]
%Restriction reindexes the underlying complexes along $f$ and
%restricts the operations along $\phi$.
\end{main}
For the endpoint-quiver inclusion of \S~\ref{sss_derived_tensor_construction},
the new-arrow component of extension is the ordinary balanced tensor
product; over the Koszul resolutions, it is the two-sided bar model of the
derived tensor product (Proposition~\ref{prop_tensor_algebra_extension}).

Compatibility of pullback with cooperads and Koszul resolutions
holds under the \textit{unique
diagram lifting} condition of
Definition~\ref{defn_unique_diagram_lifting}. This condition holds
in many cases of interest (see
Lemma~\ref{lem_change_quiver_thin}).

\begin{main}[Change of quiver; Theorem~\ref{thm_change_quiver_koszul}]
\label{thm_introduction_change_quiver}
Let $f:Q\to Q'$ have unique diagram lifting, and let $\mathscr P$ be
a quadratic Koszul dg $Q'$-operad. Then $f^*\mathscr P$ is a
quadratic Koszul dg $Q$-operad, with canonical identifications
\[
 (f^*\mathscr P)^{\antishriek}
 \cong f^*(\mathscr P^{\antishriek}),
 \qquad
 (f^*\mathscr P)_\infty\cong f^*(\mathscr P_\infty),
\]
compatible with the canonical twisting morphisms and resolution maps.
\end{main}
Proposition~\ref{prop_change_quiver_bar_cobar} supplies
the corresponding bar--cobar compatibility. The internal differentials
of the quadratic data may be nonzero.

\subsection*{Applications to categories and modules}

For the statements below, fix object sets $V,W$. We write $\mathsf{dgCat}_V$ for the quadratic quiver-operad of \S~\ref{ss_example_dg_categories} governing non-unital dg categories with object set $V$. We write $\mathscr B\mathsf{Mod}_{V,W}$ for the quadratic quiver-operad of \S~\ref{sss_dg_category_bimodules} governing pairs of non-unital dg categories with object sets $V,W$, together with a dg bimodule between them.

In each case, $X$ and $Y$ are fixed families of cochain complexes indexed by the arrows of the respective category and bimodule quivers. All structures in the correspondences below use these prescribed underlying complexes. We suppress the quiver subscripts on endomorphism operads; all displayed morphisms are morphisms of dg quiver-operads. The categories in the following two theorems are non-unital.

\begin{main}[Categories and modules]
The quiver-operads $\mathsf{dgCat}_V$ and $\mathscr B\mathsf{Mod}_{V,W}$ are Koszul. Their Koszul resolutions have the following interpretations:
\begin{enumerate}[(i)]
\itemsep 1pt
\item An $A_\infty$ category structure with object set $V$ is equivalent to a $(\mathsf{dgCat}_V)_\infty$-algebra structure on $X$, that is, a morphism
\[
\alpha:(\mathsf{dgCat}_V)_\infty
\longrightarrow \mathbf{End}(X).
\]
\item A pair of $A_\infty$ categories with object sets $V,W$, together with an $A_\infty$ bimodule over this pair, is equivalent to a $(\mathscr B\mathsf{Mod}_{V,W})_\infty$-algebra structure on $Y$, that is, a morphism
\[
\beta:(\mathscr B\mathsf{Mod}_{V,W})_\infty
\longrightarrow \mathbf{End}(Y).
\]
\end{enumerate}
\end{main}
The analogous assertions hold for the left- and right-module quiver-operads. The algebra and module cases are treated in \S~\ref{ss_example_left_module}; the categorical versions are treated in \S~\ref{ss_example_dg_categories}. The proofs identify the relevant complexes profilewise with the associative resolution, with the required object or module colors.

The derived application starts from Sagave's derived $A_\infty$ algebras \cite{sagave2007dg} and their operadic description by Livernet--Roitzheim--Whitehouse \cite{Livernet_Roitzheim_Whitehouse_derived_Ainfinity}. We extend this construction to categories and bimodules by including horizontal degrees in the quiver data. The strict structures are bidg-categories, whose morphism spaces carry compatible horizontal and vertical differentials.

We write $\mathsf{bidgCat}_V$ for the quadratic quiver-operad of \S~\ref{sss_bidg_categories} governing non-unital bidg-categories with object set $V$. We write $\mathsf{bidgBMod}_{V,W}$ for the quadratic quiver-operad of \S~\ref{sss_derived_Ainfinity_category_bimodules} governing pairs of non-unital bidg-categories with object sets $V,W$, together with a bidg-bimodule between them.

\begin{main}[Derived $A_\infty$ categories and modules]
The quiver-operads $\mathsf{bidgCat}_V$ and $\mathsf{bidgBMod}_{V,W}$ are Koszul. Their Koszul resolutions have the following interpretations:
\begin{enumerate}[(i)]
\itemsep 1pt
\item A derived $A_\infty$ category structure with object set $V$ is equivalent to a $(\mathsf{bidgCat}_V)_\infty$-algebra structure on $X$, that is, a morphism
\[
\alpha:(\mathsf{bidgCat}_V)_\infty
\longrightarrow \mathbf{End}(X).
\]
\item A pair of derived $A_\infty$ categories with object sets $V,W$, together with a derived $A_\infty$ bimodule over this pair of categories, is equivalent to a $(\mathsf{bidgBMod}_{V,W})_\infty$-algebra structure on $Y$, that is, a morphism
\[
\beta:(\mathsf{bidgBMod}_{V,W})_\infty
\longrightarrow \mathbf{End}(Y).
\]
\end{enumerate}
\end{main}
We describe these resolutions in \S~\ref{s_derived_Ainfinity_categories}. The bimodule construction extends the two-sided representation viewpoint of \cite{Aponte_Livernet_Robertson_Whitehouse_Ziegenhagen_representations} to pairs of algebras and categories.

Finally, we illustrate Theorem~\ref{thm_introduction_algebra_adjunction}
with two examples: diagonal bimodules and tensor products.
The folding quiver morphism $f:Q^{\mathrm{bi}}\to Q_V$ identifies
the two copies of $V$ and forgets the copy labels on arrows
(\S~\ref{sss_diagonal_bimodule_construction},
Figure~\ref{fig:diagonal-quiver-folding}).
The quiver morphism $i:Q^-\hookrightarrow Q^+$ adds the composite
arrow $e_{ac}:a\to c$ to the quiver for right and left modules
over a common category
(\S~\ref{sss_derived_tensor_construction},
Figure~\ref{fig:bar-tensor-quiver}).

\begin{main}
\leavevmode
\begin{enumerate}[(i),itemsep=6pt]
\item \textup{\textbf{Diagonal bimodules.}}
For $A_\infty$ categories on a fixed object set, the functor
sending $\mathcal C$ to its two labelled copies together with
the diagonal $A_\infty$ bimodule $\mathcal C_\Delta$ is right adjoint to
extension along the folding quiver morphism $f$
(Corollary~\ref{cor_diagonal_algebra_adjunction}).

\item \textup{\textbf{Tensor products.}}
Extension along the endpoint-quiver inclusion $i$ is left adjoint
to forgetting the added component and its pairing operations.
For strictly unital data $(\mathcal B,M,N)$, it retains
$\mathcal B$, its right module $M$, and its left module $N$.
The added component is the ordinary balanced tensor product
$M\otimes_{\mathcal B}N$ in the dg case and the two-sided bar
model $M\otimes_{\mathcal B}^{\mathbb L}N$ over the Koszul
resolutions
(Corollary~\ref{cor_tensor_algebra_adjunction}).
\end{enumerate}
\end{main}
We refer to \S~\ref{sss_diagonal_bimodule_pullback} for further details.

\subsection*{Further work}

We plan to develop deformation cohomology for homotopy algebras over
quivers using their twisting-morphism and coderivation descriptions,
and to compare it with Hochschild cohomology. The diagonal bimodule
and bar tensor-product constructions suggest applications to coherent
gluing of $A_\infty$ categories. A concrete goal is to express
obstructions to extending compatible gluings and their deformations
in terms of derived Hom complexes of bimodules.

Two further directions concern extensions of the quiver-operadic framework.
A cyclic version, incorporating reversal of arrow colors and cyclic
symmetry of profiles, may be motivated by cyclic $A_\infty$ categories,
homotopy inner products, and Calabi--Yau structures
\cite{TradlerInfinityInnerProducts,Brav_Dyckerhoff_relative_CY}.
Certain versions of properads and PROPs over quivers can also be
anticipated from geometric considerations analogous to those in
\cite{Yuan_higher_operad_Fuk}.

\medskip
\noindent\textit{Acknowledgment and AI disclosure.}
The author thanks Xinxing Tang, Yi Wang, and Peng Yang for helpful conversations.
The author contributed the underlying idea of this work inspired by \cite{Yuan_higher_operad_Fuk}. 
AI tools were substantially used for the drafting and revision and for most computations in \S~\ref{s_examples_applications} and \S~\ref{sss_diagonal_bimodule_pullback}.

\section{Operads and cooperads over quivers}
\label{s_operads_over}

To intuitively illustrate how a quiver records the inputs of their operations, we begin with a basic model: an associative algebra $A$ and a left $A$-module $M$.
Let $Q$ have vertices $v_a,v_m$, a loop $e_a:v_a\to v_a$, and an arrow
$e_m:v_a\to v_m$:

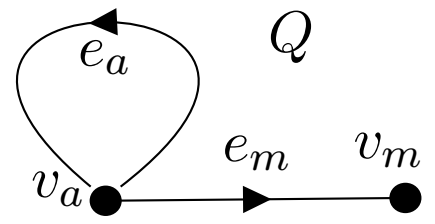
\begin{center}
\begin{tikzpicture}[x=0.6pt,y=0.6pt,yscale=-1,xscale=1]
			%uncomment if require: \path (0,116); %set diagram left start at 0, and has height of 116
			
			%Shape: Ellipse [id:dp47198104618365977] 
			\draw  [fill={rgb, 255:red, 0; green, 0; blue, 0 }  ,fill opacity=1 ] (374.53,88.7) .. controls (374.53,86.25) and (376.63,84.27) .. (379.23,84.27) .. controls (381.83,84.27) and (383.93,86.25) .. (383.93,88.7) .. controls (383.93,91.15) and (381.83,93.14) .. (379.23,93.14) .. controls (376.63,93.14) and (374.53,91.15) .. (374.53,88.7) -- cycle ;
			%Straight Lines [id:da9168862974323535] 
			\draw    (379.23,88.7) -- (473.25,87.82) ;
			\draw [shift={(431.24,88.21)}, rotate = 179.46] [fill={rgb, 255:red, 0; green, 0; blue, 0 }  ][line width=0.08]  [draw opacity=0] (8.93,-4.29) -- (0,0) -- (8.93,4.29) -- cycle    ;
			%Curve Lines [id:da3908907017334565] 
			\draw    (374.53,84.27) .. controls (289.91,14.15) and (471.37,16.81) .. (383.93,84.27) ;
			\draw [shift={(373.89,33.01)}, rotate = 357.06] [fill={rgb, 255:red, 0; green, 0; blue, 0 }  ][line width=0.08]  [draw opacity=0] (8.93,-4.29) -- (0,0) -- (8.93,4.29) -- cycle    ;
			%Shape: Ellipse [id:dp9844218967159158] 
			\draw  [fill={rgb, 255:red, 0; green, 0; blue, 0 }  ,fill opacity=1 ] (468.55,87.82) .. controls (468.55,85.37) and (470.65,83.38) .. (473.25,83.38) .. controls (475.85,83.38) and (477.95,85.37) .. (477.95,87.82) .. controls (477.95,90.27) and (475.85,92.25) .. (473.25,92.25) .. controls (470.65,92.25) and (468.55,90.27) .. (468.55,87.82) -- cycle ;
			
			% Text Node
			\draw (354.22,78.06) node [anchor=north west][inner sep=0.75pt]    {$v_{a}$};
			% Text Node
			\draw (414.27,66.8) node [anchor=north west][inner sep=0.75pt]    {$e_{m}$};
			% Text Node
			\draw (368.97,36.88) node [anchor=north west][inner sep=0.75pt]    {$e_{a}$};
			% Text Node
			\draw (455.35,65.9) node [anchor=north west][inner sep=0.75pt]    {$v_{m}$};
			% Text Node
			\draw (428.58,27.86) node [anchor=north west][inner sep=0.75pt]    {$Q$};
		\end{tikzpicture}
\end{center}

Since no arrow leaves $v_m$, a composable path can contain $e_m$ only
as its last arrow. The admissible profiles therefore have the forms
$(e_a,\ldots,e_a;e_a)$ and $(e_a,\ldots,e_a,e_m;e_m)$, pictured below
together with their dual trees.

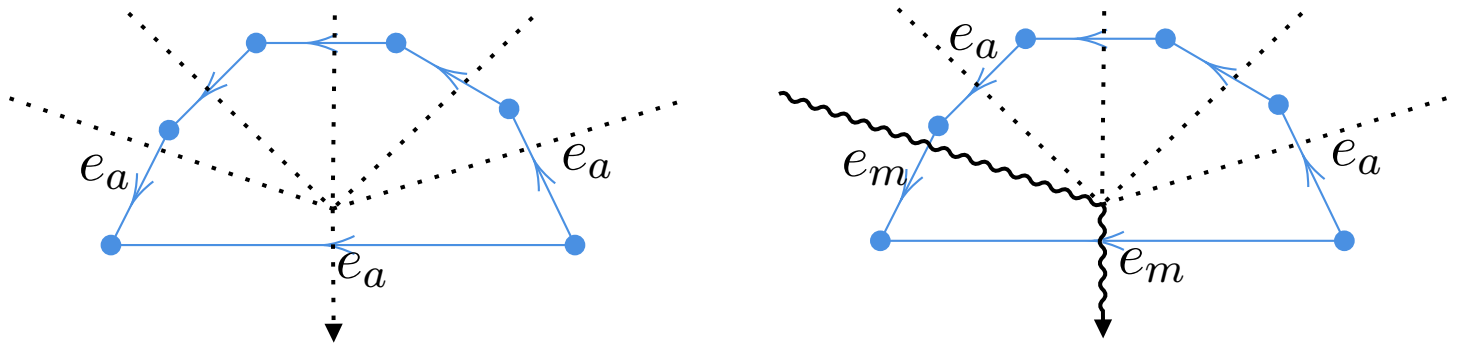
\begin{center}
\begin{tikzpicture}[x=0.5pt,y=0.5pt,yscale=-0.8,xscale=0.8]
		%uncomment if require: \path (0,187); %set diagram left start at 0, and has height of 187
		
		%Straight Lines [id:da2743514629985906] 
		\draw [color={rgb, 255:red, 74; green, 144; blue, 226 }  ,draw opacity=1 ]     (318,64) -- (349,128.07) ;
		\draw [shift={(329.67,88.11)}, rotate = 64.18] [color={rgb, 255:red, 74; green, 144; blue, 226 }  ,draw opacity=1 ]     (14.21,-4.28) .. controls (9.04,-1.82) and (4.3,-0.39) .. (0,0) .. controls (4.3,0.39) and (9.04,1.82) .. (14.21,4.28)   ;
		\draw [shift={(318,64)}, rotate = 64.18] [color={rgb, 255:red, 74; green, 144; blue, 226 }  ,draw opacity=1 ][fill={rgb, 255:red, 74; green, 144; blue, 226 }  ,fill opacity=1 ]       (0, 0) circle [x radius= 4.36, y radius= 4.36]   ;
		%Straight Lines [id:da2458406783735987] 
		\draw [color={rgb, 255:red, 74; green, 144; blue, 226 }  ,draw opacity=1 ]     (130.67,128.07) -- (349,128.07) ;
		\draw [shift={(349,128.07)}, rotate = 0] [color={rgb, 255:red, 74; green, 144; blue, 226 }  ,draw opacity=1 ][fill={rgb, 255:red, 74; green, 144; blue, 226 }  ,fill opacity=1 ]       (0, 0) circle [x radius= 4.36, y radius= 4.36]   ;
		\draw [shift={(231.04,128.07)}, rotate = 0] [color={rgb, 255:red, 74; green, 144; blue, 226 }  ,draw opacity=1 ]     (14.21,-4.28) .. controls (9.04,-1.82) and (4.3,-0.39) .. (0,0) .. controls (4.3,0.39) and (9.04,1.82) .. (14.21,4.28)   ;
		\draw [shift={(130.67,128.07)}, rotate = 0] [color={rgb, 255:red, 74; green, 144; blue, 226 }  ,draw opacity=1 ][fill={rgb, 255:red, 74; green, 144; blue, 226 }  ,fill opacity=1 ]       (0, 0) circle [x radius= 4.36, y radius= 4.36]   ;
		%Straight Lines [id:da6824309545734145] 
		\draw [color={rgb, 255:red, 74; green, 144; blue, 226 }  ,draw opacity=1 ]     (130.67,128.07) -- (158,74) ;
		\draw [shift={(158,74)}, rotate = 296.81] [color={rgb, 255:red, 74; green, 144; blue, 226 }  ,draw opacity=1 ][fill={rgb, 255:red, 74; green, 144; blue, 226 }  ,fill opacity=1 ]       (0, 0) circle [x radius= 4.36, y radius= 4.36]   ;
		\draw [shift={(140.37,108.89)}, rotate = 296.81] [color={rgb, 255:red, 74; green, 144; blue, 226 }  ,draw opacity=1 ]     (14.21,-4.28) .. controls (9.04,-1.82) and (4.3,-0.39) .. (0,0) .. controls (4.3,0.39) and (9.04,1.82) .. (14.21,4.28)   ;
		%Straight Lines [id:da5508826724035706] 
		\draw [color={rgb, 255:red, 74; green, 144; blue, 226 }  ,draw opacity=1 ]     (264.84,33) -- (318,64) ;
		\draw [shift={(283.82,44.07)}, rotate = 30.25] [color={rgb, 255:red, 74; green, 144; blue, 226 }  ,draw opacity=1 ]     (14.21,-4.28) .. controls (9.04,-1.82) and (4.3,-0.39) .. (0,0) .. controls (4.3,0.39) and (9.04,1.82) .. (14.21,4.28)   ;
		\draw [shift={(264.84,33)}, rotate = 30.25] [color={rgb, 255:red, 74; green, 144; blue, 226 }  ,draw opacity=1 ][fill={rgb, 255:red, 74; green, 144; blue, 226 }  ,fill opacity=1 ]       (0, 0) circle [x radius= 4.36, y radius= 4.36]   ;
		%Straight Lines [id:da7636222570512089] 
		\draw [color={rgb, 255:red, 74; green, 144; blue, 226 }  ,draw opacity=1 ]     (158,74) -- (199,33) ;
		\draw [shift={(172.28,59.72)}, rotate = 315] [color={rgb, 255:red, 74; green, 144; blue, 226 }  ,draw opacity=1 ]     (14.21,-4.28) .. controls (9.04,-1.82) and (4.3,-0.39) .. (0,0) .. controls (4.3,0.39) and (9.04,1.82) .. (14.21,4.28)   ;
		%Straight Lines [id:da7950693085877614] 
		\draw [color={rgb, 255:red, 74; green, 144; blue, 226 }  ,draw opacity=1 ]     (199,33) -- (264.84,33) ;
		\draw [shift={(223.12,33)}, rotate = 0] [color={rgb, 255:red, 74; green, 144; blue, 226 }  ,draw opacity=1 ]     (14.21,-4.28) .. controls (9.04,-1.82) and (4.3,-0.39) .. (0,0) .. controls (4.3,0.39) and (9.04,1.82) .. (14.21,4.28)   ;
		\draw [shift={(199,33)}, rotate = 0] [color={rgb, 255:red, 74; green, 144; blue, 226 }  ,draw opacity=1 ][fill={rgb, 255:red, 74; green, 144; blue, 226 }  ,fill opacity=1 ]       (0, 0) circle [x radius= 4.36, y radius= 4.36]   ;
		%Straight Lines [id:da483902969481506] 
		\draw [color={rgb, 255:red, 0; green, 0; blue, 0 }  ,draw opacity=1 ][fill={rgb, 255:red, 255; green, 255; blue, 255 }  ,fill opacity=1 ][line width=0.75]  [dash pattern={on 0.84pt off 2.51pt}]  (139,19) -- (235,111) ;
		%Straight Lines [id:da9043853049072024] 
		\draw [color={rgb, 255:red, 0; green, 0; blue, 0 }  ,draw opacity=1 ][fill={rgb, 255:red, 255; green, 255; blue, 255 }  ,fill opacity=1 ][line width=0.75]  [dash pattern={on 0.84pt off 2.51pt}]  (236,20) -- (235,111) ;
		%Straight Lines [id:da6006333222259239] 
		\draw [color={rgb, 255:red, 0; green, 0; blue, 0 }  ,draw opacity=1 ][fill={rgb, 255:red, 255; green, 255; blue, 255 }  ,fill opacity=1 ][line width=0.75]  [dash pattern={on 0.84pt off 2.51pt}]  (329,21) -- (235,111) ;
		%Straight Lines [id:da9346034315245484] 
		\draw [color={rgb, 255:red, 0; green, 0; blue, 0 }  ,draw opacity=1 ][fill={rgb, 255:red, 255; green, 255; blue, 255 }  ,fill opacity=1 ][line width=0.75]  [dash pattern={on 0.84pt off 2.51pt}]  (397,61) -- (235,111) ;
		%Straight Lines [id:da8403787251435826] 
		\draw [color={rgb, 255:red, 0; green, 0; blue, 0 }  ,draw opacity=1 ][fill={rgb, 255:red, 255; green, 255; blue, 255 }  ,fill opacity=1 ][line width=0.75]  [dash pattern={on 0.84pt off 2.51pt}]  (235,111) -- (235.48,171) ;
		\draw [shift={(235.5,174)}, rotate = 269.55] [fill={rgb, 255:red, 0; green, 0; blue, 0 }  ,fill opacity=1 ][line width=0.08]  [draw opacity=0] (8.93,-4.29) -- (0,0) -- (8.93,4.29) -- cycle    ;
		%Straight Lines [id:da20442311028773263] 
		\draw [color={rgb, 255:red, 0; green, 0; blue, 0 }  ,draw opacity=1 ][fill={rgb, 255:red, 255; green, 255; blue, 255 }  ,fill opacity=1 ][line width=0.75]  [dash pattern={on 0.84pt off 2.51pt}]  (83,59) -- (235,111) ;
		%Straight Lines [id:da515530565483363] 
		\draw [color={rgb, 255:red, 74; green, 144; blue, 226 }  ,draw opacity=1 ]     (680,62) -- (711,126.07) ;
		\draw [shift={(691.67,86.11)}, rotate = 64.18] [color={rgb, 255:red, 74; green, 144; blue, 226 }  ,draw opacity=1 ]     (14.21,-4.28) .. controls (9.04,-1.82) and (4.3,-0.39) .. (0,0) .. controls (4.3,0.39) and (9.04,1.82) .. (14.21,4.28)   ;
		\draw [shift={(680,62)}, rotate = 64.18] [color={rgb, 255:red, 74; green, 144; blue, 226 }  ,draw opacity=1 ][fill={rgb, 255:red, 74; green, 144; blue, 226 }  ,fill opacity=1 ]       (0, 0) circle [x radius= 4.36, y radius= 4.36]   ;
		%Straight Lines [id:da11343569772679718] 
		\draw [color={rgb, 255:red, 74; green, 144; blue, 226 }  ,draw opacity=1 ]     (492.67,126.07) -- (711,126.07) ;
		\draw [shift={(711,126.07)}, rotate = 0] [color={rgb, 255:red, 74; green, 144; blue, 226 }  ,draw opacity=1 ][fill={rgb, 255:red, 74; green, 144; blue, 226 }  ,fill opacity=1 ]       (0, 0) circle [x radius= 4.36, y radius= 4.36]   ;
		\draw [shift={(593.04,126.07)}, rotate = 0] [color={rgb, 255:red, 74; green, 144; blue, 226 }  ,draw opacity=1 ]     (14.21,-4.28) .. controls (9.04,-1.82) and (4.3,-0.39) .. (0,0) .. controls (4.3,0.39) and (9.04,1.82) .. (14.21,4.28)   ;
		\draw [shift={(492.67,126.07)}, rotate = 0] [color={rgb, 255:red, 74; green, 144; blue, 226 }  ,draw opacity=1 ][fill={rgb, 255:red, 74; green, 144; blue, 226 }  ,fill opacity=1 ]       (0, 0) circle [x radius= 4.36, y radius= 4.36]   ;
		%Straight Lines [id:da7408323817023215] 
		\draw [color={rgb, 255:red, 74; green, 144; blue, 226 }  ,draw opacity=1 ]     (492.67,126.07) -- (520,72) ;
		\draw [shift={(520,72)}, rotate = 296.81] [color={rgb, 255:red, 74; green, 144; blue, 226 }  ,draw opacity=1 ][fill={rgb, 255:red, 74; green, 144; blue, 226 }  ,fill opacity=1 ]       (0, 0) circle [x radius= 4.36, y radius= 4.36]   ;
		\draw [shift={(502.37,106.89)}, rotate = 296.81] [color={rgb, 255:red, 74; green, 144; blue, 226 }  ,draw opacity=1 ]     (14.21,-4.28) .. controls (9.04,-1.82) and (4.3,-0.39) .. (0,0) .. controls (4.3,0.39) and (9.04,1.82) .. (14.21,4.28)   ;
		%Straight Lines [id:da47657487153338407] 
		\draw [color={rgb, 255:red, 74; green, 144; blue, 226 }  ,draw opacity=1 ]     (626.84,31) -- (680,62) ;
		\draw [shift={(645.82,42.07)}, rotate = 30.25] [color={rgb, 255:red, 74; green, 144; blue, 226 }  ,draw opacity=1 ]     (14.21,-4.28) .. controls (9.04,-1.82) and (4.3,-0.39) .. (0,0) .. controls (4.3,0.39) and (9.04,1.82) .. (14.21,4.28)   ;
		\draw [shift={(626.84,31)}, rotate = 30.25] [color={rgb, 255:red, 74; green, 144; blue, 226 }  ,draw opacity=1 ][fill={rgb, 255:red, 74; green, 144; blue, 226 }  ,fill opacity=1 ]       (0, 0) circle [x radius= 4.36, y radius= 4.36]   ;
		%Straight Lines [id:da25654951253673564] 
		\draw [color={rgb, 255:red, 74; green, 144; blue, 226 }  ,draw opacity=1 ]     (520,72) -- (561,31) ;
		\draw [shift={(534.28,57.72)}, rotate = 315] [color={rgb, 255:red, 74; green, 144; blue, 226 }  ,draw opacity=1 ]     (14.21,-4.28) .. controls (9.04,-1.82) and (4.3,-0.39) .. (0,0) .. controls (4.3,0.39) and (9.04,1.82) .. (14.21,4.28)   ;
		%Straight Lines [id:da439817490286485] 
		\draw [color={rgb, 255:red, 74; green, 144; blue, 226 }  ,draw opacity=1 ]     (561,31) -- (626.84,31) ;
		\draw [shift={(585.12,31)}, rotate = 0] [color={rgb, 255:red, 74; green, 144; blue, 226 }  ,draw opacity=1 ]     (14.21,-4.28) .. controls (9.04,-1.82) and (4.3,-0.39) .. (0,0) .. controls (4.3,0.39) and (9.04,1.82) .. (14.21,4.28)   ;
		\draw [shift={(561,31)}, rotate = 0] [color={rgb, 255:red, 74; green, 144; blue, 226 }  ,draw opacity=1 ][fill={rgb, 255:red, 74; green, 144; blue, 226 }  ,fill opacity=1 ]       (0, 0) circle [x radius= 4.36, y radius= 4.36]   ;
		%Straight Lines [id:da008523490281752255] 
		\draw [color={rgb, 255:red, 0; green, 0; blue, 0 }  ,draw opacity=1 ][fill={rgb, 255:red, 255; green, 255; blue, 255 }  ,fill opacity=1 ][line width=0.75]  [dash pattern={on 0.84pt off 2.51pt}]  (501,17) -- (597,109) ;
		%Straight Lines [id:da05811564003016201] 
		\draw [color={rgb, 255:red, 0; green, 0; blue, 0 }  ,draw opacity=1 ][fill={rgb, 255:red, 255; green, 255; blue, 255 }  ,fill opacity=1 ][line width=0.75]  [dash pattern={on 0.84pt off 2.51pt}]  (598,18) -- (597,109) ;
		%Straight Lines [id:da716207197021065] 
		\draw [color={rgb, 255:red, 0; green, 0; blue, 0 }  ,draw opacity=1 ][fill={rgb, 255:red, 255; green, 255; blue, 255 }  ,fill opacity=1 ][line width=0.75]  [dash pattern={on 0.84pt off 2.51pt}]  (691,19) -- (597,109) ;
		%Straight Lines [id:da8674293461590721] 
		\draw [color={rgb, 255:red, 0; green, 0; blue, 0 }  ,draw opacity=1 ][fill={rgb, 255:red, 255; green, 255; blue, 255 }  ,fill opacity=1 ][line width=0.75]  [dash pattern={on 0.84pt off 2.51pt}]  (759,59) -- (597,109) ;
		%Straight Lines [id:da7108946363859182] 
		\draw [color={rgb, 255:red, 0; green, 0; blue, 0 }  ,draw opacity=1 ][fill={rgb, 255:red, 255; green, 255; blue, 255 }  ,fill opacity=1 ][line width=0.75]    (597,109) .. controls (598.68,110.65) and (598.69,112.32) .. (597.04,114) .. controls (595.39,115.68) and (595.4,117.35) .. (597.08,119) .. controls (598.76,120.65) and (598.77,122.32) .. (597.12,124) .. controls (595.47,125.68) and (595.48,127.35) .. (597.16,129) .. controls (598.84,130.65) and (598.85,132.32) .. (597.2,134) .. controls (595.55,135.68) and (595.56,137.35) .. (597.24,139) .. controls (598.92,140.65) and (598.93,142.32) .. (597.28,144) .. controls (595.63,145.68) and (595.64,147.35) .. (597.32,149) .. controls (599,150.65) and (599.01,152.32) .. (597.36,154) .. controls (595.71,155.68) and (595.72,157.35) .. (597.4,159) -- (597.41,161) -- (597.48,169) ;
		\draw [shift={(597.5,172)}, rotate = 269.55] [fill={rgb, 255:red, 0; green, 0; blue, 0 }  ,fill opacity=1 ][line width=0.08]  [draw opacity=0] (8.93,-4.29) -- (0,0) -- (8.93,4.29) -- cycle    ;
		%Straight Lines [id:da6056787942647931] 
		\draw [color={rgb, 255:red, 0; green, 0; blue, 0 }  ,draw opacity=1 ][fill={rgb, 255:red, 255; green, 255; blue, 255 }  ,fill opacity=1 ][line width=0.75]    (445,57) .. controls (447.12,55.96) and (448.7,56.5) .. (449.73,58.62) .. controls (450.76,60.74) and (452.34,61.28) .. (454.46,60.24) .. controls (456.58,59.2) and (458.16,59.74) .. (459.19,61.86) .. controls (460.23,63.97) and (461.81,64.51) .. (463.92,63.47) .. controls (466.04,62.43) and (467.62,62.97) .. (468.65,65.09) .. controls (469.68,67.21) and (471.26,67.75) .. (473.38,66.71) .. controls (475.5,65.67) and (477.08,66.21) .. (478.12,68.33) .. controls (479.15,70.45) and (480.73,70.99) .. (482.85,69.95) .. controls (484.97,68.91) and (486.55,69.45) .. (487.58,71.57) .. controls (488.62,73.68) and (490.2,74.22) .. (492.31,73.18) .. controls (494.43,72.14) and (496.01,72.68) .. (497.04,74.8) .. controls (498.07,76.92) and (499.65,77.46) .. (501.77,76.42) .. controls (503.89,75.38) and (505.47,75.92) .. (506.5,78.04) .. controls (507.53,80.16) and (509.11,80.7) .. (511.23,79.66) .. controls (513.35,78.62) and (514.93,79.16) .. (515.96,81.28) .. controls (516.99,83.4) and (518.57,83.94) .. (520.69,82.9) .. controls (522.8,81.86) and (524.38,82.4) .. (525.42,84.51) .. controls (526.45,86.63) and (528.03,87.17) .. (530.15,86.13) .. controls (532.27,85.09) and (533.85,85.63) .. (534.89,87.75) .. controls (535.92,89.87) and (537.5,90.41) .. (539.62,89.37) .. controls (541.74,88.33) and (543.32,88.87) .. (544.35,90.99) .. controls (545.38,93.11) and (546.96,93.65) .. (549.08,92.61) .. controls (551.19,91.57) and (552.77,92.11) .. (553.81,94.22) .. controls (554.84,96.34) and (556.42,96.88) .. (558.54,95.84) .. controls (560.66,94.8) and (562.24,95.34) .. (563.27,97.46) .. controls (564.3,99.58) and (565.88,100.12) .. (568,99.08) .. controls (570.12,98.04) and (571.7,98.58) .. (572.73,100.7) .. controls (573.76,102.82) and (575.34,103.36) .. (577.46,102.32) .. controls (579.57,101.28) and (581.15,101.82) .. (582.19,103.93) .. controls (583.22,106.05) and (584.8,106.59) .. (586.92,105.55) .. controls (589.04,104.51) and (590.62,105.05) .. (591.66,107.17) .. controls (592.69,109.29) and (594.27,109.83) .. (596.39,108.79) -- (597,109) -- (597,109) ;
		
		% Text Node
		\draw (234,129.4) node [anchor=north west][inner sep=0.75pt]    {$e_{a}$};
		% Text Node
		\draw (340,78.4) node [anchor=north west][inner sep=0.75pt]    {$e_{a}$};
		% Text Node
		\draw (113,83.4) node [anchor=north west][inner sep=0.75pt]    {$e_{a}$};
		% Text Node
		\draw (602,127.4) node [anchor=north west][inner sep=0.75pt]    {$e_{m}$};
		% Text Node
		\draw (702,76.4) node [anchor=north west][inner sep=0.75pt]    {$e_{a}$};
		% Text Node
		\draw (472,81.4) node [anchor=north west][inner sep=0.75pt]    {$e_{m}$};
		% Text Node
		\draw (522,22) node [anchor=north west][inner sep=0.75pt]    {$e_{a}$};	
	\end{tikzpicture}
\end{center}

With $e_a$ labeling $A$ and $e_m$ labeling $M$, these profiles encode operations
$A\times\cdots\times A\to A$ and $A\times\cdots\times A\times M\to M$.
The module input is last in the ordered input list, so it appears on
the left in our right-to-left drawings. The multiplication and left
action are represented by binary generators $\mu_a$ and $\mu_m$ of
profiles $(e_a,e_a;e_a)$ and $(e_a,e_m;e_m)$:

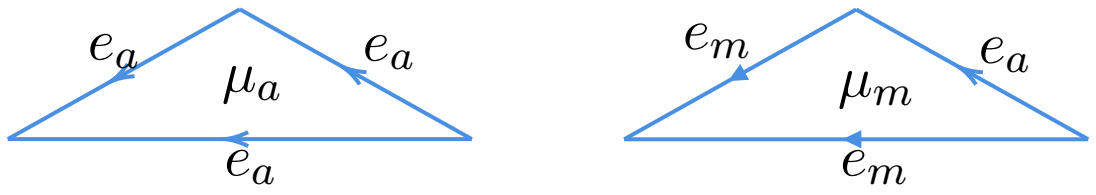
\begin{center}
\begin{tikzpicture}[x=0.5pt,y=0.5pt,yscale=-0.8,xscale=0.8]
			%uncomment if require: \path (0,116); %set diagram left start at 0, and has height of 116
			
			%Straight Lines [id:da2743514629985906] 
			\draw [color={rgb, 255:red, 74; green, 144; blue, 226 }  ,draw opacity=1 ][line width=0.75]    (274.84,21) -- (384,82.07) ;
			\draw [shift={(323.31,48.12)}, rotate = 29.23] [color={rgb, 255:red, 74; green, 144; blue, 226 }  ,draw opacity=1 ][line width=0.75]    (10.93,-3.29) .. controls (6.95,-1.4) and (3.31,-0.3) .. (0,0) .. controls (3.31,0.3) and (6.95,1.4) .. (10.93,3.29)   ;
			%Straight Lines [id:da2458406783735987] 
			\draw [color={rgb, 255:red, 74; green, 144; blue, 226 }  ,draw opacity=1 ][line width=0.75]    (165.67,82.07) -- (384,82.07) ;
			\draw [shift={(267.84,82.07)}, rotate = 0] [color={rgb, 255:red, 74; green, 144; blue, 226 }  ,draw opacity=1 ][line width=0.75]    (10.93,-3.29) .. controls (6.95,-1.4) and (3.31,-0.3) .. (0,0) .. controls (3.31,0.3) and (6.95,1.4) .. (10.93,3.29)   ;
			%Straight Lines [id:da6824309545734145] 
			\draw [color={rgb, 255:red, 74; green, 144; blue, 226 }  ,draw opacity=1 ][line width=0.75]    (165.67,82.07) -- (274.84,21) ;
			\draw [shift={(214.14,54.95)}, rotate = 330.77] [color={rgb, 255:red, 74; green, 144; blue, 226 }  ,draw opacity=1 ][line width=0.75]    (10.93,-3.29) .. controls (6.95,-1.4) and (3.31,-0.3) .. (0,0) .. controls (3.31,0.3) and (6.95,1.4) .. (10.93,3.29)   ;
			%Straight Lines [id:da9270322782894765] 
			\draw [color={rgb, 255:red, 74; green, 144; blue, 226 }  ,draw opacity=1 ][line width=0.75]    (564.84,21.08) -- (674,82.16) ;
			\draw [shift={(613.31,48.2)}, rotate = 29.23] [color={rgb, 255:red, 74; green, 144; blue, 226 }  ,draw opacity=1 ][line width=0.75]    (10.93,-3.29) .. controls (6.95,-1.4) and (3.31,-0.3) .. (0,0) .. controls (3.31,0.3) and (6.95,1.4) .. (10.93,3.29)   ;
			%Straight Lines [id:da7337750194283615] 
			\draw [color={rgb, 255:red, 74; green, 144; blue, 226 }  ,draw opacity=1 ][line width=0.75]    (455.67,82.16) -- (674,82.16) ;
			\draw [shift={(558.34,82.16)}, rotate = 0] [fill={rgb, 255:red, 74; green, 144; blue, 226 }  ,fill opacity=1 ][line width=0.08]  [draw opacity=0] (8.93,-4.29) -- (0,0) -- (8.93,4.29) -- cycle    ;
			%Straight Lines [id:da6793022159713502] 
			\draw [color={rgb, 255:red, 74; green, 144; blue, 226 }  ,draw opacity=1 ][line width=0.75]    (455.67,82.16) -- (564.84,21.08) ;
			\draw [shift={(504.58,54.79)}, rotate = 330.77] [fill={rgb, 255:red, 74; green, 144; blue, 226 }  ,fill opacity=1 ][line width=0.08]  [draw opacity=0] (8.93,-4.29) -- (0,0) -- (8.93,4.29) -- cycle    ;
			
			% Text Node
			\draw (200.87,30.46) node [anchor=north west][inner sep=0.75pt]    {$e_{a}$};
			% Text Node
			\draw (330.22,30.97) node [anchor=north west][inner sep=0.75pt]    {$e_{a}$};
			% Text Node
			\draw (265.13,84.98) node [anchor=north west][inner sep=0.75pt]    {$e_{a}$};
			% Text Node
			\draw (480.87,25.54) node [anchor=north west][inner sep=0.75pt]    {$e_{m}$};
			% Text Node
			\draw (620.04,32.05) node [anchor=north west][inner sep=0.75pt]    {$e_{a}$};
			% Text Node
			\draw (554.95,85.06) node [anchor=north west][inner sep=0.75pt]    {$e_{m}$};
			% Text Node
			\draw (263.95,45.54) node [anchor=north west][inner sep=0.75pt]    {$\mu _{a}$};
			% Text Node
			\draw (553.78,47.89) node [anchor=north west][inner sep=0.75pt]    {$\mu _{m}$};
		\end{tikzpicture}
\end{center}

The arrow shapes distinguish $e_a$ from $e_m$. Imposing associativity
and the left-module axiom gives the two relations below:

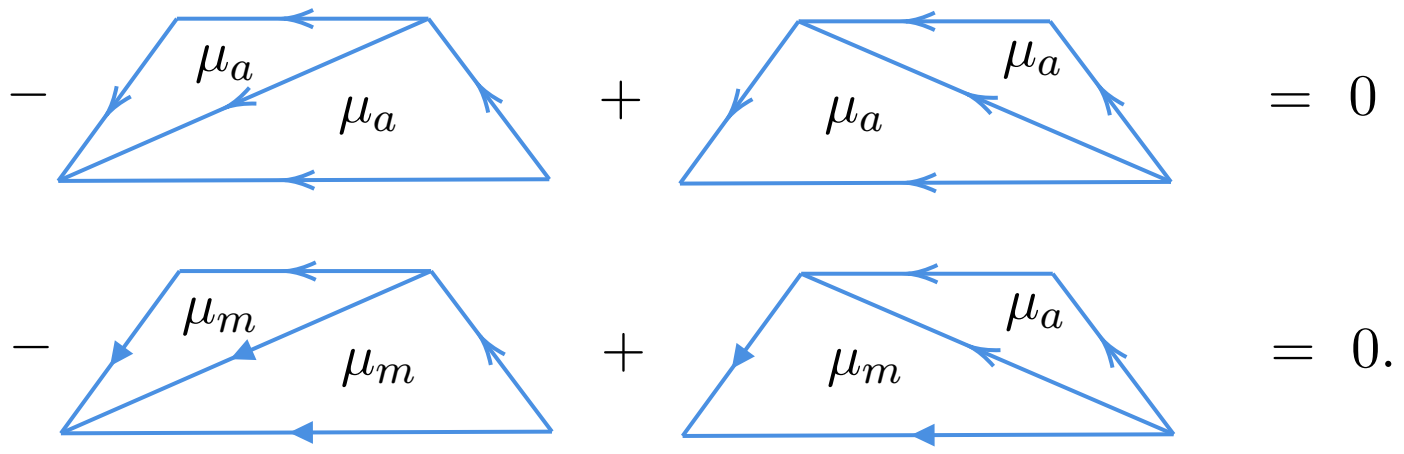
\begin{center}
\begin{tikzpicture}[x=0.5pt,y=0.5pt,yscale=-1,xscale=1]
			%uncomment if require: \path (0,207); %set diagram left start at 0, and has height of 207
			
			%Straight Lines [id:da6441544961134441] 
			\draw [color={rgb, 255:red, 74; green, 144; blue, 226 }  ,draw opacity=1 ][line width=0.75]    (207.67,87.07) -- (392.5,86.44) ;
			\draw [shift={(293.09,86.78)}, rotate = 359.8] [color={rgb, 255:red, 74; green, 144; blue, 226 }  ,draw opacity=1 ][line width=0.75]    (10.93,-3.29) .. controls (6.95,-1.4) and (3.31,-0.3) .. (0,0) .. controls (3.31,0.3) and (6.95,1.4) .. (10.93,3.29)   ;
			%Straight Lines [id:da757608279333391] 
			\draw [color={rgb, 255:red, 74; green, 144; blue, 226 }  ,draw opacity=1 ][line width=0.75]    (252.32,26) -- (347,26) ;
			\draw [shift={(292.66,26)}, rotate = 0] [color={rgb, 255:red, 74; green, 144; blue, 226 }  ,draw opacity=1 ][line width=0.75]    (10.93,-3.29) .. controls (6.95,-1.4) and (3.31,-0.3) .. (0,0) .. controls (3.31,0.3) and (6.95,1.4) .. (10.93,3.29)   ;
			%Straight Lines [id:da8978206509356398] 
			\draw [color={rgb, 255:red, 74; green, 144; blue, 226 }  ,draw opacity=1 ][line width=0.75]    (347,26) -- (392.5,86.44) ;
			\draw [shift={(365.54,50.63)}, rotate = 53.02] [color={rgb, 255:red, 74; green, 144; blue, 226 }  ,draw opacity=1 ][line width=0.75]    (10.93,-3.29) .. controls (6.95,-1.4) and (3.31,-0.3) .. (0,0) .. controls (3.31,0.3) and (6.95,1.4) .. (10.93,3.29)   ;
			%Straight Lines [id:da9179897864215004] 
			\draw [color={rgb, 255:red, 74; green, 144; blue, 226 }  ,draw opacity=1 ][line width=0.75]    (207.67,87.07) -- (252.32,26) ;
			\draw [shift={(225.86,62.19)}, rotate = 306.17] [color={rgb, 255:red, 74; green, 144; blue, 226 }  ,draw opacity=1 ][line width=0.75]    (10.93,-3.29) .. controls (6.95,-1.4) and (3.31,-0.3) .. (0,0) .. controls (3.31,0.3) and (6.95,1.4) .. (10.93,3.29)   ;
			%Straight Lines [id:da21763065939913195] 
			\draw [color={rgb, 255:red, 74; green, 144; blue, 226 }  ,draw opacity=1 ][line width=0.75]    (207.67,87.07) -- (347,26) ;
			\draw [shift={(270.92,59.35)}, rotate = 336.33] [color={rgb, 255:red, 74; green, 144; blue, 226 }  ,draw opacity=1 ][line width=0.75]    (10.93,-3.29) .. controls (6.95,-1.4) and (3.31,-0.3) .. (0,0) .. controls (3.31,0.3) and (6.95,1.4) .. (10.93,3.29)   ;
			%Straight Lines [id:da5058307364936606] 
			\draw [color={rgb, 255:red, 74; green, 144; blue, 226 }  ,draw opacity=1 ][line width=0.75]    (441.67,88.07) -- (626.5,87.44) ;
			\draw [shift={(527.09,87.78)}, rotate = 359.8] [color={rgb, 255:red, 74; green, 144; blue, 226 }  ,draw opacity=1 ][line width=0.75]    (10.93,-3.29) .. controls (6.95,-1.4) and (3.31,-0.3) .. (0,0) .. controls (3.31,0.3) and (6.95,1.4) .. (10.93,3.29)   ;
			%Straight Lines [id:da8703328585965336] 
			\draw [color={rgb, 255:red, 74; green, 144; blue, 226 }  ,draw opacity=1 ][line width=0.75]    (486.32,27) -- (581,27) ;
			\draw [shift={(526.66,27)}, rotate = 0] [color={rgb, 255:red, 74; green, 144; blue, 226 }  ,draw opacity=1 ][line width=0.75]    (10.93,-3.29) .. controls (6.95,-1.4) and (3.31,-0.3) .. (0,0) .. controls (3.31,0.3) and (6.95,1.4) .. (10.93,3.29)   ;
			%Straight Lines [id:da09263160229897238] 
			\draw [color={rgb, 255:red, 74; green, 144; blue, 226 }  ,draw opacity=1 ][line width=0.75]    (581,27) -- (626.5,87.44) ;
			\draw [shift={(599.54,51.63)}, rotate = 53.02] [color={rgb, 255:red, 74; green, 144; blue, 226 }  ,draw opacity=1 ][line width=0.75]    (10.93,-3.29) .. controls (6.95,-1.4) and (3.31,-0.3) .. (0,0) .. controls (3.31,0.3) and (6.95,1.4) .. (10.93,3.29)   ;
			%Straight Lines [id:da30590672430193944] 
			\draw [color={rgb, 255:red, 74; green, 144; blue, 226 }  ,draw opacity=1 ][line width=0.75]    (441.67,88.07) -- (486.32,27) ;
			\draw [shift={(459.86,63.19)}, rotate = 306.17] [color={rgb, 255:red, 74; green, 144; blue, 226 }  ,draw opacity=1 ][line width=0.75]    (10.93,-3.29) .. controls (6.95,-1.4) and (3.31,-0.3) .. (0,0) .. controls (3.31,0.3) and (6.95,1.4) .. (10.93,3.29)   ;
			%Straight Lines [id:da28169632237429676] 
			\draw [color={rgb, 255:red, 74; green, 144; blue, 226 }  ,draw opacity=1 ][line width=0.75]    (486.32,27) -- (626.5,87.44) ;
			\draw [shift={(549.98,54.45)}, rotate = 23.32] [color={rgb, 255:red, 74; green, 144; blue, 226 }  ,draw opacity=1 ][line width=0.75]    (10.93,-3.29) .. controls (6.95,-1.4) and (3.31,-0.3) .. (0,0) .. controls (3.31,0.3) and (6.95,1.4) .. (10.93,3.29)   ;
			%Straight Lines [id:da08834603040539857] 
			\draw [color={rgb, 255:red, 74; green, 144; blue, 226 }  ,draw opacity=1 ][line width=0.75]    (208.67,182.07) -- (393.5,181.44) ;
			\draw [shift={(294.59,181.78)}, rotate = 359.8] [fill={rgb, 255:red, 74; green, 144; blue, 226 }  ,fill opacity=1 ][line width=0.08]  [draw opacity=0] (8.93,-4.29) -- (0,0) -- (8.93,4.29) -- cycle    ;
			%Straight Lines [id:da014231528947963312] 
			\draw [color={rgb, 255:red, 74; green, 144; blue, 226 }  ,draw opacity=1 ][line width=0.75]    (253.32,121) -- (348,121) ;
			\draw [shift={(293.66,121)}, rotate = 0] [color={rgb, 255:red, 74; green, 144; blue, 226 }  ,draw opacity=1 ][line width=0.75]    (10.93,-3.29) .. controls (6.95,-1.4) and (3.31,-0.3) .. (0,0) .. controls (3.31,0.3) and (6.95,1.4) .. (10.93,3.29)   ;
			%Straight Lines [id:da6523362679556463] 
			\draw [color={rgb, 255:red, 74; green, 144; blue, 226 }  ,draw opacity=1 ][line width=0.75]    (348,121) -- (393.5,181.44) ;
			\draw [shift={(366.54,145.63)}, rotate = 53.02] [color={rgb, 255:red, 74; green, 144; blue, 226 }  ,draw opacity=1 ][line width=0.75]    (10.93,-3.29) .. controls (6.95,-1.4) and (3.31,-0.3) .. (0,0) .. controls (3.31,0.3) and (6.95,1.4) .. (10.93,3.29)   ;
			%Straight Lines [id:da3583571033623294] 
			\draw [color={rgb, 255:red, 74; green, 144; blue, 226 }  ,draw opacity=1 ][line width=0.75]    (208.67,182.07) -- (253.32,121) ;
			\draw [shift={(227.16,156.78)}, rotate = 306.17] [fill={rgb, 255:red, 74; green, 144; blue, 226 }  ,fill opacity=1 ][line width=0.08]  [draw opacity=0] (8.93,-4.29) -- (0,0) -- (8.93,4.29) -- cycle    ;
			%Straight Lines [id:da4345125206192829] 
			\draw [color={rgb, 255:red, 74; green, 144; blue, 226 }  ,draw opacity=1 ][line width=0.75]    (208.67,182.07) -- (348,121) ;
			\draw [shift={(272.38,154.15)}, rotate = 336.33] [fill={rgb, 255:red, 74; green, 144; blue, 226 }  ,fill opacity=1 ][line width=0.08]  [draw opacity=0] (8.93,-4.29) -- (0,0) -- (8.93,4.29) -- cycle    ;
			%Straight Lines [id:da834179610120251] 
			\draw [color={rgb, 255:red, 74; green, 144; blue, 226 }  ,draw opacity=1 ][line width=0.75]    (442.67,183.07) -- (627.5,182.44) ;
			\draw [shift={(528.59,182.78)}, rotate = 359.8] [fill={rgb, 255:red, 74; green, 144; blue, 226 }  ,fill opacity=1 ][line width=0.08]  [draw opacity=0] (8.93,-4.29) -- (0,0) -- (8.93,4.29) -- cycle    ;
			%Straight Lines [id:da3543124047806566] 
			\draw [color={rgb, 255:red, 74; green, 144; blue, 226 }  ,draw opacity=1 ][line width=0.75]    (487.32,122) -- (582,122) ;
			\draw [shift={(527.66,122)}, rotate = 0] [color={rgb, 255:red, 74; green, 144; blue, 226 }  ,draw opacity=1 ][line width=0.75]    (10.93,-3.29) .. controls (6.95,-1.4) and (3.31,-0.3) .. (0,0) .. controls (3.31,0.3) and (6.95,1.4) .. (10.93,3.29)   ;
			%Straight Lines [id:da2672971964339226] 
			\draw [color={rgb, 255:red, 74; green, 144; blue, 226 }  ,draw opacity=1 ][line width=0.75]    (582,122) -- (627.5,182.44) ;
			\draw [shift={(600.54,146.63)}, rotate = 53.02] [color={rgb, 255:red, 74; green, 144; blue, 226 }  ,draw opacity=1 ][line width=0.75]    (10.93,-3.29) .. controls (6.95,-1.4) and (3.31,-0.3) .. (0,0) .. controls (3.31,0.3) and (6.95,1.4) .. (10.93,3.29)   ;
			%Straight Lines [id:da04931642769782807] 
			\draw [color={rgb, 255:red, 74; green, 144; blue, 226 }  ,draw opacity=1 ][line width=0.75]    (442.67,183.07) -- (487.32,122) ;
			\draw [shift={(461.16,157.78)}, rotate = 306.17] [fill={rgb, 255:red, 74; green, 144; blue, 226 }  ,fill opacity=1 ][line width=0.08]  [draw opacity=0] (8.93,-4.29) -- (0,0) -- (8.93,4.29) -- cycle    ;
			%Straight Lines [id:da189259921127787] 
			\draw [color={rgb, 255:red, 74; green, 144; blue, 226 }  ,draw opacity=1 ][line width=0.75]    (487.32,122) -- (627.5,182.44) ;
			\draw [shift={(550.98,149.45)}, rotate = 23.32] [color={rgb, 255:red, 74; green, 144; blue, 226 }  ,draw opacity=1 ][line width=0.75]    (10.93,-3.29) .. controls (6.95,-1.4) and (3.31,-0.3) .. (0,0) .. controls (3.31,0.3) and (6.95,1.4) .. (10.93,3.29)   ;
			
			% Text Node
			\draw (257.13,34.3) node [anchor=north west][inner sep=0.75pt]    {$\mu _{a}$};
			% Text Node
			\draw (311.04,53.63) node [anchor=north west][inner sep=0.75pt]    {$\mu _{a}$};
			% Text Node
			\draw (494.04,53.63) node [anchor=north west][inner sep=0.75pt]    {$\mu _{a}$};
			% Text Node
			\draw (561.04,32.63) node [anchor=north west][inner sep=0.75pt]    {$\mu _{a}$};
			% Text Node
			\draw (186,45.4) node [anchor=north west][inner sep=0.75pt]    {$-$};
			% Text Node
			\draw (409,47.4) node [anchor=north west][inner sep=0.75pt]    {$+$};
			% Text Node
			\draw (661,46.4) node [anchor=north west][inner sep=0.75pt]    {$=\ 0$};
			% Text Node
			\draw (252.13,129.3) node [anchor=north west][inner sep=0.75pt]    {$\mu _{m}$};
			% Text Node
			\draw (312.04,148.63) node [anchor=north west][inner sep=0.75pt]    {$\mu _{m}$};
			% Text Node
			\draw (495.04,148.63) node [anchor=north west][inner sep=0.75pt]    {$\mu _{m}$};
			% Text Node
			\draw (562.04,127.63) node [anchor=north west][inner sep=0.75pt]    {$\mu _{a}$};
			% Text Node
			\draw (187,140.4) node [anchor=north west][inner sep=0.75pt]    {$-$};
			% Text Node
			\draw (410,142.4) node [anchor=north west][inner sep=0.75pt]    {$+$};
			% Text Node
			\draw (662,141.4) node [anchor=north west][inner sep=0.75pt]    {$=\ 0$.};

		\end{tikzpicture}
\end{center}

These generators and relations define the quadratic $Q$-operad
$\mathscr L\mathsf{Mod}$, whose algebras are precisely such pairs $(A,M)$.
We treat this example in detail in \S~\ref{ss_example_left_module};
related pictures appear in \cite[\S5.1]{hinich2015rectification}.
We now develop the general definitions, beginning with quivers and
their admissible profiles.

\subsection{Quivers and collections}
\subsubsection{Quivers}
\label{sss_directed_graphs}
A \textit{quiver} is a tuple $Q=(Q_0,Q,s,t)$ of vertex and arrow sets and maps
$s,t:Q\to Q_0$. We use $Q$ for both the quiver and its arrow set.
A \textit{morphism} $(\bar\phi,\phi):Q\to Q'$ consists of maps on vertices and
arrows satisfying $s'\phi=\bar\phi s$ and $t'\phi=\bar\phi t$.
A \textit{subquiver} has vertex and arrow subsets $Q'_0,Q'$ with
$s(Q'),t(Q')\subseteq Q'_0$.
Write $Q^*$ for the quiver of finite composable paths, with the same
vertices as $Q$:
\begin{equation*} Q^*=Q_0\sqcup\bigsqcup_{n\geq1}
 \{(e_1,\ldots,e_n)\in Q^n:t(e_i)=s(e_{i+1})\}.
\end{equation*}
The length-zero path $\varnothing_v$ at $v\in Q_0$ has source and target $v$;
for a nonempty path, $s^*(\vec e)=s(e_1)$ and
$t^*(\vec e)=t(e_n)$. For nonempty paths $\vec e=(e_1,\ldots,e_m)$ and
$\vec f=(f_1,\ldots,f_n)$ satisfying $t^*(\vec e)=s^*(\vec f)$, their
concatenation is denoted by
$\vec e*\vec f=(e_1,\ldots,e_m,f_1,\ldots,f_n)$.
Length-zero paths act as units for concatenation.

A \textit{profile-loop} is a pair $(\vec e;e_0)$ consisting of a
path $\vec e=(e_1,\ldots,e_n)\in Q^*$ and an arrow $e_0$
with $s^*(\vec e)=s(e_0)$ and $t^*(\vec e)=t(e_0)$.
Its \textit{arity} is $n$. The profile $(e;e)$ is the \textit{identity profile} at $e$.
The set of profile-loops is denoted by $\mathcal P_Q$.
For the \textit{terminal quiver}, with one vertex and one loop,
$\mathcal P_Q$ is the set of nonnegative integers.

\subsubsection{Collections over quivers}
\label{sss_collections_over_E}
We work in graded vector spaces or cochain complexes over the fixed
field $\Bbbk$, with the Koszul tensor symmetry. A \textit{$Q$-collection}
is a family $S=\{S(\vec e;e_0)\}_{\mathcal P_Q}$; a \textit{dg collection} has
a differential $\delta_S$ of degree $1$ on each component. \textit{Maps of graded collections} have degree zero, and \textit{maps of dg collections} are componentwise
cochain maps. \textit{Subcollections} are defined componentwise.

The \textit{internal Hom collection} is
\[
 \operatorname{Coll}_Q(S,T)(\vec e;e_0)
 =\underline{\operatorname{Hom}}_{\Bbbk}
       (S(\vec e;e_0),T(\vec e;e_0)),
\]
where, for complexes $A,B$,
\[
 \underline{\operatorname{Hom}}_{\Bbbk}(A,B)^r
 =\prod_j\operatorname{Hom}_{\Bbbk}(A^j,B^{j+r}),
 \qquad \partial f=\delta_B f-(-1)^{|f|}f\delta_A.
\]
Its degree-zero cycles are the cochain maps.

We also use \textit{suspension} and \textit{desuspension} componentwise on collections.
For a graded vector space $W$, put $(sW)^i=W^{i-1}$, so
$|sx|=|x|+1$. For a complex set
$\delta_{sW}(sx)=-s\delta_Wx$; similarly
$|s^{-1}x|=|x|-1$ and
$\delta_{s^{-1}W}(s^{-1}x)=-s^{-1}\delta_Wx$.

Over the terminal quiver, a $Q$-collection is a sequence
$\{S(n)\}_{n\geq0}$, the usual nonsymmetric collection.

\subsubsection{Substitution product and unit}
\label{sss_substitution_product}
\label{sss_unit_collection}
For $Q$-collections $S,T$, define their \textit{substitution product} by
\begin{equation}\label{eq_S_1_circ_S_2}
 (S\circ_Q T)(\vec e;e_0)=
 \bigoplus_{\substack{k\geq0,\ \vec e=\vec g_1*\cdots*\vec g_k\\
 (\vec f;e_0),(\vec g_i;f_i)\in\mathcal P_Q}}
 S(\vec f;e_0)\otimes\bigotimes_{i=1}^k T(\vec g_i;f_i),
 \qquad \vec f=(f_1,\ldots,f_k).
\end{equation}
Empty tensor products are $\Bbbk$. On dg collections the differential
is the tensor-product differential, and maps act by tensor products
on the displayed summands. Reindexing three-level substitutions and
using the Koszul symmetry gives the associator
$(S\circ_Q T)\circ_Q U\cong S\circ_Q(T\circ_Q U)$.
Both sides index the same three-level colored trees; reindexing four
levels proves the pentagon identity.

The \textit{unit collection} $\mathbb I_Q$ is concentrated in degree zero, with
$\mathbb I_Q(e;e)=\Bbbk$ for every arrow $e$ and all other components zero.
Its distinguished basis element at $(e;e)$ is $1_e$.
The only nonzero summands in $\mathbb I_Q\circ_Q S$ and
$S\circ_Q\mathbb I_Q$ give the unit isomorphisms
$\lambda_S:\mathbb I_Q\circ_Q S\cong S$ and
$\rho_S:S\circ_Q\mathbb I_Q\cong S$.

\subsection{Quiver-operads and quiver-cooperads}
\label{ss_operad_over_graphs}
\label{s_fc_comulticategories}
\label{ss_augmentation_coaugmentation}
A \textit{vertically discrete $\mathbf{fc}$-multicategory} has objects $Q_0$,
horizontal arrows $Q$, and only identity vertical arrows. Its $2$-cells
therefore have boundaries
\begin{equation*} \xymatrix@C=18pt{
 v_n&\cdots\ar[l]_{e_n}&v_1\ar[l]&
 v_0\ar[l]_{e_1}\ar@/^18pt/[lll]^{e_0}}
 \qquad s(e_0)=v_0,\quad t(e_0)=v_n.
\end{equation*}
Following \cite[\S3.2]{Yuan_higher_operad_Fuk}, we draw $v_0$ on the right
and $v_n$ on the left, so the inputs $e_1,\ldots,e_n$ run from right to
left. Profiles and tensor products retain this indexed order, and
$\circ_i$ denotes insertion at the $e_i$-slot.
With multilinear compositions, this gives the following equivalent
monoidal formulation; see \cite{leinster1999fc,leinster2004higher,
cruttwell2009unified} for the unenriched framework.

\begin{defn}\label{defn_operad_over_E}
A \textit{$Q$-operad} is a monoid $(M,\gamma,\iota)$ in graded
$Q$-collections under $\circ_Q$: its multiplication
$\gamma:M\circ_Q M\to M$ is associative and has unit
$\iota:\mathbb I_Q\to M$. We write
$\mathscr M=(Q_0,Q,M,\iota,\gamma)$ and
$\id_e=\iota(1_e)$.
\end{defn}
The multiplication has components
$ M(\vec e;e_0)\otimes\bigotimes_iM(\vec g_i;e_i)
 \to M(\vec g_1*\cdots*\vec g_n;e_0)$.
For $1\leq i\leq n$, inserting identities in all other slots gives \textit{partial compositions}
\[
 \circ_i:M(e_1,\ldots,e_n;e_0)\otimes M(\vec g;e_i)
 \longrightarrow M(e_1,\ldots,e_{i-1},\vec g,e_{i+1},\ldots,e_n;e_0).
\]
For homogeneous elements, their axioms are
\[
 (u\circ_i v)\circ_{i+j-1}w=u\circ_i(v\circ_jw),\qquad
 (u\circ_i v)\circ_{j+m-1}w
 =(-1)^{|v||w|}(u\circ_jw)\circ_i v\quad(i<j),
\]
where $m$ is the arity of $v$, together with
$\id_{e_0}\circ_1u=u=u\circ_i\id_{e_i}$.
These identities recover the monoid axioms by iterated insertion.
For $\vec h=\vec e\circ_i\vec g$, the notation on the right means
replacement of the $i$-th arrow of $\vec e$ by $\vec g$.
The sum of the partial compositions is
\begin{equation}\label{eq_gamma_(1)}
 \gamma_{(1)}:
 \bigoplus_{\substack{\vec h=\vec e\circ_i\vec g\\
 (\vec e;e_0),(\vec g;e_i)\in\mathcal P_Q}}
 M(\vec e;e_0)\otimes M(\vec g;e_i)
 \longrightarrow M(\vec h;e_0).
\end{equation}

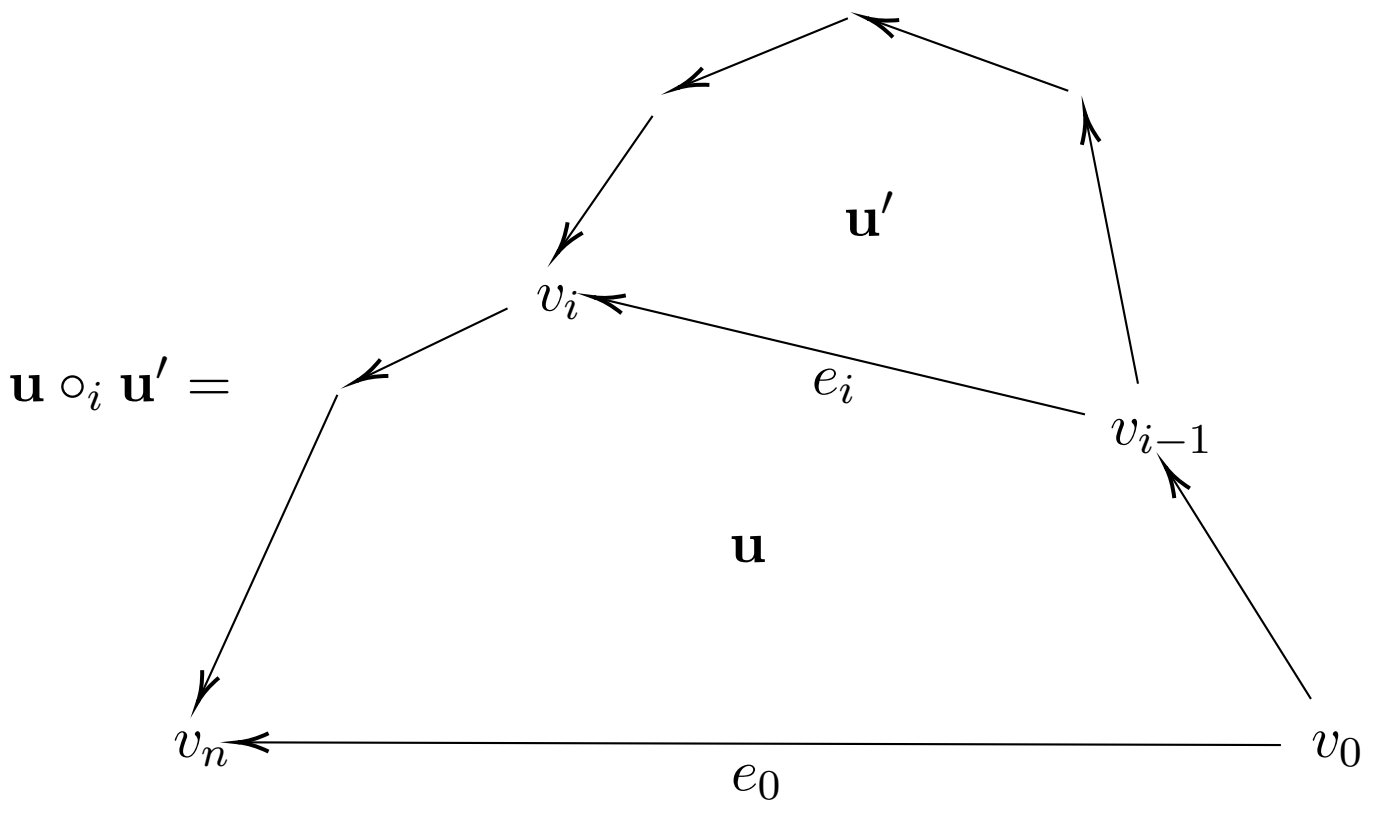
\begin{figure}[ht]
		\centering
		\begin{tikzpicture}[x=0.75pt,y=0.75pt,yscale=-0.7,xscale=0.7]
			%uncomment if require: \path (0,375); %set diagram left start at 0, and has height of 375
			
			%Straight Lines [id:da8827101392474932] 
			\draw    (508.23,306.53) -- (136.4,305.57) ;
			\draw [shift={(134.4,305.56)}, rotate = 0.15] [color={rgb, 255:red, 0; green, 0; blue, 0 }  ][line width=0.75]    (10.93,-3.29) .. controls (6.95,-1.4) and (3.31,-0.3) .. (0,0) .. controls (3.31,0.3) and (6.95,1.4) .. (10.93,3.29)   ;
			%Straight Lines [id:da40636261009415153] 
			\draw    (438,188) -- (263.95,146.46) ;
			\draw [shift={(262,146)}, rotate = 13.42] [color={rgb, 255:red, 0; green, 0; blue, 0 }  ][line width=0.75]    (10.93,-3.29) .. controls (6.95,-1.4) and (3.31,-0.3) .. (0,0) .. controls (3.31,0.3) and (6.95,1.4) .. (10.93,3.29)   ;
			%Straight Lines [id:da6001351195599821] 
			\draw    (283,81) -- (250.14,128.36) ;
			\draw [shift={(249,130)}, rotate = 304.76] [color={rgb, 255:red, 0; green, 0; blue, 0 }  ][line width=0.75]    (10.93,-3.29) .. controls (6.95,-1.4) and (3.31,-0.3) .. (0,0) .. controls (3.31,0.3) and (6.95,1.4) .. (10.93,3.29)   ;
			%Straight Lines [id:da21538714651971236] 
			\draw    (353,46) -- (293.85,70.24) ;
			\draw [shift={(292,71)}, rotate = 337.71] [color={rgb, 255:red, 0; green, 0; blue, 0 }  ][line width=0.75]    (10.93,-3.29) .. controls (6.95,-1.4) and (3.31,-0.3) .. (0,0) .. controls (3.31,0.3) and (6.95,1.4) .. (10.93,3.29)   ;
			%Straight Lines [id:da06384884716881045] 
			\draw    (432,72) -- (361.88,46.68) ;
			\draw [shift={(360,46)}, rotate = 19.86] [color={rgb, 255:red, 0; green, 0; blue, 0 }  ][line width=0.75]    (10.93,-3.29) .. controls (6.95,-1.4) and (3.31,-0.3) .. (0,0) .. controls (3.31,0.3) and (6.95,1.4) .. (10.93,3.29)   ;
			%Straight Lines [id:da33828689530558] 
			\draw    (457,177) -- (438.38,81.96) ;
			\draw [shift={(438,80)}, rotate = 78.92] [color={rgb, 255:red, 0; green, 0; blue, 0 }  ][line width=0.75]    (10.93,-3.29) .. controls (6.95,-1.4) and (3.31,-0.3) .. (0,0) .. controls (3.31,0.3) and (6.95,1.4) .. (10.93,3.29)   ;
			%Straight Lines [id:da6099773766995279] 
			\draw    (231,150) -- (178.8,175.13) ;
			\draw [shift={(177,176)}, rotate = 334.29] [color={rgb, 255:red, 0; green, 0; blue, 0 }  ][line width=0.75]    (10.93,-3.29) .. controls (6.95,-1.4) and (3.31,-0.3) .. (0,0) .. controls (3.31,0.3) and (6.95,1.4) .. (10.93,3.29)   ;
			%Straight Lines [id:da2220578526215704] 
			\draw    (170,181) -- (120.83,289.18) ;
			\draw [shift={(120,291)}, rotate = 294.44] [color={rgb, 255:red, 0; green, 0; blue, 0 }  ][line width=0.75]    (10.93,-3.29) .. controls (6.95,-1.4) and (3.31,-0.3) .. (0,0) .. controls (3.31,0.3) and (6.95,1.4) .. (10.93,3.29)   ;
			%Straight Lines [id:da8870208807355183] 
			\draw    (519,290) -- (468.06,208.69) ;
			\draw [shift={(467,207)}, rotate = 57.93] [color={rgb, 255:red, 0; green, 0; blue, 0 }  ][line width=0.75]    (10.93,-3.29) .. controls (6.95,-1.4) and (3.31,-0.3) .. (0,0) .. controls (3.31,0.3) and (6.95,1.4) .. (10.93,3.29)   ;
			
			% Text Node
			\draw (109.27,300) node [anchor=north west][inner sep=0.75pt]    {$v_{n}$};
			% Text Node
			\draw (517.22,300) node [anchor=north west][inner sep=0.75pt]    {$v_{0}$};
			% Text Node
			\draw (239.29,140) node [anchor=north west][inner sep=0.75pt]    {$v_{i}$};
			% Text Node
			\draw (445.11,188) node [anchor=north west][inner sep=0.75pt]    {$v_{i-1}$};
			% Text Node
			\draw (309.26,311.98) node [anchor=north west][inner sep=0.75pt]    {$e_{0}$};
			% Text Node
			\draw (308.6,230.07) node [anchor=north west][inner sep=0.75pt]    {$\mathbf{u}$};
			% Text Node
			\draw (349.6,106.07) node [anchor=north west][inner sep=0.75pt]    {$\mathbf{u} '$};
			% Text Node
			\draw (338.26,169.98) node [anchor=north west][inner sep=0.75pt]    {$e_{i}$};
			% Text Node
			\draw (50,165.4) node [anchor=north west][inner sep=0.75pt]    {$\mathbf{u} \circ_i \mathbf u' =$};
		\end{tikzpicture}
		\caption{The partial composition of 2-cells $\mathbf u$ and $\mathbf u'$}
		\label{fig:2_cell_composition}
	\end{figure}

%For a one-vertex quiver, all arrow sequences are composable, so
%$Q$-operads are positive-arity nonsymmetric colored operads with color
%set $Q$. The terminal quiver recovers ordinary positive-arity
%nonsymmetric operads.

In the dg setting, a \textit{dg $Q$-operad} is a monoid in dg $Q$-collections.
Equivalently, its units are closed and its differential satisfies
$ \delta(u\circ_i v)=(\delta u)\circ_i v
                  +(-1)^{|u|}u\circ_i\delta v$.
A \textit{morphism of (dg) $Q$-operads} is a map of (dg) collections preserving
units and all partial compositions. It fixes the underlying quiver;
morphisms over a quiver map are introduced in
\S~\ref{ss_change_quiver}.

\medskip
\noindent\textit{Quiver-cooperads.}
\begin{defn}\label{defn_E_cooperad}
A \textit{$Q$-cooperad} is a comonoid $(C,\Delta,\varepsilon)$ in
$Q$-collections. Thus
\begin{equation*} \Delta:C(\vec e;e_0)\longrightarrow
 \bigoplus_{\substack{\vec e=\vec g_1*\cdots*\vec g_k\\
 (\vec f;e_0),(\vec g_i;f_i)\in\mathcal P_Q}}
 C(\vec f;e_0)\otimes\bigotimes_i C(\vec g_i;f_i)
\end{equation*}
is coassociative, and $\varepsilon:C\to\mathbb I_Q$ satisfies
$(\varepsilon\circ_Q\id_C)\Delta=\lambda_C^{-1}$ and
$(\id_C\circ_Q\varepsilon)\Delta=\rho_C^{-1}$.
We write $\mathscr C=(Q_0,Q,C,\Delta,\varepsilon)$.
A \textit{dg $Q$-cooperad} is a comonoid in dg $Q$-collections;
in particular, both $\Delta$ and $\varepsilon$ are cochain maps.
\end{defn}

The unit collection is both an operad and a cooperad, with
$1_e\circ_1 1_e=1_e$, $\Delta(1_e)=1_e\otimes1_e$ and
$\varepsilon(1_e)=1$.
Projecting a component of $\Delta$ and applying the counit to every
inner factor except the $i$-th gives
$\Delta_i:C(\vec e\circ_i\vec g;e_0)
\to C(\vec e;e_0)\otimes C(\vec g;e_i)$.
Their sum is the \textit{infinitesimal coproduct}
\begin{equation}\label{eq_Delta_(1)}
 \Delta_{(1)}:C(\vec h;e_0)\longrightarrow
 \bigoplus_{\substack{\vec h=\vec e\circ_i\vec g\\
 (\vec e;e_0),(\vec g;e_i)\in\mathcal P_Q}}
 C(\vec e;e_0)\otimes C(\vec g;e_i)
 = (C\circ_{(1)}C)(\vec h;e_0).
\end{equation}
More generally, $S\circ_{(1)}T$ denotes the same direct sum with
factors $S(\vec e;e_0)\otimes T(\vec g;e_i)$.
The slot is part of the summand index, including when an inner factor
is an identity.

A \textit{morphism} $F:\mathscr C\to\mathscr D$ \textit{of (dg) $Q$-cooperads} is a
map of (dg) collections satisfying
$\Delta_DF=(F\circ_QF)\Delta_C$ and $\varepsilon_DF=\varepsilon_C$.

\begin{prop}\label{prop_extension_by_zero}
Let the arrows of $Q$ be the colors of a nonsymmetric
colored collection. Extension by zero on profiles outside $\mathcal P_Q$
is a fully faithful strong monoidal embedding of $Q$-collections into
these colored collections. Its image consists of collections supported
on $\mathcal P_Q$. It therefore identifies $Q$-operads and $Q$-cooperads
with the supported colored operads and cooperads.
\end{prop}
\begin{proof}
Admissible profiles are closed under substitution, so the colored
substitution product of extended collections has exactly the summands
of \eqref{eq_S_1_circ_S_2}. The colored unit extends $\mathbb I_Q$,
and these identifications respect associativity and units.
Restriction to admissible profiles is inverse to extension by zero on
supported collections. The operad and cooperad assertions follow by
passing to monoids and comonoids.
\end{proof}

\medskip
\noindent\textit{Augmentations and coaugmentations.}
A \textit{coaugmentation} is a dg cooperad morphism $\iota_C:\mathbb I_Q\to\mathscr C$ splitting the counit. An \textit{augmentation} is a dg operad morphism $\epsilon_M:\mathscr M\to\mathbb I_Q$ splitting the unit.
For an augmented dg operad $\mathscr M$ and a coaugmented dg cooperad $\mathscr C$, the \textit{reduced parts} are
$\overline M=\ker\epsilon_M$ and $\overline C=\ker\varepsilon_C\cong\operatorname{coker}\iota_C$.
The (co)augmentations give canonical dg splittings
$M=\mathbb I_Q\oplus\overline M$ and $C=\mathbb I_Q\oplus\overline C$, with $\overline M$ a dg operadic ideal.

Let $p:C\to\overline C$ be the projection. The \textit{reduced infinitesimal decomposition} is
$\overline\Delta_{(1)}=(p\circ_{(1)}p)\Delta_{(1)}|_{\overline C}$.

\subsection{Algebras over quiver-operads}
\label{ss_algebras}
\label{ss_endomrophism_E_operad}
Let $X=\{(X(e),\mathsf d_e^X)\}_{e\in Q}$ be a family of cochain
complexes indexed by the arrows of $Q$. We call such a family a
\textit{$Q$-indexed family}. For two such families $X,Y$, the vector space
of families of cochain maps is denoted by
\begin{equation*} \Hom_Q(X,Y)=Z^0\underline{\operatorname{Hom}}_Q(X,Y),
 \qquad
 \underline{\operatorname{Hom}}_Q(X,Y)
 =\prod_{e\in Q}\underline{\operatorname{Hom}}_{\Bbbk}(X(e),Y(e)).
\end{equation*}
We distinguish these morphisms from arbitrary homogeneous elements of
the internal Hom complex.

For a path $\vec e=(e_1,\ldots,e_n)$ set
$X^*(\vec e)=X(e_1)\otimes\cdots\otimes X(e_n)$, with differential
\begin{equation*} \mathsf d^X_{\vec e}(x_1\otimes\cdots\otimes x_n)
 =\sum_i(-1)^{\sum_{j<i}|x_j|}
 x_1\otimes\cdots\otimes\mathsf d^X_{e_i}x_i\otimes\cdots\otimes x_n.
\end{equation*}
Define
$\mathbf{End}_Q(X)(\vec e;e_0)=
\underline{\operatorname{Hom}}_{\Bbbk}(X^*(\vec e),X(e_0))$.
Its differential is
\begin{equation}\label{widehat_d_eq}
 \mathsf d\xi=\mathsf d^X_{e_0}\xi-(-1)^{|\xi|}\xi\mathsf d^X_{\vec e}.
\end{equation}
For homogeneous $f,g$, with $g$ of arity $m$ and output $e_i$, put
\begin{equation}\label{eq_End_X_hat_composition}
\begin{aligned}
 &(f\circ_i g)(x_1,\ldots,x_{i-1},y_1,\ldots,y_m,x_{i+1},\ldots,x_n)\\
 &\qquad=(-1)^{|g|\sum_{j<i}|x_j|}
 f(x_1,\ldots,x_{i-1},g(y_1,\ldots,y_m),x_{i+1},\ldots,x_n).
\end{aligned}
\end{equation}
This operation has degree $|f|+|g|$ and defines a cochain map on the
tensor product of the two Hom complexes. Associativity of composition
and the tensor sign rule give the dg operad axioms, with units
$\id_{X(e)}$. This is the \textit{endomorphism $Q$-operad}; compare
\cite[Proposition 5.10]{Yuan_higher_operad_Fuk}.

\medskip

An \textit{$\mathscr M$-algebra} is a $Q$-indexed family $X$ together
with a morphism of dg $Q$-operads
$\alpha:\mathscr M\to\mathbf{End}_Q(X)$.
For $o\in M(\vec e;e_0)$, denote the corresponding operation on $X$ by
\begin{equation}\label{eq_algebra_operation_notation}
 o_X:=\alpha_{\vec e;e_0}(o).
\end{equation}
The component maps
$
 \alpha_{\vec e;e_0}:M(\vec e;e_0)
 \to\underline{\operatorname{Hom}}_{\Bbbk}
       (X^*(\vec e),X(e_0))
$
preserve units and satisfy
\begin{equation}\label{eq_alpha_compatibility}
 \mathsf d\alpha(u)=\alpha(\delta_Mu),
 \qquad
 \alpha(u\circ_i v)=\alpha(u)\circ_i\alpha(v).
\end{equation}
Equivalently, consider the endofunctor
\begin{equation*} \mathbb T_{\mathscr M}(X)(e_0)
 =\bigoplus_{(\vec e;e_0)\in\mathcal P_Q}
     M(\vec e;e_0)\otimes X^*(\vec e),
\end{equation*}
with tensor differential and maps acting on every $X$-factor.
Operadic units and composition induce a monad structure, with
$\eta_X(x)=\id_e\otimes x$ and multiplication $\mu_X$ given by
substitution with the Koszul signs for regrouping tensor factors.
Tensor--Hom adjunction identifies an $\mathscr M$-algebra structure
with a cochain map $a:\mathbb T_{\mathscr M}(X)\to X$ satisfying
$a\eta_X=\id_X$ and $a\mathbb T_{\mathscr M}(a)=a\mu_X$:
these are the unit and composition identities, while the cochain-map
condition is the first identity in \eqref{eq_alpha_compatibility}.

\begin{defn}\label{defn_quiver_operad_algebra_morphisms}
Let $(X,a_X)$ and $(Y,a_Y)$ be algebras over a dg $Q$-operad
$\mathscr P$. A \textit{morphism} of $\mathscr P$-algebras
is a family of cochain maps preserving every operation. Equivalently,
\begin{equation}\label{eq_quiver_operad_algebra_morphisms}
 \operatorname{Hom}_{\operatorname{Alg}(\mathscr P)}(X,Y)
 =\left\{f\in\Hom_Q(X,Y)\,\middle|\,
 f\circ a_X=a_Y\circ\mathbb T_{\mathscr P}(f)\right\}.
\end{equation}
Identities and composition are inherited from $Q$-indexed families,
giving the category $\operatorname{Alg}(\mathscr P)$. The displayed
Hom is a set of algebra morphisms, not an internal Hom complex.
\end{defn}

For a $Q$-indexed family $V$, the family $\mathbb T_{\mathscr P}(V)$ with the structure map
$ \mu_V:\mathbb T_{\mathscr P} (\mathbb T_{\mathscr P}(V) )
 \to \mathbb T_{\mathscr P}(V)$
is called the \textit{free $\mathscr P$-algebra} on $V$.
Thus there is a natural bijection
\begin{equation}\label{eq_free_algebra_universal_property}
 \operatorname{Hom}_{\operatorname{Alg}(\mathscr P)}
 (\mathbb T_{\mathscr P}(V),Y)\cong\Hom_Q(V,Y).
\end{equation}

\subsection{Coalgebras over quiver-cooperads}
\label{subsec:coalgebras_E_cooperads}

\subsubsection{Coalgebras and cofree coalgebras}
\label{sss_comonad_str_induced}
For a dg $Q$-cooperad $\mathscr C$, put
\begin{equation*} \mathbb T^{\mathscr C}(X)(e_0)
 =\bigoplus_{(\vec e;e_0)\in\mathcal P_Q}
     C(\vec e;e_0)\otimes X^*(\vec e).
\end{equation*}
This carries the tensor differential and acts on cochain maps by
applying them to every $X$-factor.

The coproduct and counit induce
$\Delta_X:\mathbb T^{\mathscr C}(X)\to
\mathbb T^{\mathscr C}(\mathbb T^{\mathscr C}(X))$ and
$\varepsilon_X:\mathbb T^{\mathscr C}(X)\to X$.
If $\vec e=\vec g_1*\cdots*\vec g_k$, partition
$\vec x=x_1\otimes\cdots\otimes x_n$ into the corresponding
blocks $\vec x_i$. Writing the component of $\Delta(c)$ as
$\sum c_{(0)}\otimes c_{(1)}\otimes\cdots\otimes c_{(k)}$, one has
\begin{equation*} \Delta_X(c\otimes\vec x)
 =\sum(-1)^{\sum_{i<j}|c_{(j)}||\vec x_i|}
 c_{(0)}\otimes(c_{(1)}\otimes\vec x_1)\otimes\cdots
              \otimes(c_{(k)}\otimes\vec x_k).
\end{equation*}
The sum includes all admissible decompositions. The map
$\varepsilon_X$ is $c\otimes x\mapsto\varepsilon_e(c)x$ on
$C(e;e)\otimes X(e)$ and zero on other summands.
Coassociativity and the counit identities give the comonad identities:
expanding either composite amounts to the same three-level
decomposition.

\begin{defn}\label{defn:coalgebra_over_E_cooperad}
	Coalgebras over $\mathscr C$ are coalgebras for this comonad.
A \textit{$\mathscr C$-coalgebra} is a $Q$-indexed family $X$ and a
cochain map $\rho:X\to\mathbb T^{\mathscr C}(X)$ such that
$
 \Delta_X\rho=\mathbb T^{\mathscr C}(\rho)\rho$ and $ \varepsilon_X\rho=\id_X$.
A \textit{morphism} $f:(X,\rho_X)\to(Y,\rho_Y)$ is a cochain map satisfying
$\rho_Yf=\mathbb T^{\mathscr C}(f)\rho_X$.
\end{defn}
This is the direct-sum convention for cooperadic coalgebras.
The same comonad supplies the cofree construction: for every $X$, the coalgebra
$(\mathbb T^{\mathscr C}(X),\Delta_X)$ is the \textit{cofree $\mathscr C$-coalgebra} on $X$ for this convention:
\[
 \operatorname{Hom}_{\mathscr C\text{-coalg}}
   (Y,\mathbb T^{\mathscr C}(X))\cong\Hom_Q(Y,X).
\]
The forward map is corestriction by $\varepsilon_X$, and the inverse
sends $f$ to $\mathbb T^{\mathscr C}(f)\rho_Y$.
The comonad identities show that these maps are mutually inverse.

\subsubsection{Coderivations and codifferentials}
\label{sss_coderivation_codifferential}
For a homogeneous map $D:X\to X$ of degree $r$, let
$\mathscr C\langle D\rangle$ denote the \textit{one-factor insertion map}
\[
 \mathscr C\langle D\rangle(c\otimes x_1\otimes\cdots\otimes x_n)
 =\sum_i(-1)^{r(|c|+\sum_{j<i}|x_j|)}
 c\otimes x_1\otimes\cdots\otimes D(x_i)\otimes\cdots\otimes x_n.
\]
This is distinct from applying the functor $\mathbb T^{\mathscr C}$
to a degree-zero map, which acts on every factor.
\begin{defn}\label{defn:C_coderivation}
An \textit{ordinary coderivation} of degree $r$ is a homogeneous map
$D$ of the underlying graded family satisfying
$\rho D=\mathscr C\langle D\rangle\rho$.
A \textit{codifferential compatible with $\delta_C$} is a degree-one
map $d$ with $d^2=0$ and
\[
 \rho d=(\delta_C\otimes\id+\mathscr C\langle d\rangle)\rho.
\]
Thus the latter condition says exactly that $\rho$ is a cochain map
for the total tensor differential. When $\delta_C=0$, it reduces to
an ordinary square-zero coderivation.
\end{defn}
Ordinary coderivations form a graded Lie algebra under commutator.
Indeed, expanding two one-factor insertion maps cancels terms acting
on distinct factors in the graded commutator, leaving the one-factor
insertion of the commutator.

On the cofree coalgebra, corestriction gives an isomorphism
\begin{equation}\label{eq_Coder_Hom_Q}
 \operatorname{Coder}(\mathbb T^{\mathscr C}(X),\Delta_X)
 \xrightarrow{\ \cong\ }
 \underline{\operatorname{Hom}}_Q(\mathbb T^{\mathscr C}(X),X),
 \qquad D\longmapsto\varepsilon_XD.
\end{equation}
For a homogeneous $f$ the inverse applies $\Delta_X$, then applies $f$
to exactly one inner $\mathbb T^{\mathscr C}(X)$-factor and
$\varepsilon_X$ to all others, summing over the choice of slot, with
the graded tensor signs. Coassociativity verifies the coderivation
identity. Corestricting this formula yields $f$ by counitality.
Conversely, applying $\mathbb T^{\mathscr C}(\varepsilon_X)$ to the
coderivation identity recovers this formula, proving uniqueness.

Write $d_{C(X)}$ for the tensor differential on
$\mathbb T^{\mathscr C}(X)$. A compatible total differential is
$d_{C(X)}+D$ for an ordinary degree-one coderivation $D$, provided
$(d_{C(X)}+D)^2=0$. When $\mathscr C$ is coaugmented, its coaugmentation
embeds $X$ in the cofree coalgebra. Requiring the total differential
to restrict to the given $d_X$ is equivalent to $D|_X=0$.

\subsection{Pullbacks and change of quiver}
\label{ss_change_quiver}

We now consider how these constructions behave when the underlying quiver changes.

Let $f:Q\to Q'$ be a quiver morphism, writing $f$ for both its
vertex and arrow maps. It sends a path $\vec e$ to the path
$f\vec e$, with $f\varnothing_v=\varnothing_{f(v)}$, and induces
$f:\mathcal P_Q\to\mathcal P_{Q'}$.
For a dg $Q'$-collection $S$, define its \textit{pullback} by
\begin{equation}\label{eq_change_quiver_pullback}
 (f^*S)(\vec e;e_0)=S(f\vec e;f(e_0)).
\end{equation}
Maps and differentials are pulled back componentwise. Thus $f^*$ is
exact, preserves direct sums and quasi-isomorphisms, and commutes with
suspensions. The identities $(g\circ f)^*=f^*g^*$ and
$\id_Q^*=\id$ hold with these definitions.

Applying $f$ to the intermediate arrows of a substitution gives
natural comparison maps
\begin{equation}\label{eq_change_quiver_comparison}
 \chi_{S,T}:f^*S\circ_Q f^*T\longrightarrow
 f^*(S\circ_{Q'}T).
\end{equation}
On each summand of \eqref{eq_S_1_circ_S_2}, this is the identity on
the ordered tensor factors, placed in the summand indexed by the
image arrows. There is also a natural map
\begin{equation}\label{eq_change_quiver_unit}
 \eta_f:\mathbb I_Q\longrightarrow f^*\mathbb I_{Q'},
 \qquad 1_e\longmapsto 1_{f(e)}.
\end{equation}
Both ways of applying these maps to a three-level substitution send
each summand to the same image summand, with the same tensor order.
The unit comparisons have the same property.
In particular, a dg $Q'$-operad $\mathscr P$ pulls back to a dg
$Q$-operad $f^*\mathscr P$: its composition is $f^*\gamma\circ\chi$
and its unit is $f^*\iota\circ\eta_f$. Equivalently, its partial
compositions and identities are those of $\mathscr P$ at the image
profiles. This is the quiver version of change of colours for
operads; compare \cite[\S~1.6]{BergerMoerdijkColoured}.

The comparisons \eqref{eq_change_quiver_comparison} and
\eqref{eq_change_quiver_unit} need not be invertible, so pullback of
cooperads and augmentations is not automatic. For example, identifying
two loops $a,b$ at one vertex creates mutually inverse unary operations
in $f^*\mathbb I_{Q'}$ between $a$ and $b$, whereas $\mathbb I_Q(a;b)=0$;
hence $f^*\mathbb I_{Q'}$ admits no augmentation to $\mathbb I_Q$.
The lifting condition of \S~\ref{ss_change_quiver_trees} makes both
comparisons invertible, allowing cooperad pullback and preservation of
augmentations.

\begin{defn}\label{defn_operad_morphism_over_quiver}
A \textit{morphism of dg quiver-operads}
$(f,\phi):(Q,\mathscr O)\to(Q',\mathscr P)$ consists of a quiver
morphism $f$ and a morphism of dg $Q$-operads
$\phi:\mathscr O\to f^*\mathscr P$.
For $(g,\psi):(Q',\mathscr P)\to(Q'',\mathscr R)$, composition is
\[
 (g,\psi)\circ(f,\phi)
   =(g\circ f,\,f^*\psi\circ\phi).
\]
The identity is $(\id_Q,\id_{\mathscr O})$. This generalizes the notion
of morphism of dg $Q$-operads from \S~\ref{ss_operad_over_graphs},
which is recovered when $f=\id_Q$.
\end{defn}

For a $Q'$-indexed family $X$, define the $Q$-indexed family $f^*X$
by $(f^*X)(e)=X(f(e))$ for each arrow $e\in Q$.
The endomorphism complexes and their compositions give a natural
isomorphism of dg $Q$-operads
\begin{equation}\label{eq_change_quiver_endomorphism}
 f^*\mathbf{End}_{Q'}(X)
 \cong\mathbf{End}_Q(f^*X).
\end{equation}
Indeed, at $(e_1,\ldots,e_n;e_0)$ both sides are
$\underline{\Hom}(\bigotimes_i X(f(e_i)),X(f(e_0)))$,
with the same differential and insertion signs, including $n=0$.

\begin{prop}\label{prop_change_quiver_algebras}
Let $\Phi=(f,\phi):(Q,\mathscr O)\to(Q',\mathscr P)$ be a morphism of
dg quiver-operads. Every $\mathscr P$-algebra $(X,\alpha)$ induces
an $\mathscr O$-algebra on $f^*X$ with structure map
\[
 \mathscr O\xrightarrow{\phi}f^*\mathscr P
 \xrightarrow{f^*\alpha}f^*\mathbf{End}_{Q'}(X)
 \xrightarrow{\cong}\mathbf{End}_Q(f^*X).
\]
On algebra morphisms, set $(\Phi^*h)_e=h_{f(e)}$. This defines
the \textit{restriction functor}
\[
 \Phi^*:\operatorname{Alg}(\mathscr P)
 \longrightarrow\operatorname{Alg}(\mathscr O).
\]
These functors respect identities and composition of morphisms of
quiver-operads.
\end{prop}
\begin{proof}
Each arrow is a dg $Q$-operad morphism, so their composite is an
algebra structure. Reindexing an algebra morphism preserves
the operation identities, so it defines a functor. The asserted
compatibility follows from
$(g\circ f)^*=f^*g^*$ and Definition~\ref{defn_operad_morphism_over_quiver}.
\end{proof}

\begin{prop}[Extension and restriction of algebras]
\label{prop_change_quiver_algebra_adjunction}
For every morphism $\Phi=(f,\phi):(Q,\mathscr O)\to(Q',\mathscr P)$
of dg quiver-operads, the restriction functor $\Phi^*$ admits a left
adjoint $\Phi_!$, called \textit{extension of algebras}. Thus there
are bijections, natural in $X\in\operatorname{Alg}(\mathscr O)$ and
$Y\in\operatorname{Alg}(\mathscr P)$,
\begin{equation}\label{eq_change_quiver_algebra_adjunction}
 \operatorname{Hom}_{\operatorname{Alg}(\mathscr P)}(\Phi_!X,Y)
 \cong
 \operatorname{Hom}_{\operatorname{Alg}(\mathscr O)}(X,\Phi^*Y).
\end{equation}
\end{prop}
\begin{proof}
For the underlying $Q$-indexed family of $X$, put
\[
 (f_\#X)(e')=\bigoplus_{\substack{e\in Q\\f(e)=e'}}X(e),
 \qquad e'\in Q'.
\]
Form the free algebra $F=\mathbb T_{\mathscr P}(f_\#X)$ of
\eqref{eq_free_algebra_universal_property}, with canonical cochain maps
$\iota_e:X(e)\to F(f(e))$. Define $\Phi_!X$ as the quotient of $F$
by the dg $\mathscr P$-algebra ideal generated by
\begin{equation}\label{eq_change_quiver_algebra_relations}
 \phi(o)_F(\iota_{e_1}x_1,\ldots,\iota_{e_n}x_n)
 -\iota_{e_0}\bigl(o_X(x_1,\ldots,x_n)\bigr),
\end{equation}
for homogeneous $o\in\mathscr O(e_1,\ldots,e_n;e_0)$ and
$x_i\in X(e_i)$, including $n=0$. Here $o_X$ and $\phi(o)_F$ use
the notation \eqref{eq_algebra_operation_notation}. An ideal is a family of
subcomplexes stable under operations with at least one input in that family;
the quotient therefore inherits a dg $\mathscr P$-algebra structure.

By freeness and the bijection
$\Hom_{Q'}(f_\#X,Y)\cong\Hom_Q(X,f^*Y)$, a $\mathscr P$-algebra
map $F\to Y$ is equivalent to a family of cochain maps
$u:X\to f^*Y$. It factors through the quotient precisely when
$u$ preserves the $\mathscr O$-operations. This proves the natural
bijection and hence the adjunction; compare
\cite[\S\S~1.6 and~4]{BergerMoerdijkColoured}.
\end{proof}

The \textit{adjunction unit}
$\eta_X^\Phi:X\to\Phi^*(\Phi_!X)$ is induced componentwise by the
canonical generator maps $X(e)\to\mathbb T_{\mathscr P}(f_\#X)(f(e))$
followed by the quotient to $(\Phi_!X)(f(e))$.
The unit need not be an isomorphism:
extension does not in general preserve the original underlying family.

\section{Free quiver-operads and cofree quiver-cooperads}
\label{s_rooted_planar_diagram}

\subsection{Rooted planar diagrams over quivers}
\label{ss_rooted_planar_diagram}
\label{s_correspondence_tree_diagram}
\label{sss_gluing_planar_diagram}
\label{s_canonical_order}
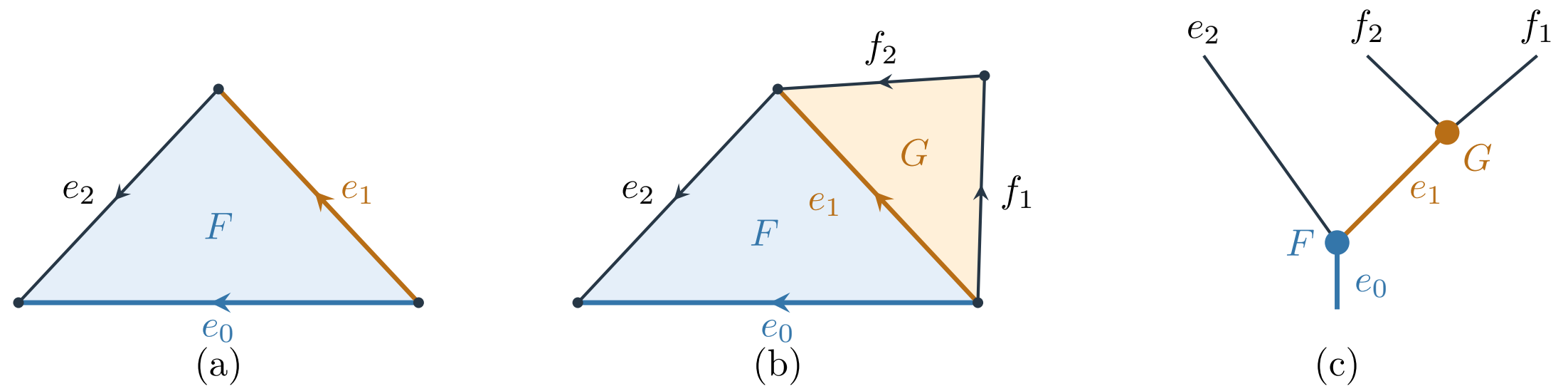
\begin{figure}[htbp]
\centering
\begin{tikzpicture}[x=.62cm,y=.62cm,font=\small,
  edge/.style={draw=planarink,line width=.8pt},
  arrow/.style={postaction={decorate},decoration={markings,
    mark=at position .52 with {\arrow{stealth}}}},
  root/.style={draw=planarblue,line width=1.3pt},
  seam/.style={draw=planarorange,line width=1.2pt},
  dot/.style={circle,fill=planarink,inner sep=1pt}]
\definecolor{planarink}{HTML}{273746}
\definecolor{planarblue}{HTML}{3476AA}
\definecolor{planarorange}{HTML}{B86E16}
\definecolor{planaroldface}{HTML}{E5EFF9}
\definecolor{planarnewface}{HTML}{FFF0D9}

\begin{scope}
\coordinate (A) at (0,0); \coordinate (B) at (6,0);
\coordinate (C) at (3,3.2);
\fill[planaroldface] (A)--(B)--(C)--cycle;
\draw[root,arrow] (B)--(A);
\draw[seam,arrow] (B)--(C);
\draw[edge,arrow] (C)--(A);
\foreach \p in {A,B,C} \node[dot] at (\p) {};
\node[text=planarblue] at (3,1.15) {$F$};
\node[right,text=planarorange] at (4.65,1.65) {$e_1$};
\node[left] at (1.35,1.65) {$e_2$};
\node[below,text=planarblue] at (3,-.08) {$e_0$};
\node at (3,-.95) {(a)};
\end{scope}

\begin{scope}[xshift=5.2cm]
\coordinate (A) at (0,0); \coordinate (B) at (6,0);
\coordinate (C) at (3,3.2); \coordinate (D) at (6.1,3.4);
\fill[planaroldface] (A)--(B)--(C)--cycle;
\fill[planarnewface] (B)--(D)--(C)--cycle;
\draw[root,arrow] (B)--(A);
\draw[edge,arrow] (C)--(A);
\draw[edge,arrow] (B)--(D);
\draw[edge,arrow] (D)--(C);
\draw[seam,arrow] (B)--(C);
\foreach \p in {A,B,C,D} \node[dot] at (\p) {};
\node[text=planarblue] at (2.8,1.05) {$F$};
\node[text=planarorange] at (5.05,2.25) {$G$};
\node[left,text=planarorange] at (4.15,1.47) {$e_1$};
\node[left] at (1.35,1.65) {$e_2$};
\node[right] at (6.15,1.65) {$f_1$};
\node[above] at (4.55,3.38) {$f_2$};
\node[below,text=planarblue] at (3,-.08) {$e_0$};
\node at (3,-.95) {(b)};
\end{scope}

\begin{scope}[xshift=10.4cm]
\coordinate (F) at (3,.9); \coordinate (G) at (4.65,2.55);
\draw[root] (3,-.1)--(F);
\draw[edge] (F)--(1,3.7);
\draw[seam] (F)--(G);
\draw[edge] (G)--(3.45,3.7);
\draw[edge] (G)--(6,3.7);
\node[circle,fill=planarblue,inner sep=2.3pt] at (F) {};
\node[circle,fill=planarorange,inner sep=2.3pt] at (G) {};
\node[left,text=planarblue] at (2.85,.9) {$F$};
\node[below right,text=planarorange] at (4.7,2.55) {$G$};
\node[right,text=planarorange] at (3.9,1.65) {$e_1$};
\node[above] at (1,3.7) {$e_2$};
\node[above] at (3.45,3.7) {$f_2$};
\node[above] at (6,3.7) {$f_1$};
\node[right,text=planarblue] at (3.08,.25) {$e_0$};
\node at (3,-.95) {(c)};
\end{scope}
\end{tikzpicture}
\caption{Constructing a rooted planar diagram. (a) An elementary face
$F$ of profile $(e_1,e_2;e_0)$. (b) Attaching a face $G$ of profile
$(f_1,f_2;e_1)$ along $e_1$ gives outer profile $(f_1,f_2,e_2;e_0)$:
$e_1$ becomes internal, while $e_0$ remains the root. (c) The
corresponding tree has vertices $F$ and $G$ joined by $e_1$, with
leaves read from right to left.}
\label{fig:planar_diagram_attachment}
\end{figure}

\noindent\textit{Definition.}
Rooted planar diagrams encode compositions by attachment along
matching edges, as illustrated in Figure~\ref{fig:planar_diagram_attachment}.
We define them as follows.

\begin{itemize}
\item \emph{Elementary faces.}
For a profile $(e_1,\ldots,e_n;e_0)\in\mathcal P_Q$ with $n\geq1$, an
\textit{elementary diagram} is a polygon bounded by a directed input path
of colors $e_1,\ldots,e_n$ and a directed output edge of color $e_0$
with the same source and target. Vertex colors are the corresponding sources
and targets in $Q_0$; equal colors do not identify distinct vertices.
For $n=1$, this is a genuine two-edge face.
For $n=0$, an elementary diagram is a formal face with one output edge
and no input edges, called a \textit{nullary face}.

\item \emph{Attachment.}
A \textit{rooted planar diagram over $Q$} is built from one elementary
diagram by finitely many attachments of fresh elementary diagrams.
Each attachment matches the new output edge to an exposed input edge
of the same color, retaining the joining edge and replacing that input
by the new ordered input list. For polygonal faces, this is exterior
gluing of matching edges and their corresponding endpoints only.

\item \emph{Root and inputs.}
The initial output edge is the \textit{root edge}, and its incident face
is the \textit{root face}. The unattached input edges form the
\textit{input path}, read from right to left with the root drawn below it.
Attachment preserves admissibility of the outer profile.

\item \emph{Identity diagrams.}
For each arrow $e$, also admit the \textit{formal identity diagram}
$D_e$: a single edge with input and output $e$ and no faces. This is
distinct from a unary face, even one of profile $(e;e)$.
\end{itemize}

Diagrams are considered up to isomorphisms of their rooted planar trees
of faces, preserving arrow colors;
attachment histories are not part of the data.

These diagrams correspond to \textit{admissible $Q$-trees}: finite planar rooted
trees whose edges are colored by arrows of $Q$ and whose ordered input
and output colors form a profile in $\mathcal P_Q$ at every vertex.
Faces correspond to tree vertices, and $D_e$ to the vertex-free tree of
color $e$. For positive-arity faces, this is planar duality
(Figure~\ref{figure:tree_to_planar_diagram}). The inverse attaches the
elementary diagrams of a tree starting at its root vertex.

\begin{figure}[htbp]
\centering
\resizebox{\linewidth}{!}{\begin{tikzpicture}[font=\normalsize,>=Stealth,
 edge/.style={line width=1.05pt},
 flow/.style={postaction={decorate},decoration={markings,mark=at position .52 with {\arrow{Stealth}}}}]
\definecolor{nullaryink}{HTML}{273746}
\definecolor{nullaryblue}{HTML}{286D9C}
\definecolor{nullaryorange}{HTML}{BE711B}
\definecolor{nullarygray}{HTML}{64717C}
\definecolor{nullaryfaceblue}{HTML}{E5EFF9}
\definecolor{nullaryfaceorange}{HTML}{FFF0D9}
\definecolor{nullaryneutral}{HTML}{F4F6F7}
\path[use as bounding box] (.4,1) rectangle (18,5.7);

\newcommand{\nullarybaseface}{%
\coordinate (P0) at (3,0);\coordinate (P1) at (3,1.8);
\coordinate (P2) at (1.8,2.6);\coordinate (P3) at (0,1.8);\coordinate (P4) at (0,0);
\fill[nullaryfaceblue] (P0)--(P1)--(P2)--(P3)--(P4)--cycle;

\draw[edge,flow,nullaryblue] (P0)--(P1);
\draw[edge,flow,nullaryorange] (P1)--(P2);
\draw[edge,flow,nullaryblue] (P2)--(P3);
\draw[edge,flow,nullaryblue] (P3)--(P4);
\draw[edge,flow,nullaryblue] (P0)--(P4);

\foreach \p in {P0,P1,P2,P3,P4} \fill[nullaryink] (\p) circle (1.4pt);
\node[nullaryink,right=4pt] at (P0) {$u$};
\node[nullaryink,right=4pt] at (P1) {$v$};
\node[nullaryink,above=4pt] at (P2) {$v$};
\node[nullaryink,left=4pt] at (P3) {$w$};
\node[nullaryink,left=4pt] at (P4) {$z$};
\node[nullaryblue,right=4pt] at (3,.9) {$e_1$};
\node[nullaryorange] at (2.35,1.94) {$e_2$};
\node[nullaryblue,above=4pt] at (.9,2.2) {$e_3$};
\node[nullaryblue,left=4pt] at (0,.9) {$e_4$};
\node[nullaryblue,below=4pt] at (1.5,0) {$e_0$};

\node[nullaryblue,font=\Large] at (1.5,1.18) {$F$};

}
\begin{scope}[shift={(1.2,2.15)}]
\nullarybaseface
\end{scope}
\node[text=nullaryink] at (2.7,1.3) {$F$};

\draw[nullarygray,-{Stealth[length=2.5mm]},line width=.9pt] (5.0,3.5)--(6.5,3.5);
\node[text=nullarygray,font=\small] at (5.75,3.85) {glue $N$};

\begin{scope}[shift={(7.6,2.15)}]
\fill[nullaryfaceorange] (3,1.8)--(1.8,2.6) .. controls (2.5,3.7) and (4.2,3.4) .. (3,1.8);
\nullarybaseface
\draw[nullarygray,densely dashed,line width=.8pt] (1.8,2.6) .. controls (2.5,3.7) and (4.2,3.4) .. (3,1.8);
\node[nullaryorange,font=\Large] at (2.95,2.86) {$N$};
\end{scope}
\node[text=nullaryink] at (9.15,1.3) {$F\circ_2N$};

\draw[nullarygray,-{Stealth[length=2.5mm]},line width=.9pt] (11.8,3.5)--(13.1,3.5);
\node[text=nullarygray,font=\small] at (12.45,3.85) {shrink $e_2$};

% Only the boundary is drawn in the contraction schematic.
\begin{scope}[shift={(14.1,2.15)}]
\coordinate (Q0) at (3,0);\coordinate (Qv) at (1.7,2.4);
\coordinate (Q3) at (0,1.8);\coordinate (Q4) at (0,0);
\fill[nullaryneutral] (Q0)--(Qv)--(Q3)--(Q4)--cycle;
\draw[edge,flow,nullaryblue] (Q0)--(Qv);
\draw[edge,flow,nullaryblue] (Qv)--(Q3);
\draw[edge,flow,nullaryblue] (Q3)--(Q4);
\draw[edge,flow,nullaryblue] (Q0)--(Q4);
\fill[nullaryfaceorange] (Qv) circle (4pt);
\draw[nullaryorange,line width=.6pt] (Qv) circle (4pt);
\foreach \p in {Q0,Qv,Q3,Q4} \fill[nullaryink] (\p) circle (1.4pt);
\node[nullaryink,right=4pt] at (Q0) {$u$};
\node[nullaryink,above=5pt] at (Qv) {$v$};
\node[nullaryink,left=4pt] at (Q3) {$w$};
\node[nullaryink,left=4pt] at (Q4) {$z$};
\node[nullaryblue,right=4pt] at (2.35,1.2) {$e_1$};
\node[nullaryblue,above=4pt] at (.85,2.1) {$e_3$};
\node[nullaryblue,left=4pt] at (0,.9) {$e_4$};
\node[nullaryblue,below=4pt] at (1.5,0) {$e_0$};

\end{scope}
\node[text=nullaryink] at (15.6,1.3) {Outer boundary};
\end{tikzpicture}}
\caption{Nullary attachment. Gluing $N:(\varnothing_v;e_2)$ to
$F:(e_1,e_2,e_3,e_4;e_0)$ gives outer profile $(e_1,e_3,e_4;e_0)$,
retaining $N$ and the internal edge $e_2$. The right panel illustrates
only the boundary change by shrinking $e_2$; the dashed cap outline
is not a diagram edge.}
\label{fig:nullary_face_attachment}
\end{figure}
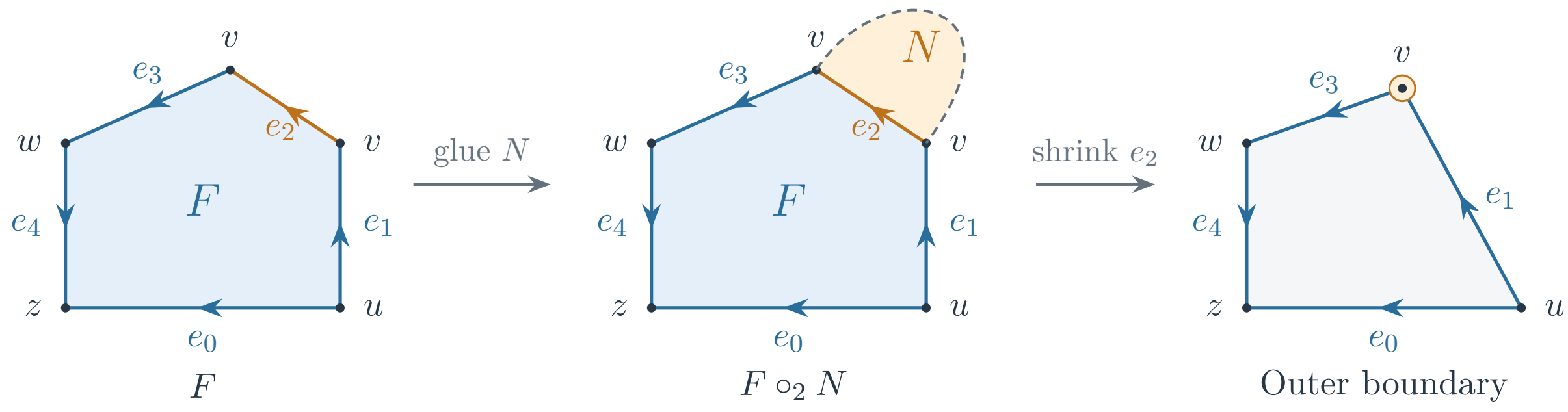

\medskip
\noindent\textit{Gluing.}
Let $F$ and $G$ be rooted planar diagrams. If the output color of $G$
matches the $i$-th input color of $F$, write $F\circ_iG$ for their
attachment along these edges;
Figure~\ref{fig:planar_diagram_attachment}(b) illustrates $F\circ_1G$.
The faces of $F\circ_iG$ are the disjoint union of those of $F$
and $G$; its outer profile is obtained by replacing the $i$-th input of
$F$ by the input path of $G$. If $G$ is nullary, the matched input
disappears from the outer profile while the joining edge remains
internal (Figure~\ref{fig:nullary_face_attachment}).
For $F$ with $n$ inputs and diagrams $G_i$
whose output colors match the respective input colors of $F$, write
$F\circ(G_1,\ldots,G_n)$ for simultaneous attachment at all inputs.
Under the tree correspondence these operations are graftings, hence
are associative with units $D_e$.

\Needspace{5\baselineskip}
\medskip
\noindent\textit{Canonical ordering.}
Order tree vertices by root-first depth-first traversal, visiting input
branches in their indexed order, from right to left in the drawings.
This gives the \textit{canonical face order}
\[
 \operatorname{Ord}(D)=
 (F_0,\operatorname{Ord}(D_1),\ldots,\operatorname{Ord}(D_k)),
\]
where $F_0$ is the root face and $D_i$ is the subdiagram attached to its
$i$-th input;
a bare input gives $D_i=D_{e_i}$ and contributes no face. This recursive
rule is unambiguous and invariant under diagram isomorphisms.

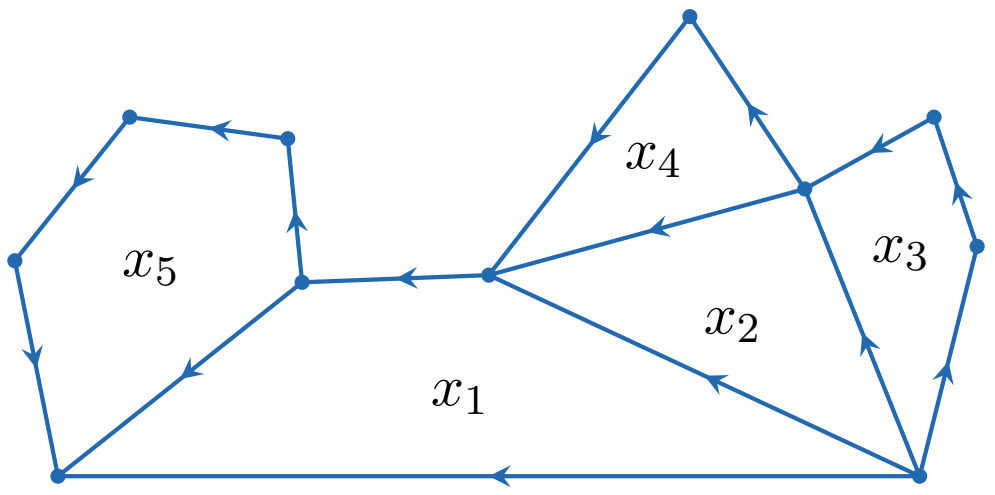
\begin{figure}[htbp]
\centering
\begin{tikzpicture}[scale=.95,xscale=-1,>=stealth,
  midarrow/.style={postaction={decorate},
    decoration={markings,mark=at position .5 with {\arrow{stealth}}}}]
\definecolor{diagramblue}{RGB}{35,105,175}
\coordinate (A) at (0,0); \coordinate (B) at (3,1.4);
\coordinate (C) at (6,0); \coordinate (D) at (.8,2);
\coordinate (E) at (-.4,1.6); \coordinate (F) at (1.6,3.2);
\coordinate (G) at (4.4,2.35); \coordinate (H) at (4.3,1.35);
\coordinate (I) at (-.1,2.5); \coordinate (J) at (5.5,2.5);
\coordinate (K) at (6.3,1.5);
\foreach \a/\b in {A/C,A/B,B/H,H/C,A/D,D/B,A/E,E/I,I/D,D/F,F/B,H/G,G/J,J/K,K/C}
  \draw[diagramblue,thick,midarrow] (\a)--(\b);
\foreach \v in {A,B,C,D,E,F,G,H,I,J,K}\fill[diagramblue] (\v) circle(1.5pt);
\node at (3.2,.55) {$x_1$}; \node at (1.3,1.05) {$x_2$};
\node at (.13,1.55) {$x_3$}; \node at (1.85,2.2) {$x_4$};
\node at (5.35,1.45) {$x_5$};
\end{tikzpicture}
\caption{Canonical face order in a five-face diagram: visit the root face first, then its input branches from right to left.}
\label{figure_decorated_diagram}
\end{figure}

\subsection{Free quiver-operads}
\label{sss_free_Q_operads}
Let $S$ be a graded $Q$-collection. For a diagram $D$ with ordered faces
$F_1,\ldots,F_r$ and their profiles $(\vec f_j;f_{0j})$, put
\begin{equation}\label{eq_decorated_diagram_collection}
 H_{S,D}=\bigotimes_{j=1}^r S(\vec f_j;f_{0j}),\qquad
 H_S(\vec e;e_0)=\bigoplus_{[D]\text{ of profile }(\vec e;e_0)}H_{S,D}.
\end{equation}
Set $H_{S,D_e}=\Bbbk$. A homogeneous \textit{decorated diagram} is written
$(D;x_1\otimes\cdots\otimes x_r)$, subject to multilinearity in its labels.

Define the \textit{free $Q$-operad} $\mathsf F(S)$ on the collection $H_S$ by grafting and graded
reordering of labels. Explicitly, if the faces of $D'$ are inserted
after the first $p$ faces of $D$, then
\[
\begin{aligned}
 &(D;x_1,\ldots,x_r)\circ_i(D';y_1,\ldots,y_q)\\
 &\quad=(-1)^{(\sum_b|y_b|)(\sum_{a>p}|x_a|)}
 (D\circ_iD';x_1,\ldots,x_p,y_1,\ldots,y_q,x_{p+1},\ldots,x_r).
\end{aligned}
\]

% BEGIN TEMPORARY FIGURE: Free quiver-operad composition.
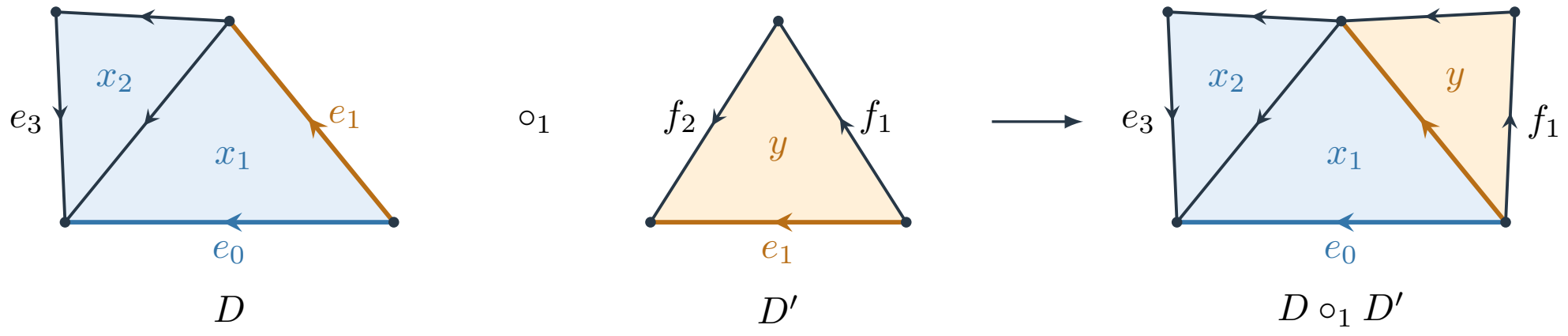
\begin{figure}[htbp]
\centering
\scalebox{0.9}{% Temporary illustration for Section 3.3.a; shared by the paper and preview.
\begin{tikzpicture}[x=.85cm,y=.85cm,font=\small,
  edge/.style={draw=freeink,line width=.8pt},
  midarrow/.style={postaction={decorate},decoration={markings,
    mark=at position .52 with {\arrow{stealth}}}},
  root/.style={draw=freeblue,line width=1.2pt},
  seam/.style={draw=freeorange,line width=1.2pt},
  dot/.style={circle,fill=freeink,inner sep=1pt}]
\definecolor{freeink}{HTML}{273746}
\definecolor{freeblue}{HTML}{3476AA}
\definecolor{freeorange}{HTML}{B86E16}
\definecolor{freeoldface}{HTML}{E5EFF9}
\definecolor{freenewface}{HTML}{FFF0D9}

% The right input of the root face is input 1. The old second
% branch carries x_2, so grafting at input 1 precedes x_2 in face order.
\begin{scope}
\coordinate (A) at (0,0); \coordinate (B) at (3.6,0);
\coordinate (C) at (1.8,2.2); \coordinate (L) at (-.1,2.3);
\fill[freeoldface] (A)--(B)--(C)--cycle;
\fill[freeoldface] (A)--(C)--(L)--cycle;
\draw[root,midarrow] (B)--(A);
\draw[seam,midarrow] (B)--(C);
\draw[edge,midarrow] (C)--(A);
\draw[edge,midarrow] (C)--(L);
\draw[edge,midarrow] (L)--(A);
\foreach \p in {A,B,C,L} \node[dot] at (\p) {};
\node[text=freeblue] at (1.85,.7) {$x_1$};
\node[text=freeblue] at (.55,1.55) {$x_2$};
\node[right,text=freeorange] at (2.75,1.15) {$e_1$};
\node[above] at (.85,2.3) {$e_2$};
\node[left] at (-.1,1.1) {$e_3$};
\node[below,text=freeblue] at (1.8,-.08) {$e_0$};
\node at (1.8,-.95) {$D$};
\end{scope}

\node at (5.15,1.1) {$\circ_1$};

\begin{scope}[xshift=5.45cm]
\coordinate (A) at (0,0); \coordinate (B) at (2.8,0);
\coordinate (C) at (1.4,2.2);
\fill[freenewface] (A)--(B)--(C)--cycle;
\draw[seam,midarrow] (B)--(A);
\draw[edge,midarrow] (B)--(C);
\draw[edge,midarrow] (C)--(A);
\foreach \p in {A,B,C} \node[dot] at (\p) {};
\node[text=freeorange] at (1.4,.8) {$y$};
\node[right] at (2.15,1.15) {$f_1$};
\node[left] at (.65,1.15) {$f_2$};
\node[below,text=freeorange] at (1.4,-.08) {$e_1$};
\node at (1.4,-.95) {$D'$};
\end{scope}

\draw[->,freeink,line width=.7pt] (10.15,1.1)--(11.15,1.1);

\begin{scope}[xshift=10.35cm]
\coordinate (A) at (0,0); \coordinate (B) at (3.6,0);
\coordinate (C) at (1.8,2.2); \coordinate (L) at (-.1,2.3);
\coordinate (R) at (3.7,2.3);
\fill[freeoldface] (A)--(B)--(C)--cycle;
\fill[freeoldface] (A)--(C)--(L)--cycle;
\fill[freenewface] (B)--(R)--(C)--cycle;
\draw[root,midarrow] (B)--(A);
\draw[edge,midarrow] (C)--(A);
\draw[edge,midarrow] (C)--(L);
\draw[edge,midarrow] (L)--(A);
\draw[seam,midarrow] (B)--(C);
\draw[edge,midarrow] (B)--(R);
\draw[edge,midarrow] (R)--(C);
\foreach \p in {A,B,C,L,R} \node[dot] at (\p) {};
\node[text=freeblue] at (1.85,.7) {$x_1$};
\node[text=freeblue] at (.55,1.55) {$x_2$};
\node[text=freeorange] at (3.05,1.55) {$y$};
\node[above] at (.85,2.3) {$e_2$};
\node[left] at (-.1,1.1) {$e_3$};
\node[right] at (3.7,1.1) {$f_1$};
\node[above] at (2.75,2.3) {$f_2$};
\node[below,text=freeblue] at (1.8,-.08) {$e_0$};
\node at (1.8,-.95) {$D\circ_1D'$};
\end{scope}

\end{tikzpicture}}
\caption{Partial composition in the free $Q$-operad, with homogeneous
labels. Glue the output of $D'$ to the first input of $D$ along the
highlighted edge of color $e_1$. The canonical face order is then
$x_1,y,x_2$, giving
$(D;x_1,x_2)\circ_1(D';y)=(-1)^{|y||x_2|}(D\circ_1D';x_1,y,x_2)$.
The outer profile changes from
$(e_1,e_2,e_3;e_0)$ to $(f_1,f_2,e_2,e_3;e_0)$.}
\label{fig:free-operad-composition-draft}
\end{figure}
% END TEMPORARY FIGURE.

For $D$ with $n$ inputs, full composition is
\begin{equation*} \gamma_S\bigl((D;\vec x);(D_1;\vec y_1),\ldots,(D_n;\vec y_n)\bigr)
 =\epsilon(\sigma;\vec x,\vec y_1,\ldots,\vec y_n)
   (D\circ(D_1,\ldots,D_n);\sigma(\vec x,\vec y_1,\ldots,\vec y_n)).
\end{equation*}
Here $\sigma$ puts \emph{all} labels in canonical face order and
$\epsilon(\sigma;-)$ is its \textit{Koszul sign}: each interchange of homogeneous
labels $a,b$ contributes $(-1)^{|a||b|}$.
The unit is $(D_e;1)$. Associativity follows because successive
reorderings have the same total permutation and Koszul sign as the
single final reordering; the vertex-free trees give both unit identities.

The \textit{weight} of a diagram is its number of faces. Thus
$\mathsf F(S)^{(0)}=\mathbb I_Q$, $\mathsf F(S)^{(1)}=S$, and partial
composition adds weights. Only the weight-zero part of $H_S(e;e)$ is
canonically $\Bbbk$; arbitrary unary composites are retained.

If $S$ is dg, define
\begin{equation}\label{eq:free_diff_order}
 \delta(D;x_1,\ldots,x_r)
 =\sum_j(-1)^{\sum_{a<j}|x_a|}
 (D;x_1,\ldots,\delta_Sx_j,\ldots,x_r),\qquad\delta(D_e;1)=0.
\end{equation}
The tensor differential squares to zero, and graded permutations commute
with it. Hence the signed composition satisfies the Leibniz identity,
so $\mathsf F(S)$ is a dg $Q$-operad.

\begin{prop}\label{prop_universal_property_free_dg_Q_operad}
Restriction to weight one gives a natural bijection
\[
 \operatorname{Hom}_{\mathrm{dgOp}_Q}(\mathsf F(S),\mathscr M)
 \cong\operatorname{Hom}_{\mathrm{dgColl}_Q}(S,M).
\]
\end{prop}
\begin{proof}
Given a cochain map $f:S\to M$, evaluate each decorated tree by applying
$f$ to its vertex labels and composing in $\mathscr M$. The operad axioms
and Koszul interchange identity make this independent of the order of
partial graftings. It respects the signed composition and units, and
commutes with differentials by the Leibniz rule. Every tree is an iterated
grafting of corollas, so the extension is unique.
\end{proof}

\subsection{Cofree quiver-cooperads}
\label{s_cofree}
Let $\mathscr C$ be a coaugmented dg $Q$-cooperad, with reduced
part $\overline C$ and projection $p:C\to\overline C$.
For a finite admissible $Q$-tree $T$, let $\overline\Delta_T$ be its full iterated decomposition followed by $p$ at every vertex; this is well defined by coassociativity. The cooperad is \textit{conilpotent} if each $c\in\overline C$ admits $N(c)$ such that $\overline\Delta_T(c)=0$ whenever $T$ has more than $N(c)$ vertices, including unary vertices.

For a graded $Q$-collection $S$, use the decorated-diagram collection
$H_S$ of \S~\ref{sss_free_Q_operads} to define the \textit{cofree $Q$-cooperad} $\mathsf F_c(S)$ by
cutting trees. An \textit{admissible two-level decomposition} is
$D=D_0\circ(D_1,\ldots,D_k)$ with $k\geq0$,
where $D_0$ is a rooted subtree containing the root and retaining all
input edges of its vertices, and the $D_i$ are
the full branches above its inputs; formal identity diagrams are allowed.
Also include $D_0=D_{e_0},D_1=D$. Let $u_j$ denote the tensor of labels
belonging to $D_j$ in their inherited canonical order, with $u_j=1$
when $D_j$ is an identity diagram. Define
\begin{equation}\label{eq_comult_full}
 \Delta(D;\vec x)=\sum_{D=D_0\circ(D_1,\ldots,D_k)}
 \epsilon(\sigma;\vec x)\,(D_0;u_0)\otimes\cdots\otimes(D_k;u_k),
\end{equation}
where $\sigma$ reorders the original face list into the concatenated
lists for $D_0,D_1,\ldots,D_k$. The counit projects onto the
$D_e$-summands, and the coaugmentation includes them.
Figure~\ref{fig:cofree-cooperad-decomposition} illustrates one summand
and its Koszul sign.

\begin{figure}[htbp]
\centering
% Two-level decomposition for Section 3.3.b.
\begingroup
\definecolor{cofreeink}{HTML}{273746}
\definecolor{cofreebasefill}{HTML}{E5EFF9}
\definecolor{cofreebranchfill}{HTML}{FFF0D9}
\definecolor{cofreeleaffill}{HTML}{E4F0E8}
\definecolor{cofreebaseink}{HTML}{3476AA}
\definecolor{cofreebranchink}{HTML}{A86619}
\definecolor{cofreeleafink}{HTML}{397455}
\definecolor{cofreecutink}{HTML}{9D4162}
\begin{tikzpicture}[x=.49cm,y=.49cm,font=\footnotesize,
 edge/.style={draw=cofreeink,line width=.65pt},
 directed/.style={postaction={decorate},decoration={markings,
   mark=at position .52 with {\arrow{stealth}}}},
 root/.style={line width=1.05pt},
 cut/.style={draw=cofreecutink,line width=1.15pt,dash pattern=on 2.5pt off 1.4pt},
 dot/.style={circle,fill=cofreeink,inner sep=.8pt},
 seamlabel/.style={fill=white,inner sep=.8pt,text=cofreecutink,font=\scriptsize}]
% D: canonical face order x_1, x_2, x_3, x_4, x_5.
\begin{scope}
\coordinate (A) at (0,0); \coordinate (B) at (6,0);
\coordinate (C) at (3,2.2); \coordinate (E) at (.2,3);
\coordinate (R) at (6,3); \coordinate (T) at (8,2.2);
\coordinate (L) at (2.1,4.5);
\fill[cofreebasefill] (A)--(B)--(C)--cycle;
\fill[cofreebasefill] (A)--(C)--(E)--cycle;
\fill[cofreebranchfill] (B)--(R)--(C)--cycle;
\fill[cofreebranchfill] (B)--(T)--(R)--cycle;
\fill[cofreeleaffill] (C)--(L)--(E)--cycle;
\draw[edge,root,directed] (B)--(A);
\foreach \p/\q in {B/T,T/R,R/C,C/L,L/E,E/A,C/A,B/R}
 \draw[edge,directed] (\p)--(\q);
\draw[cut,directed] (B)--(C);
\draw[cut,directed] (C)--(E);
\foreach \p in {A,B,C,E,R,T,L} \node[dot] at (\p) {};
\node[text=cofreebaseink] at (3,.72) {$x_1$};
\node[text=cofreebranchink] at (4.9,1.88) {$x_2$};
\node[text=cofreebranchink] at (6.75,1.85) {$x_3$};
\node[text=cofreebaseink] at (1.12,1.7) {$x_4$};
\node[text=cofreeleafink] at (1.83,3.28) {$x_5$};
\node[seamlabel] at (4.38,1.2) {$a$};
\node[seamlabel] at (1.53,2.62) {$b$};
\node[left,font=\scriptsize] at (.1,1.4) {$e$};
\node at (3.8,-1.15) {$D$};
\end{scope}
\draw[->,cofreeink,line width=.8pt] (8.7,1.9)--(10.3,1.9);
\node[above,font=\scriptsize,text=cofreecutink] at (9.5,2.0) {cut};
% Bottom/root factor D_0.
\begin{scope}[xshift=5.39cm]
\coordinate (A) at (0,0); \coordinate (B) at (6,0);
\coordinate (C) at (3,2.2); \coordinate (E) at (.2,3);
\fill[cofreebasefill] (A)--(B)--(C)--(E)--cycle;
\draw[edge,root,directed] (B)--(A);
\foreach \p/\q in {B/C,C/E,E/A,C/A} \draw[edge,directed] (\p)--(\q);
\foreach \p in {A,B,C,E} \node[dot] at (\p) {};
\node[text=cofreebaseink] at (3,.72) {$x_1$};
\node[text=cofreebaseink] at (1.12,1.7) {$x_4$};
\node[above right,font=\scriptsize,text=cofreecutink] at (4.5,1.1) {$a$};
\node[above,font=\scriptsize,text=cofreecutink] at (1.6,2.65) {$b$};
\node[left,font=\scriptsize] at (.1,1.4) {$e$};
\node at (3,-1.15) {$D_0$};
\node[text=cofreebaseink] at (3,-1.95) {$u_0=x_1\otimes x_4$};
\end{scope}
\node at (18,1.25) {$;$};
% Upper factor D_1, rotated to put its root below.
\begin{scope}[xshift=9.31cm]
\coordinate (A) at (0,0); \coordinate (B) at (3.8,0);
\coordinate (C) at (1.7,1.8); \coordinate (R) at (3.7,2.8);
\fill[cofreebranchfill] (A)--(B)--(R)--(C)--cycle;
\draw[edge,root,directed] (B)--(A);
\foreach \p/\q in {B/R,R/C,C/A,B/C} \draw[edge,directed] (\p)--(\q);
\foreach \p in {A,B,C,R} \node[dot] at (\p) {};
\node[text=cofreebranchink] at (1.8,.58) {$x_2$};
\node[text=cofreebranchink] at (3.04,1.8) {$x_3$};
\node[below,font=\scriptsize,text=cofreecutink] at (1.9,0) {$a$};
\node at (1.9,-1.15) {$D_1$};
\node[text=cofreebranchink] at (1.9,-1.95) {$u_1=x_2\otimes x_3$};
\end{scope}
\node at (23.6,1.25) {$,$};
% Upper factor D_2.
\begin{scope}[xshift=12.005cm]
\coordinate (A) at (0,0); \coordinate (B) at (2.8,0);
\coordinate (C) at (1.4,2.35);
\fill[cofreeleaffill] (A)--(B)--(C)--cycle;
\draw[edge,root,directed] (B)--(A);
\draw[edge,directed] (B)--(C);
\draw[edge,directed] (C)--(A);
\foreach \p in {A,B,C} \node[dot] at (\p) {};
\node[text=cofreeleafink] at (1.4,.8) {$x_5$};
\node[below,font=\scriptsize,text=cofreecutink] at (1.4,0) {$b$};
\node at (1.4,-1.15) {$D_2$};
\node[text=cofreeleafink] at (1.4,-1.95) {$u_2=x_5$};
\end{scope}
\node at (28.2,1.25) {$,$};
% Formal identity diagram: no faces.
\begin{scope}[xshift=14.21cm]
\draw[edge,root,directed] (1.8,0)--(0,0);
\node[dot] at (0,0) {}; \node[dot] at (1.8,0) {};
\node[above,font=\scriptsize] at (.9,0) {$e$};
\node at (.9,-1.15) {$D_3=D_e$};
\node at (.9,-1.95) {$u_3=1$};
\end{scope}

\node[font=\small] at (15.1,-3.1) {$
 (x_1,x_2,x_3,x_4,x_5)
 \longmapsto (x_1,x_4)\mid(x_2,x_3)\mid(x_5)\mid(),
 \qquad
 \epsilon(\sigma;\vec x)=(-1)^{|x_4|(|x_2|+|x_3|)}.$};
\end{tikzpicture}
\endgroup
\caption{An admissible two-level decomposition
$D=D_0\circ(D_1,D_2,D_e)$ in the full coproduct, with homogeneous
face labels. Cutting along the dashed edges of colors $a,b$ leaves
the root factor $D_0$ and branches $D_1,D_2$; the third input of $D_0$
carries the identity diagram $D_e$, contributing the factor $(D_e;1)$.
The bars separate the canonical label lists defining $u_0,u_1,u_2,u_3$.
Moving $x_4$ past $x_2,x_3$ gives the displayed Koszul sign.
Unmarked edge and vertex colors are suppressed.}
\label{fig:cofree-cooperad-decomposition}
\end{figure}

For the infinitesimal coproduct the formula is
\begin{equation*} \Delta_{(1)}(D;\vec x)=
 \sum_{D=D_0\circ_iD_1}\epsilon(\sigma;\vec x)\,(D_0;u_0)\otimes_i(D_1;u_1).
\end{equation*}
The slot $i$ is retained in the target summand.
The \textit{reduced infinitesimal coproduct} further projects both displayed factors away
from the identity summands; its terms have two positive weights.

Each coproduct is a finite sum because a finite tree has finitely many
cuts. Two successive cuts and a three-level cut index the same data;
their Koszul signs are the sign of the same final permutation. This
proves coassociativity. Cutting immediately below or above all vertices
gives the two counit identities. The differential
\eqref{eq:free_diff_order} commutes with these graded permutations,
making $\mathsf F_c(S)$ a dg cooperad when $S$ is dg. Its reduced
iterated tree decompositions have at most as many nonidentity vertices
as the weight, so it is conilpotent.

\begin{prop}\label{prop_cofree_Q_cooperad}
The cooperad $\mathsf F_c(S)$ is cofree among coaugmented conilpotent
dg $Q$-cooperads. More precisely, corestriction to its weight-one
cogenerators gives a natural bijection
\[
 \operatorname{Hom}_{\mathrm{dgCoOp}^{\mathrm{conil}}_Q}
   (\mathscr D,\mathsf F_c(S))
 \cong\operatorname{Hom}_{\mathrm{dgColl}_Q}(\overline D,S),
\]
where morphisms on the left preserve the coaugmentation and
$\overline D=\ker\varepsilon_D$.
\end{prop}
\begin{proof}
We can adopt the standard cofree construction
\cite[Theorem~5.8.3 and \S~5.9.17]{loday2012algebraic} in planar
diagrams. Given a cochain map $f:\overline D\to S$, preserve
coaugmentations and, for each nonidentity admissible diagram $A$
with $r$ faces, set
$ F_A=f^{\otimes r}\overline\Delta_A$,
where $\overline\Delta_A$ is the iterated decomposition along $A$,
projected to $\overline D$ at every face and written in canonical
face order.
Conilpotence and finite support of iterated coproducts ensure that
$F$ lands in $H_S$. Coassociativity makes it compatible with diagram
decompositions, and it commutes with differentials. Conversely,
decomposing into individual faces and projecting each factor to
weight one forces this formula for any extension of $f$, proving
uniqueness.
\end{proof}

\subsection{Diagram lifting under a change of quiver}
\label{ss_change_quiver_trees}
\label{ss_change_quiver_diagrams}
Let $f:Q\to Q'$ be a quiver morphism, writing $f$ for its maps on
vertices and arrows. Applying $f$ to the colors of a rooted planar
diagram over $Q$ gives a rooted planar diagram over $Q'$. This
changes only the colors: the faces, their incidences, the root, and
the input order are retained.
For a profile $\pi=(e_1,\ldots,e_n;e_0)$, write
$f\pi=(fe_1,\ldots,fe_n;fe_0)$, and let
$\operatorname{Diag}_Q(\pi)$ be the set of isomorphism classes of
rooted planar diagrams of outer profile $\pi$ in the sense of
\S~\ref{ss_rooted_planar_diagram}.

\begin{defn}\label{defn_unique_tree_lifting}
\label{defn_unique_diagram_lifting}
We say that $f$ has \textit{unique diagram lifting} if, for every
$\pi\in\mathcal P_Q$, the relabeling map
\begin{equation}\label{eq_unique_diagram_lifting}
 f_*:\operatorname{Diag}_Q(\pi)
 \longrightarrow\operatorname{Diag}_{Q'}(f\pi),
 \qquad [D]\longmapsto[fD],
\end{equation}
is bijective. In other words, every target diagram with boundary $f\pi$ has
exactly one lift, up to diagram isomorphism, with the same rooted
planar shape and the prescribed source boundary $\pi$.
\end{defn}

Recall that a quiver is \textit{transitive} if arrows $a\to b$ and
$b\to c$ imply the existence of an arrow $a\to c$, \textit{thin}
if each ordered pair of vertices admits at most one arrow, and
\textit{reflexive} if every vertex has a loop.

\begin{lem}\label{lem_change_quiver_thin}
Let $f:Q\to Q'$ be a quiver morphism with $Q$ transitive and thin
and $Q'$ thin. Then $f$ has unique diagram lifting for diagrams
without nullary faces. If $Q$ is also reflexive, this holds for all
diagrams.
\end{lem}
\begin{proof}
Fix $\pi=(e_1,\ldots,e_n;e_0)$ and a diagram $D'$ with boundary
$(f(e_1),\ldots,f(e_n);f(e_0))$. To lift $D'$, label each edge by
the unique arrow between the endpoints in $\pi$ of its consecutive
input block. This arrow exists by transitivity and thinness of $Q$,
using reflexivity for empty blocks. The labels give admissible
face profiles and boundary $\pi$, and map to the prescribed colours
by thinness of $Q'$. Since every label is forced, the lift is unique,
including for identity diagrams.
\end{proof}

\begin{ex}\label{ex_change_quiver_terminal_folding}
\begin{enumerate}[(i)]
\item If $Q$ is transitive and thin, the map $q:Q\to\ast$ to the
terminal quiver has unique diagram lifting for diagrams without
nullary faces, by Lemma~\ref{lem_change_quiver_thin}. If $Q$ is also
reflexive, this holds for all diagrams.
\item Let $S_A,S_B$ be disjoint copies of a set $S$. The triangular
quiver with vertex set $S_A\sqcup S_B$ and arrows
$ (S_A\times S_A)\sqcup(S_A\times S_B)\sqcup(S_B\times S_B)$
maps to the complete quiver on $S$ by forgetting the copy labels.
The source is reflexive, transitive, and thin, so the lemma gives
unique diagram lifting for this folding.
\item Let $Q$ be reflexive, transitive, and thin, and form $Q^{\mathbb N}$ by
replacing each arrow $e_{ab}$ with arrows $e^p_{ab}$, $p\geq0$.
Let $\ast^{\mathbb N}$ have one vertex and loops $e^p$, $p\geq0$.
The map
\[
 q:Q^{\mathbb N}\longrightarrow\ast^{\mathbb N},
 \qquad e^p_{ab}\longmapsto e^p,
\]
has unique diagram lifting: each edge of a lift is uniquely determined
by its consecutive input block's endpoints and its prescribed label $p$;
reflexivity supplies the loops for empty blocks.
\end{enumerate}
\end{ex}

\begin{prop}\label{prop_change_quiver_strong_monoidal}
If $f$ has unique diagram lifting, the comparison maps
\eqref{eq_change_quiver_comparison} and
\eqref{eq_change_quiver_unit} are isomorphisms:
\begin{equation}\label{eq_change_quiver_strong_monoidal}
 \chi_{S,T}: f^*S\circ_Q f^*T\cong f^*(S\circ_{Q'}T),
 \qquad \quad \eta_f: \mathbb I_Q\cong f^*\mathbb I_{Q'}.
\end{equation}
\end{prop}
\begin{proof}
Fix an outer profile $\pi=(\vec e;e_0)$. A summand of
$(f^*S\circ_Q f^*T)(\pi)$ is specified by a decomposition
$\vec e=\vec g_1*\cdots*\vec g_k$ and intermediate arrows
$a_1,\ldots,a_k$. It is
$ S(fa_1,\ldots,fa_k;fe_0)
 \otimes\bigotimes_{i=1}^k T(f\vec g_i;fa_i)$.
The same tensor product occurs in $(S\circ_{Q'}T)(f\pi)$,
indexed by the image diagram.
On this summand,
$ \chi_{S,T}$ is given by $x\otimes y_1\otimes\cdots\otimes y_k
\mapsto x\otimes y_1\otimes\cdots\otimes y_k$.
Unique lifting pairs every target summand with exactly one source
summand, so $\chi_{S,T}$ is invertible.
For a unary profile $(e_1;e_0)$, if $fe_1=fe_0$, the target identity
diagram lifts to an identity diagram with boundary $(e_1;e_0)$,
forcing $e_1=e_0$. Thus the two unit collections have the same
nonzero components, on which $\eta_f$ is $\id_{\Bbbk}$.
\end{proof}

\begin{prop}\label{prop_change_quiver_augmented_pullback}
If $f$ has unique diagram lifting, $f^*$ pulls back augmented dg
$Q'$-operads and coaugmented conilpotent dg $Q'$-cooperads to
augmented dg $Q$-operads and coaugmented conilpotent dg $Q$-cooperads,
respectively. For the underlying collection $X$
of either such object, we have
$\overline{f^*X}\cong f^*\overline X$.
\end{prop}
\begin{proof}
First, for an augmented dg $Q'$-operad $\mathscr M$, use the pullback
operad $f^*\mathscr M$ defined in \S~\ref{ss_change_quiver}.
The map $\eta_f$ is the unit of $f^*\mathbb I_{Q'}$, hence an
operad map, and Proposition~\ref{prop_change_quiver_strong_monoidal}
makes it an isomorphism. Thus
$\epsilon_{f^*M}=\eta_f^{-1}\circ f^*\epsilon_M:
f^*M\to\mathbb I_Q$
is an operad morphism splitting the pulled-back unit, hence an augmentation.
Second, for a cooperad with underlying collection $C$, use
Proposition~\ref{prop_change_quiver_strong_monoidal} to define the coproduct
$\Delta_{f^*C}:f^*C\to f^*C\circ_Q f^*C$ by
$\Delta_{f^*C}=\chi_{C,C}^{-1}\circ f^*\Delta_C$.
Transport the counit and coaugmentation through $\eta_f$. The associativity and unit compatibilities of \(\chi\) and \(\eta_f\) transfer the coassociativity and counit identities from \(C\) to \(f^*C\), and preserve the coaugmentation identities.
For each rooted planar diagram $D$, the reduced iterated
coproducts satisfy, under the tensor identifications,
$ \overline\Delta_D^{f^*C}=\overline\Delta_{fD}^{C}$.
Relabeling preserves the number of faces and conilpotence.
In both cases, pulling back $X=\mathbb I_{Q'}\oplus\overline X$
gives $f^*X\cong\mathbb I_Q\oplus f^*\overline X$, the desired
reduced-part identity.
\end{proof}

In the next statement we display the quiver in the notation for
the free operad and cofree cooperad, writing $\mathsf F_Q$ and
$\mathsf F_{c,Q}$, respectively.

\begin{prop}\label{prop_change_quiver_free_cofree}
If $f$ has unique diagram lifting, every dg $Q'$-collection $S$ has
natural weight-preserving isomorphisms
\begin{equation}\label{eq_change_quiver_free_cofree}
 \mathsf F_Q(f^*S)\cong f^*\mathsf F_{Q'}(S),
 \qquad
 \mathsf F_{c,Q}(f^*S)\cong f^*\mathsf F_{c,Q'}(S),
\end{equation}
of dg operads and coaugmented conilpotent dg cooperads, respectively.
On weight one they are the identity of $f^*S$.
\end{prop}
\begin{proof}
With the notation of \eqref{eq_decorated_diagram_collection}, define
$ \Phi_S:H_{f^*S} \to  f^*H_S$, $ (D;\vec x)\mapsto(fD;\vec x)$.
Unique lifting and $H_{f^*S,D}=H_{S,fD}$ make $\Phi_S$ a natural,
weight-preserving isomorphism of complexes. Relabeling preserves
grafting, cuts, and their Koszul signs, so
$ \Phi_S(a\circ_i b)=\Phi_S(a)\circ_i\Phi_S(b)$ and $ \Delta'\Phi_S=(\Phi_S\circ_Q\Phi_S)\Delta$,
where $\Delta'$ is the pulled-back coproduct. Finally,
$\Phi_S(D_e;1)=(D_{f(e)};1)$ and
$\Phi_S^{(1)}=\id_{f^*S}$, giving compatibility with units,
counits, and coaugmentations and the weight-one assertion.
\end{proof}

For collections vanishing in arity zero,
Propositions~\ref{prop_change_quiver_strong_monoidal}--\ref{prop_change_quiver_free_cofree}
remain valid when unique
lifting is required only for diagrams without nullary faces,
including identity diagrams.

\section{Generalized constructions and theories}
\label{s_generalized}

\subsection{Convolution and twisting morphisms}
\label{ss_convolution_twisting}

Let $\mathscr C=(C,\Delta,\varepsilon;\delta_C)$ be a dg $Q$-cooperad and $\mathscr M=(M,\gamma,\iota;\delta_M)$ a dg $Q$-operad. Their \textit{total convolution complex} is the product of the profilewise Hom complexes:
\begin{equation}\label{eq_total_convolution_complex}
\mathfrak g_{\mathscr C,\mathscr M}
 =\prod_{(\vec e;e_0)\in\mathcal P_Q}
 \underline{\Hom}\bigl(C(\vec e;e_0),M(\vec e;e_0)\bigr).
\end{equation}
Its differential is $
\partial\Phi=\delta_M\Phi-(-1)^{|\Phi|}\Phi\delta_C$.
By \eqref{eq_gamma_(1)} and \eqref{eq_Delta_(1)}, we set
$\Phi\star\Psi=\gamma_{(1)}(\Phi\otimes\Psi)\Delta_{(1)}$.
This is the infinitesimal product induced by the convolution operad structure on the internal Hom collection \cite{berger2003axiomatic,loday2012algebraic}.
Thus a term $c_0\otimes c_1$ in the infinitesimal decomposition associated to $\vec e=\vec h\circ_i\vec g$ contributes
$(-1)^{|\Psi||c_0|}\Phi(c_0)\circ_i\Psi(c_1)$.
The expression is defined even for infinitely many profiles: on each element the cooperadic decomposition has finite support.

\begin{defn}\label{def:fc_twisting_morphism}
Let $\mathscr C$ be coaugmented and let $\mathscr M$ be a dg $Q$-operad, not necessarily augmented. A \textit{twisting morphism} is a degree-one element $\tau\in\mathfrak g_{\mathscr C,\mathscr M}$ such that $\tau\iota_C=0$ and $\partial\tau+\tau\star\tau=0$.
Write $\operatorname{Tw}_Q(\mathscr C,\mathscr M)$ for this set. If $\mathscr M$ is augmented, write $\operatorname{Tw}^{\mathrm{aug}}_Q(\mathscr C,\mathscr M)$ for its subset defined by $\epsilon_M\tau=0$.
\end{defn}
%For the terminal quiver these are the usual normalized twisting morphisms for nonsymmetric operads. The distinction between the two sets is needed when the target is an endomorphism operad.

Let $\tau\in\operatorname{Tw}_Q(\mathscr C,\mathscr M)$ and let $d_0$ be the ordinary tensor differential on $C\circ_Q M$.
For a homogeneous element $c\otimes(m_1,\ldots,m_k)$ and each infinitesimal decomposition of $c$ with inner factor occupying the consecutive block $i,\ldots,i+\ell-1$ (possibly empty, with $\ell=0$), write $\Delta_i(c)=\sum c_0\otimes c_1$. Define $\vec m= (m_1,\dots, m_k)$ and
\begin{equation*} 
 d_\tau^r\bigl(c\otimes \vec m \bigr) =\sum(-1)^{|c_0|+(|c_1|+1)\sum_{a<i}|m_a|}\,c_0\otimes
 \bigl(m_1,\ldots,m_{i-1}, \gamma\bigl(\tau(c_1);m_i,\ldots,m_{i+\ell-1}\bigr),
 m_{i+\ell},\ldots,m_k\bigr).
\end{equation*}
The sum includes all homogeneous terms of admissible decompositions, with their coefficients. The sign comes from moving $\tau$ past $c_0$ and $\tau(c_1)$ past the preceding $m_a$. Set
\begin{equation*}\mathscr C\circ_\tau\mathscr M:=(C\circ_Q M,d_0+d_\tau^r)
\end{equation*}
and call it the \textit{right twisted composite product}. In the expansion of $(d_0+d_\tau^r)^2$, terms from insertions in distinct outer slots cancel; the remaining terms are induced by $\partial\tau+\tau\star\tau$ and vanish by the Maurer--Cartan equation.

\subsection{Bar and cobar constructions over quivers}
\label{ss_bar_cobar}
\subsubsection{The cobar construction}
\label{sss_cobar}
Let $\mathscr C$ be a coaugmented conilpotent dg $Q$-cooperad. Its \textit{cobar construction} is the augmented dg operad, denoted by $\Omega_Q\mathscr C=(\mathsf F(s\overline C),d_\Omega)$.
The differential is the derivation $d_\Omega=d_1+d_2$, with $d_1(sc)=-s\delta_Cc$ and quadratic part
\begin{equation}\label{eq_cobar_d2_explicit}
 d_2(sc)=-\sum_{i,(c)}(-1)^{|c_0|}(sc_0)\circ_i(sc_1),
 \qquad \overline\Delta_{(1)}(c)=\sum_{i,(c)}c_0\otimes c_1.
\end{equation}
The minus sign is part of our cohomological convention. The augmentation kills the generators. Compatibility of $\delta_C$ with decomposition gives $d_1d_2+d_2d_1=0$. In $d_2^2$, the two orders of cutting any two internal edges cancel with their tensor signs, so $d_\Omega^2=0$.
The degree-one map
\begin{equation*}\iota_{\mathscr C}:\mathscr C\longrightarrow\Omega_Q\mathscr C,
\qquad \iota_{\mathscr C}|_{\mathbb I_Q}=0,\quad \iota_{\mathscr C}(c)=sc
\end{equation*}
is a twisting morphism: \eqref{eq_cobar_d2_explicit} gives
$\partial\iota_{\mathscr C}=-\iota_{\mathscr C}\star\iota_{\mathscr C}$ on generators.

\subsubsection{The bar construction}
\label{subsec:fc_bar}
For an augmented dg $Q$-operad $\mathscr M$, define its \textit{bar construction} by $B_Q\mathscr M=(\mathsf F_c(s^{-1}\overline M),d_B)$, with differential $d_B=d_1+d_2$.
Here $d_1$ is induced by $d_1(s^{-1}m)=-s^{-1}\delta_Mm$. The coderivation $d_2$ has only a weight-two corestriction. On the two-vertex tree with outer decoration $s^{-1}m$ and inner decoration $s^{-1}m'$ in input $i$, it is
\begin{equation}\label{eq_bar_d2_corestriction}
 d_2\bigl((s^{-1}m)\circ_i(s^{-1}m')\bigr)
 =(-1)^{|m|}s^{-1}(m\circ_i m').
\end{equation}
The notation on the left denotes a decorated tree, not a cooperadic multiplication. On larger trees $d_2$ sums contractions of internal edges, with the tensor signs from the canonical ordering. Associativity pairs the two orders of contraction, proving $d_2^2=0$; the dg Leibniz rule gives $d_1d_2+d_2d_1=0$. Thus $d_B^2=0$.
Projection to cogenerators followed by suspension defines
\begin{equation}\label{eq_canonical_bar_map}
\pi_{\mathscr M}:B_Q\mathscr M\longrightarrow s^{-1}\overline M
 \xrightarrow{s}\overline M\hookrightarrow M.
\end{equation}
It belongs to $\operatorname{Tw}^{\mathrm{aug}}_Q(B_Q\mathscr M,\mathscr M)$: the weight-one equation follows from $d_1$, the weight-two equation from the defining formula for $d_2$, and both sides vanish on all other weights.

\subsubsection{Representability of twisting morphisms}
Write $\mathrm{dgOp}^{\mathrm{aug}}_Q$ for augmented dg operads with augmentation-preserving maps, and $\mathrm{dgCoOp}^{\mathrm{conil}}_Q$ for coaugmented conilpotent dg cooperads with coaugmentation-preserving maps.
\begin{thm}[Twisting morphisms and the bar--cobar adjunction]\label{thm_bar_cobar_adjunction}
Let $\mathscr C\in\mathrm{dgCoOp}^{\mathrm{conil}}_Q$.
\begin{enumerate}[(i)]
\item For any dg $Q$-operad $\mathscr M$, there is a natural bijection
\begin{equation}\label{eq:fc_cobar_representability}
\Hom_{\mathrm{dgOp}_Q}(\Omega_Q\mathscr C,\mathscr M)
 \cong\operatorname{Tw}_Q(\mathscr C,\mathscr M),\qquad F\longmapsto F\iota_{\mathscr C}.
\end{equation}
\item If $\mathscr M$ is augmented, this restricts to augmentation-preserving maps and gives natural bijections
\begin{equation*}\Hom_{\mathrm{dgOp}^{\mathrm{aug}}_Q}(\Omega_Q\mathscr C,\mathscr M)
 \cong\operatorname{Tw}^{\mathrm{aug}}_Q(\mathscr C,\mathscr M)
 \cong\Hom_{\mathrm{dgCoOp}^{\mathrm{conil}}_Q}(\mathscr C,B_Q\mathscr M).
\end{equation*}
Consequently, $\Omega_Q$ is left adjoint to $B_Q$.
\end{enumerate}
\end{thm}
\begin{proof}
Freeness identifies a graded operad map $F$ with a degree-one map $\tau:\overline C\to M$ via $F(sc)=\tau(c)$. Extend $\tau$ by zero on $\mathbb I_Q$. On a generator, $\delta_MF=Fd_\Omega$ is exactly
$\delta_M\tau+\tau\delta_C+\gamma_{(1)}(\tau\otimes\tau)\overline\Delta_{(1)}=0$.
The derivation rule then gives the chain-map identity everywhere, proving (i). For augmented $\mathscr M$, augmentation preservation is equivalent to $\epsilon_M\tau=0$.

For the second bijection in (ii), the cofree property identifies a coaugmentation-preserving graded cooperad map $G$ with its corestriction $\overline C\to s^{-1}\overline M$. Write this map as $s^{-1}\tau$. The two maps $d_BG$ and $G\delta_C$ are coderivations along $G$, so their equality is determined by corestriction. The defining formula for $d_2$ and \eqref{eq_canonical_bar_map} turn this equality into $\partial\tau+\tau\star\tau=0$. Conilpotence ensures that the tree expansion of $G$ is finite on each element. Both correspondences are natural in $\mathscr C$ and $\mathscr M$, and hence give the adjunction.
\end{proof}

\subsection{Quadratic quiver-(co)operads}
\label{ss_quadratic}

A \textit{quadratic datum} is a graded $Q$-collection $S$ and a graded subcollection $R\subseteq\mathsf F(S)^{(2)}$. A \textit{dg quadratic datum} additionally has a degree-one differential $\delta_S:S\to S$ whose induced differential on $\mathsf F(S)$ preserves $R$. This differential preserves the weight grading.

A \textit{$Q$-ideal} $J\subseteq M$ is a subcollection such that $u\circ_i v\in J$ whenever either factor lies in $J$. A \textit{dg ideal} is also stable under the differential.
The componentwise kernel of a dg $Q$-operad morphism is a dg $Q$-ideal.
For any dg ideal $J$, the induced compositions and differential make
$M/J$ a dg $Q$-operad, and every dg $Q$-operad map annihilating $J$
factors uniquely through $M/J$.

For a quadratic datum $(S,R)$, let $(R)$ be the ideal generated by $R$. The \textit{quadratic $Q$-operad} of $(S,R)$ is $\mathscr M(S,R)=\mathsf F(S)/(R)$. It is weight graded, with components $\mathbb I_Q$ and $S$ in weights zero and one. Its canonical augmentation is projection to weight zero. For a dg datum the induced differential preserves $(R)$.
\begin{prop}\label{prop_dg_quadratic_operad_generators}
For a dg quadratic datum, dg operad maps $\mathscr M(S,R)\to\mathscr N$ correspond bijectively to degree-zero cochain maps of collections $\varphi:S\to N$ whose free extensions $\widehat\varphi:\mathsf F(S)\to\mathscr N$ vanish on $R$.
\end{prop}
\begin{proof}
Apply the free universal property and the quotient universal property.
\end{proof}

The \textit{quadratic cooperad} $\mathscr C(S,R)$ is the largest weight-graded subcooperad of $\mathsf F_c(S)$ on which
\begin{equation}\label{eq:quadratic_corelation_condition}
\mathsf F_c(S)\xrightarrow{\operatorname{pr}^{(2)}}\mathsf F_c(S)^{(2)}
 \longrightarrow\mathsf F_c(S)^{(2)}/R
\end{equation}
vanishes. It is coaugmented and conilpotent. Its low weights are $\mathscr C(S,R)^{(0)}=\mathbb I_Q$, $\mathscr C(S,R)^{(1)}=S$, and $\mathscr C(S,R)^{(2)}=R$.
For $r\geq3$ it can be computed by
\begin{equation}\label{eq_quadratic_co_recursive_full}
\mathscr C(S,R)^{(r)}=
\left\{u\in\mathsf F_c(S)^{(r)}:
\overline\Delta_{(1)}u\in
\bigoplus_{\substack{p+q=r\\p,q>0}}
\mathscr C(S,R)^{(p)}\circ_{(1)}\mathscr C(S,R)^{(q)}\right\}.
\end{equation}
The sum includes all admissible profiles and insertion positions. Condition \eqref{eq_quadratic_co_recursive_full} is equivalent to requiring the decoration of every adjacent pair isolated by tree cuts to lie in $R$. These conditions are stable under further cuts and characterize the largest subcooperad satisfying \eqref{eq:quadratic_corelation_condition}. For dg data, they are also preserved by the differential.

\subsection{Koszul duality over quivers}
\label{ss_koszul_duality}
Throughout the Koszul theory below, $(S,R)$ is a dg quadratic datum as in \S~\ref{ss_quadratic}. The internal differential $\delta_S$ may be nonzero, and cohomological degrees of generators may be arbitrary. The weight grading is independent of the cohomological grading and is preserved by the differential.

The signed identification of weight-two collections $j_2:\mathsf F(S)^{(2)}\longrightarrow\mathsf F_c(s^{-1}S)^{(2)}$ is given by $j_2(x\circ_i y)=(-1)^{|x|}(s^{-1}x)\circ_i(s^{-1}y)$.
It has degree $-2$ and is the tensor product of the two desuspension maps, in outer--inner order. Equip $s^{-1}S$ with its desuspension differential. Then $j_2$ commutes with the induced differentials, so $s^{-2}R:=j_2(R)$ is a subcomplex. The \textit{Koszul dual dg cooperad} of $\mathscr M=\mathscr M(S,R)$ is
\begin{equation}\label{eq_koszul_dual_definition}
 \mathscr C=\mathscr M^{\antishriek}:=\mathscr C(s^{-1}S,s^{-2}R),
\end{equation}
with the differential $\delta_C$ induced from $\mathsf F_c(s^{-1}S)$. In particular,
\begin{equation}\label{eq_koszul_dual_internal_differential}
 \delta_C(s^{-1}x)=-s^{-1}\delta_Sx
 \qquad (x\in S).
\end{equation}
Projection to weight one defines the canonical degree-one map
\begin{equation*}\kappa:\mathscr M^{\antishriek}\longrightarrow s^{-1}S
 \xrightarrow{s}S\hookrightarrow\mathscr M.
\end{equation*}
\begin{prop}\label{prop:koszul_twisting_fc}
The map $\kappa$ is in $\operatorname{Tw}^{\mathrm{aug}}_Q(\mathscr M^{\antishriek},\mathscr M)$.
\end{prop}
\begin{proof}
It vanishes off weight one, so $\kappa\star\kappa$ vanishes off weight two. On $j_2(x\circ_i y)$ its value is $-x\circ_i y$. Thus its restriction to $j_2(R)$ is minus the inclusion of $R$ followed by the quadratic quotient, and is zero. Moreover, $\partial\kappa=\delta_M\kappa+\kappa\delta_C=0$: on $s^{-1}x$, the two terms are $\delta_Sx$ and $-\delta_Sx$, and both vanish in other weights because the internal differentials preserve weight.
\end{proof}

The \textit{right Koszul complex} is
$K^r(\mathscr M)=\mathscr M^{\antishriek}\circ_\kappa\mathscr M
 = (\mathscr M^{\antishriek}\circ_Q\mathscr M,d_0+d_\kappa^r )$,
where $d_0$ is the tensor differential induced by $\delta_C$ and $\delta_M$, as in \S~\ref{ss_convolution_twisting}.
The counit and augmentation induce a chain map
$\epsilon_\kappa:K^r(\mathscr M)\to\mathbb I_Q$.
Here and below, \textit{acyclicity} of an augmented complex means that this map is a quasi-isomorphism, not that its weight-zero homology vanishes.

A weight-graded (co)operad is \textit{connected} if its weights are nonnegative and its weight-zero component is $\mathbb I_Q$. The differential preserves weight; composition adds weights, and decomposition preserves total weight.

\begin{prop}\label{thm_koszul_criterion_E_operad}
For the dg quadratic operad $\mathscr M=\mathscr M(S,R)$, the following are equivalent:
\begin{enumerate}[(i)]
\item $\epsilon_\kappa:\mathscr M^{\antishriek}\circ_\kappa\mathscr M\to\mathbb I_Q$ is a quasi-isomorphism;
\item the map $\iota_\kappa:\mathscr M^{\antishriek}\to B_Q\mathscr M$ corresponding to $\kappa$ is a quasi-isomorphism;
\item the map $\rho_\kappa:\Omega_Q\mathscr M^{\antishriek}\to\mathscr M$ corresponding to $\kappa$ is a quasi-isomorphism.
\end{enumerate}
\end{prop}
\begin{proof}
This is the quiver version of the operadic twisting morphism theorem
\cite[Theorem~6.6.2]{loday2012algebraic}, applied to $\kappa$.
By Proposition~\ref{prop_extension_by_zero}, extension by zero identifies
$Q$-(co)operads with nonsymmetric colored (co)operads supported on
admissible profiles and preserves the bar, cobar, and twisted composite
constructions. The proof of the cited theorem carries over to these
colored planar trees: its tree operations preserve admissibility, and
its homological comparisons use exact tensor products and direct sums
over $\Bbbk$, with no symmetric-group coinvariants.

For convergence, work at fixed total weight $N$, filtering by outer
cooperad weight in the comparisons and by the number of nonunit
decorations in the universal resolutions. Connectedness bounds both
indices by $N$, and the internal differentials preserve these
filtrations. Thus the argument applies in arbitrary cohomological
degrees and yields precisely (i)$\Leftrightarrow$(ii)$\Leftrightarrow$(iii).
\end{proof}
\begin{defn}\label{defn_koszul_fc_multicategory}
The dg quadratic operad $\mathscr M$ is \textit{Koszul} if these equivalent conditions hold. Its \textit{Koszul resolution} is $\rho_\kappa:\mathscr M_\infty\to\mathscr M$, where $\mathscr M_\infty:=\Omega_Q\mathscr M^{\antishriek}$.
\end{defn}
The differential of $\mathscr M_\infty$ is $d_\Omega=d_1+d_2$ as in \S~\ref{sss_cobar}, with $d_1(sc)=-s\delta_Cc$. If $\delta_S=0$, then $\delta_C=0$ and $d_\Omega=d_2$.

\subsection{Change of quiver and Koszul resolutions}
\label{ss_change_quiver_koszul}

Throughout this subsection, let $f:Q\to Q'$ have unique diagram lifting
in the sense of Definition~\ref{defn_unique_diagram_lifting}.
Proposition~\ref{prop_change_quiver_augmented_pullback} supplies the
pullback of augmented operads and coaugmented conilpotent cooperads,
with their reduced parts identified by
$\overline{f^*M}=f^*\overline M$ and
$\overline{f^*C}=f^*\overline C$.
For collections vanishing in arity zero, it suffices to require unique
diagram lifting for diagrams without nullary faces.

\begin{prop}[Bar--cobar compatibility]
\label{prop_change_quiver_bar_cobar}
For an augmented dg $Q'$-operad $\mathscr M$ and a coaugmented
conilpotent dg $Q'$-cooperad $\mathscr C$, there are natural
isomorphisms
\begin{equation}\label{eq_change_quiver_bar_cobar}
\begin{aligned}
 B_Q(f^*\mathscr M)&\cong f^*(B_{Q'}\mathscr M),\\
 \Omega_Q(f^*\mathscr C)&\cong f^*(\Omega_{Q'}\mathscr C).
\end{aligned}
\end{equation}
They preserve the (co)augmentations and identify the canonical
twisting maps $\pi_{f^*\mathscr M}$ and $\iota_{f^*\mathscr C}$ with
$f^*\pi_{\mathscr M}$ and $f^*\iota_{\mathscr C}$, respectively.
\end{prop}
\begin{proof}
Apply Proposition~\ref{prop_change_quiver_free_cofree} to
$s^{-1}\overline M$ and $s\overline C$, obtaining the graded
comparisons $\Phi_B$ and $\Phi_\Omega$ in
\eqref{eq_change_quiver_bar_cobar}.
Write $d_1+d_2$ and $d_1'+d_2'$ for the source and target
differentials, with $d_1$ the internal part, and let $p'$ project
to the target bar cogenerators.
Naturality of $d_1$ and
\eqref{eq_cobar_d2_explicit}--\eqref{eq_bar_d2_corestriction}
for the pullback (co)operations give, for $c\in\overline{f^*C}$,
\[
\begin{aligned}
 \Phi_Bd_1&=d_1'\Phi_B,
 &\Phi_\Omega d_1&=d_1'\Phi_\Omega,\\
 p'd_2'\Phi_B&=p'\Phi_Bd_2,
 &d_2'\Phi_\Omega(sc)&=\Phi_\Omega d_2(sc).
\end{aligned}
\]
Uniqueness of coderivations along $\Phi_B$ extends the bar
corestriction identity; the derivation rule extends the cobar
generator identity. Thus both comparisons commute with $d_1+d_2$.
They preserve the (co)augmentations and identify weight-one
decorations, so
\[
 (f^*\pi_{\mathscr M})\Phi_B=\pi_{f^*\mathscr M},
 \qquad
 \Phi_\Omega\iota_{f^*\mathscr C}=f^*\iota_{\mathscr C}.
 \qedhere
\]
\end{proof}

When the quiver must be specified, write $\mathscr M_Q(S,R)$ and
$\mathscr C_Q(S,R)$ for the quadratic constructions of
\S~\ref{ss_quadratic}. For a dg quadratic datum $(S,R)$ over $Q'$,
we regard $f^*R$ as a subcomplex of
$\mathsf F_Q(f^*S)^{(2)}$ through the weight-two isomorphism of
Proposition~\ref{prop_change_quiver_free_cofree}.

\begin{prop}[Quadratic constructions under change of quiver]
\label{prop_change_quiver_quadratic}
Let $(S,R)$ be a dg quadratic datum over $Q'$. Then
$(f^*S,f^*R)$ is a dg quadratic datum over $Q$, and there are
natural weight-preserving isomorphisms
\begin{equation}\label{eq_change_quiver_quadratic}
\begin{aligned}
 \mathscr M_Q(f^*S,f^*R)&\cong f^*\mathscr M_{Q'}(S,R),\\
 \mathscr C_Q(f^*S,f^*R)&\cong f^*\mathscr C_{Q'}(S,R).
\end{aligned}
\end{equation}
of dg operads and dg cooperads, respectively.
\end{prop}
\begin{proof}
Use the weight-preserving dg identifications of
Proposition~\ref{prop_change_quiver_free_cofree}.
Contracting a distinguished relation, lifting the contracted tree
uniquely, and expanding at the lifted vertex gives
$(f^*R)=f^*((R))$. Exactness of $f^*$ then gives
$\mathsf F_Q(f^*S)/(f^*R)
\cong f^*\bigl(\mathsf F_{Q'}(S)/(R)\bigr)$.
Put $C=\mathscr C_Q(f^*S,f^*R)$ and
$C'=\mathscr C_{Q'}(S,R)$. Then $C^{(r)}=f^*(C'^{(r)})$
for $r=0,1,2$. For $r\geq3$, exactness, compatibility with
infinitesimal substitution, and
\eqref{eq_quadratic_co_recursive_full} give inductively, in
$\mathsf F_{c,Q}(f^*S)^{(r)}$,
\[
 f^*(C'^{(r)})
 =\overline\Delta_{(1)}^{-1}\!\left(
   \bigoplus_{\substack{p+q=r\\p,q>0}}
   f^*(C'^{(p)})\circ_{(1)}f^*(C'^{(q)})\right)
 =C^{(r)}.
 \qedhere
\]
\end{proof}

For a dg quadratic $Q'$-operad $\mathscr P=\mathscr M_{Q'}(S,R)$,
recall its Koszul dual dg cooperad
$\mathscr P^{\antishriek}=\mathscr C_{Q'}(s^{-1}S,s^{-2}R)$
from \eqref{eq_koszul_dual_definition}.
Recall also the canonical
twisting morphism $\kappa_{\mathscr P}:\mathscr P^{\antishriek}\to\mathscr P$
of Proposition~\ref{prop:koszul_twisting_fc} and the corresponding
bar and cobar maps of
Proposition~\ref{thm_koszul_criterion_E_operad}(ii)--(iii):
\[
 \iota_{\kappa_{\mathscr P}}:
 \mathscr P^{\antishriek}\longrightarrow B_{Q'}\mathscr P,
 \qquad
 \rho_{\kappa_{\mathscr P}}:
 \Omega_{Q'}\mathscr P^{\antishriek}\longrightarrow\mathscr P.
\]
As in Definition~\ref{defn_koszul_fc_multicategory}, write
$\mathscr P_\infty:=\Omega_{Q'}\mathscr P^{\antishriek}$.

\begin{thm}[Change of quiver and Koszul duality]
\label{thm_change_quiver_koszul}
Let $\mathscr P$ be a dg quadratic $Q'$-operad.
\begin{enumerate}[(i)]
\item There is a natural weight-preserving isomorphism of dg cooperads
\begin{equation}\label{eq_change_quiver_koszul_dual}
 (f^*\mathscr P)^{\antishriek}
 \cong f^*(\mathscr P^{\antishriek}),
\end{equation}
under which $\kappa_{f^*\mathscr P}=f^*\kappa_{\mathscr P}$.
Under \eqref{eq_change_quiver_bar_cobar} and
\eqref{eq_change_quiver_koszul_dual},
\[
 \iota_{\kappa_{f^*\mathscr P}}
 =f^*(\iota_{\kappa_{\mathscr P}}),
 \qquad
 \rho_{\kappa_{f^*\mathscr P}}
 =f^*(\rho_{\kappa_{\mathscr P}}).
\]
\item If $\mathscr P$ is Koszul, then $f^*\mathscr P$ is Koszul,
and its Koszul resolution is canonically identified with the
pullback resolution:
\begin{equation}\label{eq_change_quiver_koszul_resolution}
 (f^*\mathscr P)_\infty\cong f^*(\mathscr P_\infty).
\end{equation}
\end{enumerate}
\end{thm}
\begin{proof}
(i) Write $\mathscr P=\mathscr M_{Q'}(S,R)$ and apply
Proposition~\ref{prop_change_quiver_quadratic} to $(S,R)$ and
$(s^{-1}S,s^{-2}R)$.
Changing quivers preserves the order and degrees of the two
decorations, hence commutes with
$j_2(x\circ_i y)=(-1)^{|x|}(s^{-1}x)\circ_i(s^{-1}y)$.
Consequently $f^*(s^{-2}R)=s^{-2}(f^*R)$. By the definition
\eqref{eq_koszul_dual_definition}, the cooperad comparison yields
the dg isomorphism \eqref{eq_change_quiver_koszul_dual}.
Projection to weight one, suspension, and inclusion
of generators all commute with pullback, giving the equality of
canonical twisting maps.
The cobar comparison map is determined by its values on generators,
and the bar comparison map by its corestriction to cogenerators.
These values are determined by the canonical twisting map, so
Proposition~\ref{prop_change_quiver_bar_cobar} identifies both
comparison maps with their pullbacks.

(ii) If $\mathscr P$ is Koszul, $\rho_{\kappa_{\mathscr P}}$ is a
quasi-isomorphism
in every $Q'$-profile. Its pullback is a quasi-isomorphism in every
$Q$-profile. The preceding identifications and
Proposition~\ref{thm_koszul_criterion_E_operad}(iii) prove Koszulness
of $f^*\mathscr P$ and \eqref{eq_change_quiver_koszul_resolution}.
\end{proof}

\subsection{Higher operadic homotopy algebras}
Let $\mathscr M$ be a Koszul dg operad as above, set $\mathscr C=\mathscr M^{\antishriek}$, and fix a family of cochain complexes $(X,d_X)$. An $\mathscr M_\infty$-algebra structure on $X$ is called a \textit{homotopy $\mathscr M$-algebra structure}. Write
$C(X)=\mathbb T^{\mathscr C}(X)$ for the cofree coalgebra of \S~\ref{sss_comonad_str_induced}, with coaction $\Delta_X$, projection $\varepsilon_X:C(X)\to X$, and canonical differential $d_{C(X)}$ induced by $\delta_C$ and $d_X$. Explicitly, for homogeneous $c$ and $\vec x=x_1\otimes\cdots\otimes x_n$,
\[
 d_{C(X)}(c\otimes\vec x)
 =\delta_Cc\otimes\vec x+(-1)^{|c|}c\otimes d_X^{\otimes}\vec x,
\]
where $d_X^{\otimes}$ denotes the tensor-product differential on the $X$-factors.
\Needspace{7\baselineskip}
\begin{thm}\label{thm_three_homotopy_M_algebra_definitions}
The following data are naturally equivalent:
\begin{enumerate}[(i)]
\item an $\mathscr M_\infty$-algebra structure on the fixed complex $X$;
\item a normalized twisting morphism $\theta\in\operatorname{Tw}_Q(\mathscr C,\mathbf{End}_Q(X))$;
\item a differential $d=d_{C(X)}+D$ on $C(X)$, where $D$ is a degree-one ordinary coderivation of the underlying graded $\mathscr C$-coalgebra and $D|_X=0$.
\end{enumerate}
In (iii), $d^2=0$ is required; in particular $d|_X=d_X$. The total differential $d$ is a codifferential compatible with $\delta_C$ in the sense of Definition~\ref{defn:C_coderivation}.
\end{thm}
\begin{proof}
Cobar representability \eqref{eq:fc_cobar_representability}, which does not require an augmented target, gives (i)$\Leftrightarrow$(ii). For the other equivalence, evaluation gives an isomorphism of graded vector spaces
\begin{equation*}\mathfrak g_{\mathscr C,\mathbf{End}_Q(X)}
 \cong\prod_{e\in Q}\underline{\Hom}\bigl(C(X)(e),X(e)\bigr),
 \qquad\theta\longmapsto\ev_\theta.
\end{equation*}
Here the underline denotes all homogeneous linear maps, and
$\ev_\theta(c\otimes x_1\otimes\cdots\otimes x_n)=\theta(c)(x_1,\ldots,x_n)$.
By \eqref{eq_Coder_Hom_Q}, $\ev_\theta$ extends uniquely to a coderivation $D_\theta$ with $\varepsilon_XD_\theta=\ev_\theta$.
Let $c\in C(e_1,\ldots,e_n;e_0)$ and $x_j\in X(e_j)$ be homogeneous. For a term $c_0\otimes c_1$ of $\Delta_{(1)}c$ whose inner block begins at $i$ and has length $r$, the extension is
\[
\begin{aligned}
 &D_\theta(c\otimes x_1\otimes\cdots\otimes x_n) =\sum(-1)^{|c_0|+(|c_1|+1)\sum_{j<i}|x_j|}\,
 c_0\otimes x_1\otimes\cdots\otimes x_{i-1}\\[-2pt]
 &\hspace{50mm}\otimes\ev_\theta(c_1\otimes x_i\otimes\cdots\otimes x_{i+r-1})
 \otimes x_{i+r}\otimes\cdots\otimes x_n.
\end{aligned}
\]
Coefficients of the coproduct are included in the sum. Normalization $\theta\iota_C=0$ is equivalent to $D_\theta|_X=0$.
The Hom differential, tensor rule, and displayed extension formula give
\[ \ev_{\partial\theta}=d_X\ev_\theta+\ev_\theta d_{C(X)}, \qquad \ev_{\theta\star\theta}=\ev_\theta D_\theta.
\]
Put $d=d_{C(X)}+D_\theta$. The coaction satisfies
\[
 \Delta_Xd=(\delta_C\otimes\id+\mathscr C\langle d\rangle)\Delta_X,
\]
because $d_{C(X)}$ satisfies this compatibility and $D_\theta$ is an ordinary coderivation. Hence $d^2$ is an ordinary coderivation: the terms involving both $\delta_C$ and an input differential cancel, as do the two contributions acting on distinct input factors, while $\delta_C^2=0$. This argument uses no division by $2$. Since $\varepsilon_Xd_{C(X)}=d_X\varepsilon_X$, these identities yield
\[
\varepsilon_Xd^2
=d_X\ev_\theta+\ev_\theta d_{C(X)}+\ev_\theta D_\theta
=\ev_{\partial\theta+\theta\star\theta}.
\]
A coderivation of a cofree coalgebra is determined by its corestriction, so $d^2=0$ precisely when $\theta$ satisfies the Maurer--Cartan equation. Conversely, corestriction recovers $\theta$ from any $D$ in (iii), proving bijectivity.
\end{proof}

\section{Koszul resolutions for categories and modules}
	\label{s_examples_applications}
Algebras and categories in this section are non-unital, and modules
are not required to satisfy a unit-action axiom. All operads have positive
arities; their unary identities act as identity maps, not as algebra units.
The quadratic data in the examples of this section have zero internal
differential. Thus their Koszul dual cooperads have zero differential,
although the underlying complexes and the cobar resolutions may have
nonzero differentials.

	\subsection{Koszul resolution: from dg modules to \texorpdfstring{$A_\infty$}{A-infinity}-type modules}
	\label{ss_example_left_module}

\subsubsection{A reinterpretation of ordinary modules}
Recall the left-module quiver $Q_1$ of \S~\ref{s_operads_over}: it has vertices
$v_a,v_m$, a loop $e_a:v_a\to v_a$, and an arrow $e_m:v_a\to v_m$.
The quadratic operad $\mathscr L\mathsf{Mod}=\mathscr M(S,R)$ has
degree-zero generators $\mu_a\in S(e_a,e_a;e_a)$ and
$\mu_m\in S(e_a,e_m;e_m)$, zero differential, and relations
$\mu_a\circ_1\mu_a-\mu_a\circ_2\mu_a=0$ and
$\mu_m\circ_1\mu_a-\mu_m\circ_2\mu_m=0$.

\begin{prop}
\label{prop_LMod}
An algebra over $\mathscr L\mathsf{Mod}$ is equivalent to a non-unital
dg algebra $A$ together with a left dg $A$-module $M$.
\end{prop}
\begin{proof}
Put $A=X(e_a)$ and $M=X(e_m)$. A morphism to
$\mathbf{End}_{Q_1}(X)$ assigns degree-zero cochain maps
$A\otimes A\to A$ and $A\otimes M\to M$ to the generators.
The two relations give associativity and the left module axiom, while
the cochain-map conditions give the Leibniz rules. Conversely, these
data determine the required morphism.
\end{proof}

For bimodules, let $Q$ have vertices $v_a,v_b$, loops
$e_a:v_a\to v_a$ and $e_b:v_b\to v_b$, and an arrow $e_m:v_a\to v_b$.
Let $S$ have degree-zero generators
\[
\begin{aligned}
S(e_a,e_a;e_a)&=\Bbbk\mu_a,& S(e_b,e_b;e_b)&=\Bbbk\mu_b,\\
S(e_a,e_m;e_m)&=\Bbbk\mu_l,& S(e_m,e_b;e_m)&=\Bbbk\mu_r,
\end{aligned}
\]
with zero differential and all other components zero. Let
$R\subset\mathsf F(S)^{(2)}$ be generated by
\begin{equation}
\label{eq_R_BiMod}
\begin{aligned}
r_a	&= \mu_a\circ_1\mu_a - \mu_a\circ_2\mu_a,\\
r_b	&= \mu_b\circ_1\mu_b-\mu_b\circ_2\mu_b,\\
r_l	&= \mu_l\circ_1\mu_a-\mu_l\circ_2\mu_l,\\
r_r	&= \mu_r\circ_1\mu_r-\mu_r\circ_2\mu_b,\\
r_{lr}	&= \mu_r\circ_1\mu_l- \mu_l\circ_2\mu_r.
\end{aligned}
\end{equation}

Define $\mathscr B\mathsf{Mod}:=\mathscr M(S,R)$.

\begin{prop}
\label{prop_bimod}
An algebra over $\mathscr B\mathsf{Mod}$ is equivalent to a pair of
non-unital dg algebras $A,B$ together with a dg $(A,B)$-bimodule $M$.
\end{prop}
\begin{proof}
Put $A=X(e_a)$, $B=X(e_b)$, and $M=X(e_m)$.
A morphism to $\mathbf{End}_Q(X)$ assigns two products and two actions
to the four generators. The five relations express associativity of
the products, the left and right module axioms, and compatibility of
the actions. The images are degree-zero cochain maps, giving the
Leibniz identities; conversely, these data determine the morphism.
\end{proof}

The defining relations identify all planar binary composites with a
fixed admissible profile. Denoting their common classes by
$\mu_n^a$, $\mu_n^b$, and $\mu_{p,q}^m$, respectively, we obtain
\[
\begin{aligned}
\mathscr B\mathsf{Mod}(e_a^n;e_a)&=\Bbbk\mu_n^a,&
\mathscr B\mathsf{Mod}(e_b^n;e_b)&=\Bbbk\mu_n^b,
\qquad n\ge1,\\
\mathscr B\mathsf{Mod}(e_a^p,e_m,e_b^q;e_m)&=\Bbbk\mu_{p,q}^m,
&&p,q\ge0.
\end{aligned}
\]
All other components are zero, and the unary classes are the operadic
identities. These elements encode the iterated products and the
bimodule operations $A^{\otimes p}\otimes M\otimes B^{\otimes q}\to M$.

\subsubsection{Calculation of the Koszul dual $Q$-cooperad}
\label{sss_calculate_koszul_dual_case_module}
In the common underlying $Q$-collection of $\mathsf F(s^{-1}S)$ and
$\mathsf F_c(s^{-1}S)$, set $\nu_1^a=\id_{e_a}$ and
$\nu_2^a=s^{-1}\mu_a$, and for $n\ge3$ define
\begin{equation}
\label{eq_inductive_nu_a}
\nu_n^a
:=
\sum_{\substack{p+q=n\\p,q\geq1}}
(-1)^{p-1}
\gamma
\bigl(
\nu_2^a;
\nu_p^a,\nu_q^a
\bigr),
\end{equation}
where $\gamma$ is the full composition in $\mathsf F(s^{-1}S)$.
Define $\nu_n^b$ by replacing $a$ with $b$.

For a mixed profile, mark leaf $p+1$ of a planar binary tree with
$p+q+1$ leaves, colouring the preceding leaves by $a$, the marked
leaf by $m$, and the following leaves by $b$. Every subtree has a
unique output colour: $a$ or $b$ if all its leaves have that colour,
and $m$ if it contains the marked leaf. Thus each binary vertex is
uniquely labelled by $s^{-1}\mu_a$, $s^{-1}\mu_b$,
$s^{-1}\mu_l$, or $s^{-1}\mu_r$.
Define $\nu_{p,q}^m$ by applying this colouring to each tree in
$\nu_{p+q+1}^a$, retaining its coefficient. In particular,
$\nu_{0,0}^m=\id_{e_m}$, $\nu_{1,0}^m=s^{-1}\mu_l$, and
$\nu_{0,1}^m=s^{-1}\mu_r$.

\begin{prop}
\label{prop_koszul_dual_BMod}
The nonzero components of the Koszul dual $Q$-cooperad are
\begin{equation}
\label{eq_BMod_cooperad}
\begin{aligned}
\mathscr B\mathsf{Mod}^{\antishriek}(e_a^n;e_a)
=
\Bbbk\nu_n^a,
\qquad
\mathscr B\mathsf{Mod}^{\antishriek}(e_b^n;e_b)
=
\Bbbk\nu_n^b,
\qquad n\geq1,
\\
\mathscr B\mathsf{Mod}^{\antishriek}
(e_a^p,e_m,e_b^q;e_m)
=
\Bbbk\nu_{p,q}^m,
\qquad p,q\geq0.
\end{aligned}
\end{equation}
The pure generators have degree $1-n$ and weight $n-1$; the mixed
generators have degree $-p-q$ and weight $p+q$. The full coproduct is
\begin{equation}
\label{eq_full_comultiplication_nu_a}
\Delta(\nu_n^a)
=
\sum_{\substack{
k\geq1\\
i_1+\cdots+i_k=n\\
i_1,\ldots,i_k\geq1
}}
(-1)^{	\sum_{j=1}^k
(i_j-1)(k-j) }
\,
\nu_k^a
\otimes
\nu_{i_1}^a
\otimes\cdots\otimes
\nu_{i_k}^a,
\end{equation}
and the same formula with $a$ replaced by $b$, while
\begin{equation}
\label{eq_full_comultiplication_nu_m}
\begin{aligned}
\Delta(\nu_{p,q}^m)
={}&
\sum_{\substack{
r,s\geq0,\ u,v\geq0\\
i_1+\cdots+i_r+u=p\\
v+j_1+\cdots+j_s=q\\
i_1,\ldots,i_r\geq1\\
j_1,\ldots,j_s\geq1
}}
(-1)^{\epsilon}
\,
\nu_{r,s}^m
\otimes
\nu_{i_1}^a\otimes\cdots\otimes\nu_{i_r}^a
\otimes
\nu_{u,v}^m
\otimes
\nu_{j_1}^b\otimes\cdots\otimes\nu_{j_s}^b,
\end{aligned}
\end{equation}
where
\begin{equation}
\label{eq_epsilon_comultiplication_nu_m}
\epsilon
=
\sum_{\alpha=1}^{r}
(i_\alpha-1)(r+s+1-\alpha)
+
(u+v)s
+
\sum_{\beta=1}^{s}
(j_\beta-1)(s-\beta).
\end{equation}
\end{prop}

\begin{proof}
The degrees and weights follow by counting binary vertices.
For the pure coproduct, use the associative calculation
\cite[Lemma 9.1.2]{loday2012algebraic}; we record its sign check in our
conventions. The cases $n=1,2$ are immediate. A cut of a binary-rooted
tree either lies below the root
or is a pair of cuts of its subtrees. Fix $i_1+\cdots+i_k=n$ and split
$k=r+s$, with $r,s\ge1$. Put $p=i_1+\cdots+i_r$ and
\[
E'=\sum_{j=1}^r(i_j-1)(r-j),\qquad
E''=\sum_{j=r+1}^k(i_j-1)(k-j),\qquad
E=\sum_{j=1}^k(i_j-1)(k-j).
\]
Induction gives exponents $E',E''$ for the two subtree cuts. Moving
the right outer factor past the left inner factors contributes
$(1-s)(r-p)$, so the total exponent, including the recursion sign, is
\[
p-1+E'+E''+(1-s)(r-p)\equiv r-1+E\pmod2,
\qquad E=s(p-r)+E'+E''.
\]
Since $E$ is independent of the split, summing over $r+s=k$ recovers
$\nu_k^a$ by \eqref{eq_inductive_nu_a}. The cut below the root supplies
$\nu_1^a\otimes\nu_n^a$, proving \eqref{eq_full_comultiplication_nu_a}.

The mixed colouring preserves cuts, planar vertex order, and tensor
signs. In each decomposition the unique block containing the marked
leaf has length $u+v+1$; the remaining blocks have lengths
$i_1,\ldots,i_r,j_1,\ldots,j_s$. Substitution in the pure coproduct
gives \eqref{eq_full_comultiplication_nu_m} and
\eqref{eq_epsilon_comultiplication_nu_m}, including $p=0$ or $q=0$.

The displayed spans form a subcooperad with weight-zero part
$\mathbb I_Q$, weight-one part $s^{-1}S$, and weight-two part
$s^{-2}R$, so they lie in $\mathscr C(s^{-1}S,s^{-2}R)$.
Conversely, the corelation condition at adjacent vertices fixes the
ratio of the coefficients of two binary trees related by an
associativity rotation. All binary trees with the fixed ordered leaves
are connected by these rotations, and their colours are uniquely
determined. Thus each component has dimension at most one. The
constructed elements have nonzero coefficients, proving
\eqref{eq_BMod_cooperad} and the asserted cooperad structure.
\end{proof}

\subsubsection{The explicit cobar construction}
\label{sss_bimodule_explicit_cobar}
The cobar construction $\Omega_Q\mathscr B\mathsf{Mod}^{\antishriek}$
is freely generated, as a graded $Q$-operad, by
$\mathbf m_n^\xi:=s\nu_n^\xi$ of degree $2-n$
($\xi=a,b$, $n\ge2$) and
$\mathbf n_{p,q}^m:=s\nu_{p,q}^m$ of degree $1-p-q$
($p,q\ge0$, $p+q\ge1$). The unary identities belong to the
coaugmentation and supply no generators. Since the dual cooperad has
zero internal differential, $d_\Omega=d_2$.

	\begin{lem}
		\label{lem_cobar_BMod_differential}
		The differential of
		$\Omega_Q\mathscr B\mathsf{Mod}^{\antishriek}$
		is the unique derivation determined on the generators as follows.
		For $\xi=a,b$ and $n\geq2$,
		\begin{equation*}
			d_\Omega(\mathbf m_n^\xi)
			=-
			\sum_{\substack{
					r+s+t=n\\
					s\geq2,\ r+t\geq1
			}}
			(-1)^{r+st}
			\,
			\mathbf m_{r+1+t}^\xi
			\circ_{r+1}
			\mathbf m_s^\xi .
		\end{equation*}
		For $p,q\geq0$ with $p+q\geq1$,
		\begin{align}
			\label{eq_cobar_BMod_differential_m}
			d_\Omega(\mathbf n_{p,q}^m)
			={}&-
			\sum_{\substack{
					r+s+t=p\\
					s\geq2
			}}
			(-1)^{\,r+s(t+q+1)}
			\,
			\mathbf n_{r+1+t,q}^m
			\circ_{r+1}
			\mathbf m_s^a
			\nonumber\\
			&-
			\sum_{\substack{
					r+s+t=q\\
					s\geq2
			}}
			(-1)^{\,p+r+1+st}
			\,
			\mathbf n_{p,r+1+t}^m
			\circ_{p+r+2}
			\mathbf m_s^b
			\nonumber\\
			&-
			\sum_{\substack{
					p'+p''=p,\ q'+q''=q\\
					p',q',p'',q''\geq0\\
					p'+q'\geq1,\ p''+q''\geq1
			}}
			(-1)^{
				p'+q'+(p''+q'')q'
			}
			\,
			\mathbf n_{p',q'}^m
			\circ_{p'+1}
			\mathbf n_{p'',q''}^m .
		\end{align}
	\end{lem}

	\begin{proof}
Only decompositions with one nonidentity inner factor enter the reduced
infinitesimal coproduct. Formula \eqref{eq_full_comultiplication_nu_a} gives
\begin{equation}
\label{eq:reduced_inf_comult_nu_a}
\overline\Delta_{(1)}(\nu_n^a)
=\sum_{\substack{r+s+t=n\\s\ge2,\ r+t\ge1}}
(-1)^{(s-1)t}\nu_{r+1+t}^a\otimes_{r+1}\nu_s^a.
\end{equation}
By \eqref{eq_cobar_d2_explicit}, the extra exponent is
$1+|\nu_{r+1+t}^a|=1-r-t$; hence the total exponent is
$1+r+st$ modulo $2$. The same argument applies to $b$.
For the mixed profile, \eqref{eq_full_comultiplication_nu_m} has three
possibilities: an $a$-block, a $b$-block, or the distinguished $m$-block.
Their coproduct exponents are respectively
\[
(s-1)(t+q+1),\qquad (s-1)t,\qquad (p''+q'')q'.
\]
Adding one plus the outer degree gives respectively
\[
1+r+s(t+q+1),\qquad 1+p+r+1+st,\qquad
1+p'+q'+(p''+q'')q'.
\]
These are the three negative sums in
\eqref{eq_cobar_BMod_differential_m}.
\end{proof}

	Let $A,B,M$ be cochain complexes, and denote the images of the
above generators by $m_n^a,m_n^b,n_{p,q}$. Set
$m_1^a=d_A$, $m_1^b=d_B$, and $n_{0,0}=d_M$.
The dg-operad map condition is the negative of each displayed quadratic
sum. Moving the internal differential to the same side supplies exactly
the terms with a unary differential, so the resulting identities are the
usual unshifted $A_\infty$ algebra and bimodule identities. Explicitly,
the mixed identity is the sum of the three expressions in
\eqref{eq_cobar_BMod_differential_m} without their leading minus signs,
now allowing $s=1$ and allowing either module factor to be $n_{0,0}$.
For $p=q=0$ it reads $d_M^2=0$.

For comparison with shifted conventions, put $\overline A=s^{-1}A$,
$\overline B=s^{-1}B$, and $\overline M=s^{-1}M$. For any operation
$f$ with $n$ ordered inputs, use the degree-one map on these shifted
spaces defined on homogeneous inputs by
\begin{equation}\label{eq_bimodule_suspension_convention}
 \widetilde f(s^{-1}x_1,\ldots,s^{-1}x_n)
 =(-1)^{\sum_{j=1}^n(n-j)|x_j|}\,s^{-1}f(x_1,\ldots,x_n).
\end{equation}
The unshifted relations become the square-zero coderivation equations
on the tensor coalgebras and on the tensor bicomodule with one marked
$\overline M$ input. This describes explicitly the change to the shifted
bimodule conventions; compare \cite[\S 6.5]{Yuan_higher_operad_Fuk}.

Let $\mathscr R\mathsf{Mod}$ be the right-module analogue of
$\mathscr L\mathsf{Mod}$, obtained by reversing the module arrow.
We obtain the following resolutions.

\begin{thm}
\label{thm_Koszulness_BMod}
The operads $\mathscr B\mathsf{Mod}$, $\mathscr L\mathsf{Mod}$, and
$\mathscr R\mathsf{Mod}$ are Koszul. Their cobar resolutions govern,
respectively, $A_\infty$ bimodules over pairs of $A_\infty$ algebras,
and left and right $A_\infty$ modules over an $A_\infty$ algebra.
\end{thm}
\begin{proof}
Let $\mathsf{As}$ be the positive-arity nonsymmetric associative operad,
which is Koszul \cite[\S 9.1]{loday2012algebraic}. For each of the three
quivers, let $q$ be its terminal map. The binary generators and relations
are pullbacks of multiplication and associativity; in particular,
$\mathscr B\mathsf{Mod}\cong q^*\mathsf{As}$ by \eqref{eq_R_BiMod}.
The quivers are transitive and thin, so
Example~\ref{ex_change_quiver_terminal_folding}(i) and
Theorem~\ref{thm_change_quiver_koszul}(ii) identify each canonical cobar
map with $q^*\rho_{\kappa_{\mathsf{As}}}$ and prove Koszulness.
The bimodule identities were established above. Placing the marked
input last or first gives the left or right module identities, respectively; compare
\cite[\S 6.4]{Yuan_higher_operad_Fuk}.
\end{proof}

The categorical versions are given in
Theorem~\ref{thm_Koszulness_dg_category_modules} below.

	\subsection{Koszul resolution: from dg categories to \texorpdfstring{$A_\infty$}{A-infinity} categories}
	\label{ss_example_dg_categories}

	Let $Q=(Q_0,Q,s,t)$ be a \textit{complete quiver} in the following sense: for every ordered pair of distinct vertices $(a,b)\in Q_0\times Q_0$, there is exactly
	one arrow $e_{ab}:a\to b$, and for every vertex $a\in Q_0$, there is exactly
	one loop $e_a:a\to a$. We also write $e_{aa}=e_a$.
	We refer to Remark \ref{rmk_quiver_not_complete} for the situation when $Q$ is not complete.

	Note that every admissible profile-loop in
	$Q$ is of the form
	$\pi_{\underline a}	:= (e_{a_0a_1},e_{a_1a_2},\ldots,e_{a_{n-1}a_n};e_{a_0a_n})$ for a tuple of vertices or an \textit{object string} $\underline a=(a_0,\ldots,a_n)$.

	Let $S$ be the $Q$-collection such that, for every $a,b,c\in Q_0$,
	\begin{equation}
		\label{eq_mu_abc}
		S(e_{ab},e_{bc};e_{ac})=\Bbbk\mu_{abc},
	\end{equation}
	where $\mu_{abc}$ is a degree-zero generator, and let all other components
	of $S$ be zero. We can regard $S$ as a dg $Q$-collection with zero differential.
	Let $R$ be the $Q$-subcollection of $\mathsf F(S)^{(2)}$ generated by
	\begin{equation}
		\label{eq_relations_dgCat}
		r_{abcd}
		:=
		\mu_{acd}\circ_1\mu_{abc}
		-
		\mu_{abd}\circ_2\mu_{bcd},
		\qquad a,b,c,d\in Q_0.
	\end{equation}
	Define the quadratic dg $Q$-operad
	$\mathsf{dgCat}
	:=
	\mathscr M(S,R)
	=
	\mathsf F(S)/(R)$

	\begin{prop}
		\label{prop_dgCat_algebras}
		An algebra over $\mathsf{dgCat}$ is equivalent to a non-unital dg
		category whose set of objects is $Q_0$.
	\end{prop}

	\begin{proof}
		Let $(X, \alpha)$ be a $\mathsf{dgCat}$-algebra in the sense of \S \ref{ss_algebras}. Namely,
		$X=\{X(e)\}_{e\in Q}$ is a $Q$-indexed family of cochain complexes, and	$ \alpha:	\mathsf{dgCat} \to \mathbf{End}_Q(X)$ is a morphism of dg $Q$-operads. Set
		\[
		\operatorname{Hom}(a,b):=X(e_{ab}).
		\]
		By Proposition~\ref{prop_dg_quadratic_operad_generators}, since $\mathsf{dgCat}$ is generated by the elements $\mu_{abc}$, the morphism $ \alpha$ is determined by their images $\mu_{abc}^X
		:=
		\alpha(\mu_{abc})
		\colon
		X(e_{ab})\otimes X(e_{bc})
		\longrightarrow
		X(e_{ac})$ for $a,b,c\in Q_0$.
		For homogeneous morphisms
		$f\in\operatorname{Hom}(a,b)$ and
		$g\in\operatorname{Hom}(b,c)$, write
		$f\cdot g:=\mu_{abc}^X(f,g)$.
		Since $\delta\mu_{abc}=0$ and $\alpha$ commutes with
		differentials, $\mu_{abc}^X$ is a degree-zero cochain map. Hence
		$d(f\cdot g)
		=
		(df)\cdot g
		+
		(-1)^{|f|}f\cdot(dg)$.
		Now let
		$h\in\operatorname{Hom}(c,d)$ be homogeneous. Applying
		$\alpha$ to the defining relation $r_{abcd}=0$ gives
		$
		\mu_{acd}^X\bigl(\mu_{abc}^X(f,g),h\bigr)
		=
		\mu_{abd}^X\bigl(f,\mu_{bcd}^X(g,h)\bigr)$.
		Equivalently,
		$	(f\cdot g)\cdot h=f\cdot(g\cdot h)$.
		We have therefore obtained a non-unital dg category
		with object set $Q_0$.
		Conversely, suppose that cochain complexes
		$\operatorname{Hom}(a,b)$ are equipped with composition maps
		\[
		\operatorname{Hom}(a,b)
		\otimes
		\operatorname{Hom}(b,c)
		\longrightarrow
		\operatorname{Hom}(a,c),
		\qquad
		f\otimes g\longmapsto f\cdot g,
		\]
		satisfying the displayed Leibniz rule. Set
		$X(e_{ab})=\operatorname{Hom}(a,b)$ and send
		$\mu_{abc}$ to the map $(f,g)\mapsto f\cdot g$. The Leibniz rule says
		that this image is a degree-zero cocycle in $\mathbf{End}_Q(X)$.
		The universal property of the free $Q$-operad gives a
		unique morphism of dg $Q$-operads
		$\mathsf F(S) \to \mathbf{End}_Q(X)$ (Proposition \ref{prop_universal_property_free_dg_Q_operad}).
		Associativity says precisely that every generator $r_{abcd}$ of $R$
		is sent to zero. The morphism thus factors uniquely through
		$\mathsf{dgCat}=\mathsf F(S)/(R)$, giving a
		$\mathsf{dgCat}$-algebra structure on $X$.
	\end{proof}

	\begin{rmk}
		\label{rmk_quiver_not_complete}
The same construction works for a thin quiver whose arrow relation is
\textit{transitive}: whenever $a\to b$ and $b\to c$ exist, so does $a\to c$.
Every composable input string then has all the intermediate arrows
needed for every parenthesization, so the proofs below apply unchanged.
For an arbitrary thin quiver one must also impose zero relations when
one parenthesization passes through a missing intermediate arrow.
For example, if $ab,bc,cd,ac,ad$ exist but $bd$ does not, category
associativity requires $\mu_{acd}\circ_1\mu_{abc}=0$.
We do not assert the subsequent component formulas or Koszulness
results for that more general presentation.
\end{rmk}

	For later use, we describe the underlying $Q$-collection of
	$\mathsf{dgCat}$. Set
	$\mu_{a_0a_1}:=\id_{e_{a_0a_1}}$
	and, for $n\geq2$, define recursively
	$\mu_{a_0\cdots a_n}:=\mu_{a_0a_1a_n}\circ_2\mu_{a_1\cdots a_n}$.
	The relations~\eqref{eq_relations_dgCat} imply that every composite with the same outer profile $\pi_{\underline a}$ represents the same
	element. Specifically,
	\begin{equation*}
		\mathsf{dgCat}(\pi_{\underline a}) = 	\mathsf{dgCat}(e_{a_0a_1},\dots, e_{a_{n-1}a_n}; e_{a_0a_n})
		=
		\Bbbk\mu_{a_0\cdots a_n}
	\end{equation*}
	for every object string $\underline a=(a_0, a_1,\dots, a_n)$ with $n\geq1$,
	and all other components are zero. $\mathsf{dgCat}$ is the object-coloured associative operad.

	\subsubsection{The Koszul dual $Q$-cooperad}

	We next calculate the Koszul dual $Q$-cooperad
	$\mathsf{dgCat}^{\antishriek} =\mathscr C(s^{-1}S,s^{-2}R)$.
	For a fixed object string, forgetting the uniquely determined internal colours reduces the calculation to \S\ref{sss_calculate_koszul_dual_case_module}.

	For $a_0,a_1\in Q_0$, set
	$\nu_{a_0a_1}:=\id_{e_{a_0a_1}}$ and for $a_0,a_1,a_2\in Q_0$, set $\nu_{a_0a_1a_2}:=s^{-1}\mu_{a_0a_1a_2}$.
	For $n\geq3$, define inductively
	\begin{equation}
		\label{eq_inductive_nu_dgCat}
		\nu_{a_0\cdots a_n}
		:=
		\sum_{\ell=1}^{n-1}
		(-1)^{\ell-1}
		\gamma\bigl(
		\nu_{a_0a_\ell a_n};
		\nu_{a_0\cdots a_\ell},
		\nu_{a_\ell\cdots a_n}
		\bigr),
	\end{equation}
	where $\gamma$ denotes the full composition in
	$\mathsf F(s^{-1}S)$.

	\begin{lem}
		\label{lem_comultiplication_nu_dgCat}
		For every object string $\underline a=(a_0,\ldots,a_n)$ with $n\geq1$, one has
		$|\nu_{a_0\cdots a_n}|=1-n$.
		Moreover, the full comultiplication is
		\begin{equation}
			\label{eq_full_comultiplication_nu_dgCat}
			\begin{aligned}
				\Delta(\nu_{a_0\cdots a_n})
				={}&
				\sum_{\substack{
						1\leq k\leq n\\
						0=i_0<i_1<\cdots<i_k=n
				}}
				(-1)^{\epsilon}
				\,
				\nu_{a_{i_0}a_{i_1}\cdots a_{i_k}}
				\otimes
				\nu_{a_{i_0}\cdots a_{i_1}}
				\otimes\cdots\otimes
				\nu_{a_{i_{k-1}}\cdots a_{i_k}},
			\end{aligned}
		\end{equation}
		where $\epsilon=\sum_{j=1}^{k}(i_j-i_{j-1}-1)(k-j)$.
	\end{lem}

	\begin{proof}
		After fixing the object string $(a_0,\ldots,a_n)$, the colours of all
		internal edges of every decorated planar binary diagram are uniquely
		determined. Forgetting these colours identifies the relevant diagrams
		with the usual planar binary trees occurring in the associative
		operad. Under this identification,
		\eqref{eq_inductive_nu_dgCat} is exactly the inductive formula
		\eqref{eq_inductive_nu_a}, and
		\eqref{eq_full_comultiplication_nu_dgCat} is exactly
		\eqref{eq_full_comultiplication_nu_a}. Thus the argument of
		Proposition~\ref{prop_koszul_dual_BMod} applies directly.
	\end{proof}

	\begin{lem}
		\label{lem_koszul_dual_dgCat}
		As a $Q$-collection, the Koszul dual of $\mathsf{dgCat}$ is given by
		\begin{equation*}
			\mathsf{dgCat}^{\antishriek}(\pi_{\underline a})
			=
			\Bbbk\nu_{a_0\cdots a_n},
			\qquad \underline a =(a_0,\dots, a_n)
		\end{equation*}
		for $n\geq1$, and all other components are zero. Its $Q$-cooperad structure is
		determined by
		\eqref{eq_full_comultiplication_nu_dgCat}. The element
		$\nu_{a_0\cdots a_n}$ has weight $n-1$ and degree $1-n$.
	\end{lem}

	\begin{proof}
		In weights zero and one, the assertion follows from the definition of
		$\mathscr C(s^{-1}S,s^{-2}R)$. In weight two,
		\eqref{eq_inductive_nu_dgCat} gives, up to an overall sign, the
		desuspension of the relation $r_{a_0a_1a_2a_3}$.
		For higher weights, use the recursive characterization
		\eqref{eq_quadratic_co_recursive_full}. For a fixed object string,
		forgetting the uniquely determined edge colours identifies this
		recursion with the corresponding recursion for the Koszul dual of the
		nonsymmetric associative operad. Lemma~\ref{lem_comultiplication_nu_dgCat}
		then completes the induction.
	\end{proof}

	\subsubsection{The explicit cobar construction}

	We now describe the Koszul resolution
	$\mathsf{dgCat}_\infty:=\Omega_Q\mathsf{dgCat}^{\antishriek}$.
	The elements $\nu_{a_0a_1}=\id_{e_{a_0a_1}}$ belong to the
	coaugmentation and hence do not give generators of the cobar
	construction. For every $n\geq2$, set
	$\mathbf m_{a_0\cdots a_n}:=s\nu_{a_0\cdots a_n}$.
	Then $|\mathbf m_{a_0\cdots a_n}|=2-n$
	and as a non-dg $Q$-operad,
	$\mathsf{dgCat}_\infty$ is freely generated by
	$\mathbf m_{a_0\cdots a_n}$.
	The differential is as follows.

	\begin{lem}
		\label{lem_cobar_dgCat_differential}
		The differential of $\mathsf{dgCat}_\infty$ is the unique derivation
		determined on its generators by
		\begin{equation}
			\label{eq_cobar_dgCat_differential}
			\begin{aligned}
				d_\Omega(\mathbf m_{a_0\cdots a_n})
				={}&-
				\sum_{\substack{
						r+s+t=n\\
						s\geq2,\ r+t\geq1
				}}
				(-1)^{r+st}
				\,
				\mathbf m_{a_0,\ldots,a_r,a_{r+s},\ldots,a_n}
				\circ_{r+1}
				\mathbf m_{a_r,\ldots,a_{r+s}}.
			\end{aligned}
		\end{equation}
	\end{lem}

	\begin{proof}
The unique object colouring transports the reduced infinitesimal
coproduct \eqref{eq:reduced_inf_comult_nu_a} to the object-labelled generators.
The outer degree is $-r-t$, so the sign in
\eqref{eq_cobar_d2_explicit} gives total exponent
$1-r-t+(s-1)t\equiv1+r+st\pmod2$.
\end{proof}

	\subsubsection{\texorpdfstring{$A_\infty$}{A-infinity} categories as algebras over the resolution}

	Let $X=\{X(e_{ab})\}_{a,b\in Q_0}$ be a $Q$-indexed family of
	cochain complexes, and write
	$	\mathcal A(a,b):=X(e_{ab})$ with
	$	d_{a,b}:=\mathsf d^X_{e_{ab}}$.
	An $\mathsf{dgCat}_\infty$-algebra structure on $X$ is by definition a map of dg $Q$-operads $\mathsf{dgCat}_\infty\to \mathbf{End}_Q(X)$ and is determined by the
	images of the free generators $\mathbf m_{a_0\cdots a_n}$ due to Proposition \ref{prop_universal_property_free_dg_Q_operad}. Thus, for every $n\geq2$ and every object
	string $(a_0,\ldots,a_n)$, it gives a map
	\[
	\begin{aligned}
		m_{a_0\cdots a_n}:{}
		&\mathcal A(a_0,a_1)\otimes\cdots\otimes
		\mathcal A(a_{n-1},a_n)
		\longrightarrow
		\mathcal A(a_0,a_n),
	\end{aligned}
	\]
	of degree $2-n$.
	Let $\partial$ denote the differential on the corresponding complex. By \eqref{widehat_d_eq}, this means
	$
	\partial m_{a_0\cdots a_n}
	=
	d_{a_0,a_n}\circ m_{a_0\cdots a_n}
	-
	(-1)^{2-n}
	m_{a_0\cdots a_n}\circ d_{\otimes},
	$
	where $d_{\otimes}$ is the tensor-product differential on the source tensor product.
	The condition that the structure map commute with differentials is obtained by
	applying it to \eqref{eq_cobar_dgCat_differential}, replacing each generator
	$\mathbf m$ by its image $m$ and $d_\Omega\mathbf m$ by $\partial m$.
	Here $\circ_{r+1}$ is the partial composition in the endomorphism
	$Q$-operad, so it already includes the Koszul sign from moving the inner
	map past the first $r$ inputs.

	Set $m_{a_0a_1}:=d_{a_0,a_1}$.
	Adjoining these unary differential terms, the dg-map condition becomes
	\begin{equation}
		\label{eq_Ainfty_category_Stasheff_relations}
		\sum_{\substack{
				r+s+t=n\\
				s\geq1
		}}
		(-1)^{r+st}
		\,
		m_{a_0,\ldots,a_r,a_{r+s},\ldots,a_n}
		\circ_{r+1}
		m_{a_r,\ldots,a_{r+s}}
		=0.
	\end{equation}
	When $s=1$, the inner operation is the differential $d_{a_r,a_{r+1}}$ on
	$\mathcal A(a_r,a_{r+1})$; when $r=t=0$, the outer operation is the
	differential $d_{a_0,a_n}$ on $\mathcal A(a_0,a_n)$. These are precisely the
	$A_\infty$-category relations.
	By \eqref{eq_End_X_hat_composition}, evaluating these relations on
	homogeneous elements $x_i\in\mathcal A(a_{i-1},a_i)$ adds
	$(2-s)\sum_{i=1}^{r}|x_i|$ to the exponent $r+st$ in each summand.
	Now, we further define
	$$\lambda_n(x_1,\ldots,x_n)
	:=
	\sum_{i=1}^{n}(n-i)|x_i|$$
	and set
	\begin{equation*}
			\mu^{\mathcal A}_n(x_1,\ldots,x_n)
			={}  \mu^{\mathcal A}_{a_0\cdots a_n} (x_1,\dots, x_n) 
			={}
			(-1)^{\lambda_n(x_1,\ldots,x_n)}
			m_{a_0\cdots a_n}(x_1,\ldots,x_n).
	\end{equation*}
	After this rephasing, the relations on homogeneous elements become $$
	\sum (-1)^\delta \mu_{r+t+1}^{\mathcal A} (x_1,\dots, x_r, \mu_s^{\mathcal A}(x_{r+1},\dots, x_{r+s}), \dots, x_n) =0$$
	where the sign is
	\begin{align*}
		\delta &= r+st+ (2-s)\sum_{i=1}^r  |x_i|  + \sum_{i=r+1}^{r+s}(s-(i-r))|x_i|
		\\
		&\qquad \qquad + \sum_{i=1}^r (r+t+1-i) |x_i| + t \Big(|\mu_s^{\mathcal A}|+\sum_{i=r+1}^{r+s}|x_i| \Big) + \sum_{i=r+s+1}^n (r+t+1-(i-s+1)) |x_i|  \\
		&=
		r+st + t(2-s)+\sum_{i=1}^r (1+n-i) |x_i| + \sum_{i=r+1}^{r+s} (n-i)|x_i| + \sum_{i=r+s+1}^n (n-i) |x_i| \\
		&= r + \sum_{i=1}^r |x_i| + \lambda_n(x_1,\dots, x_n)  \pmod 2
	\end{align*}
	Since the last term $\lambda_n(x_1,\dots,x_n)$ does not rely on $r,s,t$, we obtain the $A_\infty$ equation with the standard sign convention (cf. \cite{SeidelBook,FOOOBookOne}):
	\begin{equation*}
		\sum_{\substack{
				r+s+t=n\\
				s\geq1
		}}
		(-1)^{\sum_{i=1}^r(|x_i|-1)} \mu_{r+t+1}^{\mathcal A} (x_1,\dots, x_r, \mu_s^{\mathcal A}(x_{r+1},\dots, x_{r+s}), \dots, x_n) =0.
	\end{equation*}

	In summary, we have established the following result:

	\begin{thm}
		\label{thm_dgCat_resolution_AinftyCat}
		The Koszul resolution of the quadratic $Q$-operad
		$\mathsf{dgCat}$ governing non-unital dg categories with object set
		$Q_0$ is the above
		\[
		\mathsf{dgCat}_\infty
		=
		\Omega_Q\mathsf{dgCat}^{\antishriek}.
		\]
		Its algebras are precisely non-unital $A_\infty$ categories with object set $Q_0$.
	\end{thm}
\begin{proof}
The binary presentation gives $\mathsf{dgCat}\cong q^*\mathsf{As}$
for the terminal map $q:Q\to\ast$.
Example~\ref{ex_change_quiver_terminal_folding}(i) and
Theorem~\ref{thm_change_quiver_koszul}(ii) identify the displayed
cobar construction as its Koszul resolution, also for the transitive
thin quivers of Remark~\ref{rmk_quiver_not_complete}.
The algebra description follows from
\eqref{eq_Ainfty_category_Stasheff_relations}.
\end{proof}

	\subsubsection{Koszul resolution for modules over dg categories}
	\label{sss_dg_category_bimodules}
	The preceding constructions for modules over dg algebras in \S \ref{ss_example_left_module} admit a natural generalization to modules over dg categories in a similar manner as \S \ref{ss_example_dg_categories}.  We only briefly describe it below.

	Let $V$ and $W$ be two disjoint sets.  Define a quiver
	$Q^{\mathrm{bi}}_{V,W}$ by
	\[
	\left(Q^{\mathrm{bi}}_{V,W}\right)_0=V\sqcup W,
	\qquad
	Q^{\mathrm{bi}}_{V,W}
	=
	(V\times V)\sqcup(V\times W)\sqcup(W\times W).
	\]
	Namely, we have three types of arrows denoted by
	$e^A_{aa'}:a \to a'$, $e^M_{ab}:a\to  b$, $e^B_{bb'}:b\to  b'$, where $a,a'\in V$ and $b,b'\in W$.
	In particular, there are no arrows
	from $W$ to $V$.
	Similar to Remark \ref{rmk_quiver_not_complete}, we may also replace the quiver \(Q^{\mathrm{bi}}_{V,W}\) by a transitive thin subquiver. We retain the complete blocks for the formulas below.

	\begin{rmk}[Admissible profiles]\label{rmk_bimodule_admissible_profiles}
	The positive-arity profile-loops of $Q^{\mathrm{bi}}_{V,W}$ must have types
	$ (A^n;A)$, $(B^n;B)$, and $(A^p,M,B^q;M)$ for
	 $n\ge1$, $p,q\ge0$,
	where $A,B,M$ denote the three arrow types. The input arrows form
	a composable path, and the output is the unique arrow with the same
	endpoints. Insertion is allowed when the inserted output arrow equals
	the selected input arrow, including object labels; the resulting
	profile is again admissible.
	\end{rmk}

	We define a quadratic datum
	$(S^{\mathrm{bi}}_{V,W},R^{\mathrm{bi}}_{V,W})$ over
	$Q^{\mathrm{bi}}_{V,W}$; the subscripts $V,W$ record the two object sets.
	The $Q^{\mathrm{bi}}_{V,W}$-collection $S^{\mathrm{bi}}_{V,W}$ has the
	following degree-zero generators:
	\[
	\mu^A_{a_0a_1a_2}
	\in
	S^{\mathrm{bi}}_{V,W}
	(e^A_{a_0a_1},e^A_{a_1a_2};e^A_{a_0a_2}),
	\qquad
	\mu^B_{b_0b_1b_2}
	\in
	S^{\mathrm{bi}}_{V,W}
	(e^B_{b_0b_1},e^B_{b_1b_2};e^B_{b_0b_2}),
	\]
	and
	\[
	\lambda_{a_0a_1;b}
	\in
	S^{\mathrm{bi}}_{V,W}
	(e^A_{a_0a_1},e^M_{a_1b};e^M_{a_0b}),
	\qquad
	\rho_{a;b_0b_1}
	\in
	S^{\mathrm{bi}}_{V,W}
	(e^M_{ab_0},e^B_{b_0b_1};e^M_{ab_1}).
	\]
	All other components of $S^{\mathrm{bi}}_{V,W}$ are zero. The relation
	subcollection $R^{\mathrm{bi}}_{V,W}\subseteq
	\mathsf F(S^{\mathrm{bi}}_{V,W})^{(2)}$ is generated by the following
	five families:
	\begin{equation*}
		\begin{aligned}
			r^A_{a_0a_1a_2a_3}
			&=
			\mu^A_{a_0a_2a_3}\circ_1\mu^A_{a_0a_1a_2}
			-
			\mu^A_{a_0a_1a_3}\circ_2\mu^A_{a_1a_2a_3},
			\\
			r^B_{b_0b_1b_2b_3}
			&=
			\mu^B_{b_0b_2b_3}\circ_1\mu^B_{b_0b_1b_2}
			-
			\mu^B_{b_0b_1b_3}\circ_2\mu^B_{b_1b_2b_3},
			\\
			r^l_{a_0a_1a_2;b}
			&=
			\lambda_{a_0a_2;b}\circ_1\mu^A_{a_0a_1a_2}
			-
			\lambda_{a_0a_1;b}\circ_2\lambda_{a_1a_2;b},
			\\
			r^r_{a;b_0b_1b_2}
			&=
			\rho_{a;b_1b_2}\circ_1\rho_{a;b_0b_1}
			-
			\rho_{a;b_0b_2}\circ_2\mu^B_{b_0b_1b_2},
			\\
			r^{lr}_{a_0a_1;b_0b_1}
			&=
			\rho_{a_0;b_0b_1}\circ_1\lambda_{a_0a_1;b_0}
			-
			\lambda_{a_0a_1;b_1}\circ_2\rho_{a_1;b_0b_1}.
		\end{aligned}
	\end{equation*}
	Set
	\begin{equation}
		\label{eq_dg_category_bimodule_operad}
		\mathscr B\mathsf{Mod}_{V,W}
		:=
		\mathscr M(S^{\mathrm{bi}}_{V,W},R^{\mathrm{bi}}_{V,W}).
	\end{equation}

	\begin{prop}
		\label{prop_dg_category_bimodule_algebras}
		An algebra over $\mathscr B\mathsf{Mod}_{V,W}$ is equivalent to a pair of
		non-unital dg categories $\mathcal A$ and $\mathcal B$, with object sets
		$V$ and $W$, respectively, together with a dg
		$(\mathcal A,\mathcal B)$-bimodule $M$.
	\end{prop}

	By construction, such a bimodule consists
	of cochain complexes $M(a,b)$ and degree-zero cochain maps
	$\mathcal A(a_0,a_1)\otimes M(a_1,b)
	\to M(a_0,b)$ and $M(a,b_0)\otimes\mathcal B(b_0,b_1)
	\to M(a,b_1)$
	which are associative and whose left and right actions commute.
	The five relations are exactly the two category associativity axioms, the two module axioms, and their compatibility. Applying Proposition~\ref{prop_bimod} to each composable object string proves the assertion.

	There are also analogous one-sided constructions.  Let $Q^L_V$ have vertex
	set $V\sqcup\{\ast_L\}$ and arrows
	$\{e_{aa'}:a\to a'\mid a,a'\in V\}
	\sqcup
	\{\ell_a:a\to\ast_L\mid a\in V\}$,
	with no loop at $\ast_L$.  Define $\mathscr L\mathsf{Mod}_V$ using the category generators $\mu_{a_0a_1a_2}$ as \eqref{eq_mu_abc} together with generators
	$\lambda_{a_0a_1}
	\in S^L_V(e_{a_0a_1},\ell_{a_1};\ell_{a_0})$
	and the relations
	$\lambda_{a_0a_2}\circ_1\mu_{a_0a_1a_2}-\lambda_{a_0a_1}\circ_2\lambda_{a_1a_2}=0$ together with the category associativity relations.
	Similarly, let
	$Q^R_V$ have vertex set $\{\ast_R\}\sqcup V$ and arrows $e_{aa'}$ and
	$r_a:\ast_R\to a$.  The operad $\mathscr R\mathsf{Mod}_V$ is generated
	by the category generators and
	$\rho_{a_0a_1}
	\in S^R_V(r_{a_0},e_{a_0a_1};r_{a_1})$
	subject to category associativity and
	$\rho_{a_1a_2}\circ_1\rho_{a_0a_1}-\rho_{a_0a_2}\circ_2\mu_{a_0a_1a_2}=0$.
	We then set
	\[
	\mathscr L\mathsf{Mod}_V=\mathscr M(S^L_V,R^L_V),
	\qquad
	\mathscr R\mathsf{Mod}_V=\mathscr M(S^R_V,R^R_V),
	\]
	where $R^L_V$ and $R^R_V$ are generated by the above displayed module
	relations and the obvious category associativity relations.
	The same argument as above identifies their algebras with a dg category
	$\mathcal A$ with object set $V$ together with, respectively, cochain
	complexes $L(a)$ and structure maps
	$\mathcal A(a_0,a_1)\otimes L(a_1)\to L(a_0)$,
	or cochain complexes $R(a)$ and structure maps
	$R(a_0)\otimes\mathcal A(a_0,a_1)\longrightarrow R(a_1)$.
	When $V$ and $W$ are one-point sets, these three constructions recover
	$\mathscr L\mathsf{Mod}$, $\mathscr R\mathsf{Mod}$, and
	$\mathscr B\mathsf{Mod}$ from \S~\ref{ss_example_left_module}.
	The corresponding Koszul resolutions are:

	\begin{thm}
		\label{thm_Koszulness_dg_category_modules}
		The Koszul resolutions 	\begin{align*}
			(\mathscr L\mathsf{Mod}_V)_\infty
			&=
			\Omega_{Q^L_V}\mathscr L\mathsf{Mod}_V^{\antishriek},
			&
			(\mathscr R\mathsf{Mod}_V)_\infty
			&=
			\Omega_{Q^R_V}\mathscr R\mathsf{Mod}_V^{\antishriek},
			&
			(\mathscr B\mathsf{Mod}_{V,W})_\infty
			&=
			\Omega_{Q^{\mathrm{bi}}_{V,W}}
			\mathscr B\mathsf{Mod}_{V,W}^{\antishriek}.
		\end{align*}
		of the quadratic quiver-operads
		$\mathscr L\mathsf{Mod}_V$,
		$\mathscr R\mathsf{Mod}_V$, and
		$\mathscr B\mathsf{Mod}_{V,W}$ are the quiver-operads governing, respectively,
		$A_\infty$ left modules and $A_\infty$ right modules over non-unital
		$A_\infty$ categories with object set $V$, and
		$A_\infty$ bimodules over pairs of non-unital $A_\infty$ categories with
		object sets $V$ and $W$.
	\end{thm}
\begin{proof}
The displayed quadratic data are pullbacks of the associative datum
along the terminal maps of their transitive thin quivers.
Example~\ref{ex_change_quiver_terminal_folding}(i) and
Theorem~\ref{thm_change_quiver_koszul}(ii) therefore give the three
Koszul resolutions.
Their algebra identities are the object-labelled versions of the
identities already proved.
\end{proof}

\subsection{Derived \texorpdfstring{$A_\infty$}{A-infinity} categories and bimodules}
\label{s_derived_Ainfinity_categories}

We extend Sagave's derived $A_\infty$ algebras \cite{sagave2007dg}
and their operadic description in
\cite{Livernet_Roitzheim_Whitehouse_derived_Ainfinity} to categories and
bimodules. The computations follow the derived associative calculation,
with object labels and a distinguished module input as in the preceding
examples; we record the resulting structures and omit the parallel
calculations.

\subsubsection{Conventions and bidg-categories}
\label{sss_bidg_categories}
A \textit{vertical bicomplex} is an $(\mathbb N,\mathbb Z)$-bigraded
vector space $A=\bigoplus_{p\geq0,q\in\mathbb Z}A_p^q$ with a differential
$d^{\mathrm v}$ of bidegree $(0,1)$. A map of bidegree $(i,j)$ lowers
horizontal degree by $i$ and raises vertical degree by $j$. We use the
sign convention
$(f\otimes g)(x\otimes y)=(-1)^{i_gp_x+j_gq_x}f(x)\otimes g(y)$,
where $g$ has bidegree $(i_g,j_g)$ and $x\in A_{p_x}^{q_x}$.
Throughout, negative horizontal components are zero, and summands with
a negative intermediate colour are omitted.

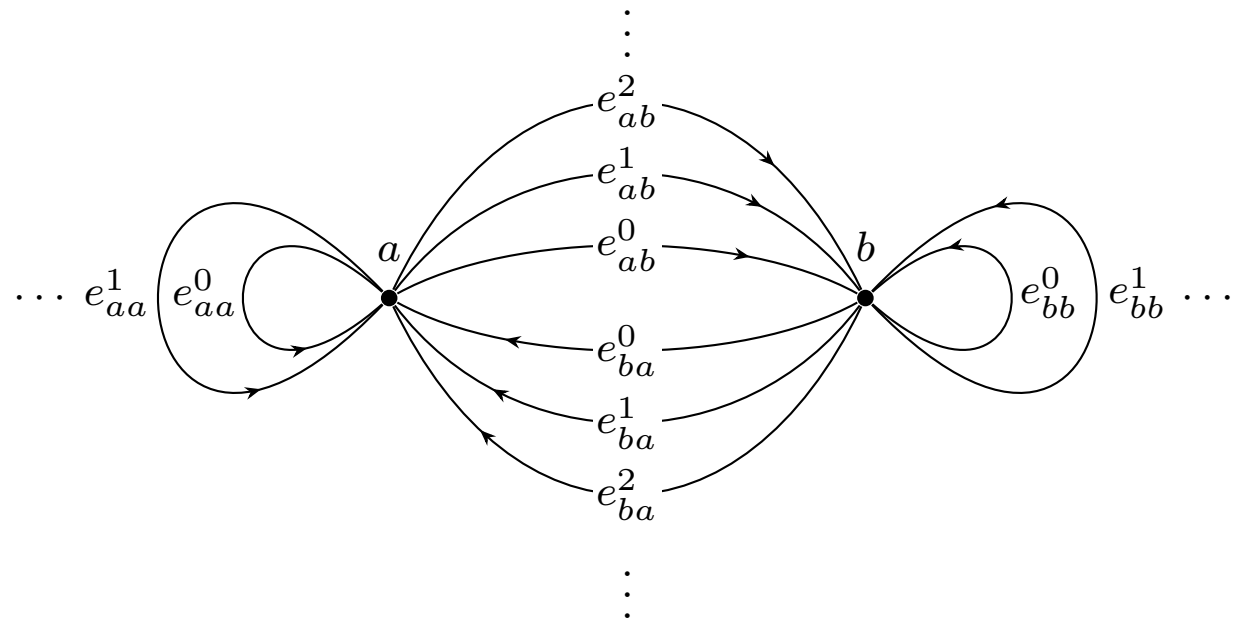
\begin{figure}[htbp]
\centering
\begin{tikzpicture}[
  x=.75cm,y=.7cm,font=\scriptsize,line width=.4pt,
  arr/.style={postaction={decorate},decoration={markings,
    mark=at position .76 with {\arrow{Stealth[length=3pt,width=2.4pt]}}}},
  vertex/.style={circle,fill=black,inner sep=1.1pt},
  lab/.style={fill=white,inner sep=.6pt}
]
\node[vertex,label={[label distance=2pt]90:$a$}] (a) at (-2.1,0) {};
\node[vertex,label={[label distance=2pt]90:$b$}] (b) at (2.1,0) {};

\draw[arr] (a) .. controls (-1.0,.65) and (1.0,.65) ..
  node[lab,pos=.5] {$e^0_{ab}$} (b);
\draw[arr] (a) .. controls (-1.0,1.55) and (1.0,1.55) ..
  node[lab,pos=.5] {$e^1_{ab}$} (b);
\draw[arr] (a) .. controls (-1.0,2.45) and (1.0,2.45) ..
  node[lab,pos=.5] {$e^2_{ab}$} (b);
\node at (0,2.65) {$\vdots$};

\draw[arr] (b) .. controls (1.0,-.65) and (-1.0,-.65) ..
  node[lab,pos=.5] {$e^0_{ba}$} (a);
\draw[arr] (b) .. controls (1.0,-1.55) and (-1.0,-1.55) ..
  node[lab,pos=.5] {$e^1_{ba}$} (a);
\draw[arr] (b) .. controls (1.0,-2.45) and (-1.0,-2.45) ..
  node[lab,pos=.5] {$e^2_{ba}$} (a);
\node at (0,-2.65) {$\vdots$};

\draw[arr] (a) .. controls (-3.8,1.6) and (-3.8,-1.6) .. (a);
\draw[arr] (a) .. controls (-4.8,3.0) and (-4.8,-3.0) .. (a);
\node at (-3.72,.05) {$e^0_{aa}$};
\node at (-4.5,.05) {$e^1_{aa}$};
\node at (-5.15,0) {$\cdots$};
\draw[arr] (b) .. controls (3.8,-1.6) and (3.8,1.6) .. (b);
\draw[arr] (b) .. controls (4.8,-3.0) and (4.8,3.0) .. (b);
\node at (3.72,.05) {$e^0_{bb}$};
\node at (4.5,.05) {$e^1_{bb}$};
\node at (5.15,0) {$\cdots$};
\pgfresetboundingbox
\path[use as bounding box] (-5.45,-2.95) rectangle (5.45,2.95);
\end{tikzpicture}
\caption{A two-vertex portion of $Q$, with arrows and loops indexed by
horizontal degree $p\geq0$.}
\label{fig:derived-category-quiver}
\end{figure}

Let $Q$ have vertices $V$ and arrows $e^p_{ab}:a\to b$ for
$a,b\in V$, $p\geq0$; see Figure~\ref{fig:derived-category-quiver}.

\Needspace{4\baselineskip}
The grading of the $Q$-collection records vertical degree; horizontal
degrees are encoded by the arrow labels.
Define $\mathsf{bidgCat}_V=\mathscr M(S,R)$ by degree-zero generators
$\lambda^p_{ab}\in S(e^p_{ab};e^{p-1}_{ab})$ for $p\geq1$ and
$\mu^{p_1,p_2}_{abc}\in
S(e^{p_1}_{ab},e^{p_2}_{bc};e^{p_1+p_2}_{ac})$ for $p_1,p_2\geq0$,
subject to the relations, with object labels suppressed,
\[
\begin{aligned}
 &\lambda^{p-1}\circ_1\lambda^p=0,\qquad p\geq2,\\
 &\mu^{p_1+p_2,p_3}\circ_1\mu^{p_1,p_2}
   =\mu^{p_1,p_2+p_3}\circ_2\mu^{p_2,p_3},\\
 &\lambda^{p_1+p_2}\circ_1\mu^{p_1,p_2}
   =\mu^{p_1-1,p_2}\circ_1\lambda^{p_1}
    +(-1)^{p_1}\mu^{p_1,p_2-1}\circ_2\lambda^{p_2}.
\end{aligned}
\]
For an algebra $(X,\alpha)$ over $\mathsf{bidgCat}_V$, put
$X_p(a,b)=X(e^p_{ab})$. The internal differential of each cochain
complex supplies $d^{\mathrm v}$, while the unary generator supplies
$d^{\mathrm h}$:
\[
\begin{aligned}
 d^{\mathrm v}&=d_{X_p(a,b)}:
 X_p(a,b)\longrightarrow X_p(a,b),\\
 d^{\mathrm h}&=\alpha(\lambda^p_{ab}):
 X_p(a,b)\longrightarrow X_{p-1}(a,b)\qquad(p\geq1),
\end{aligned}
\]
with $d^{\mathrm h}=0$ on $X_0(a,b)$. The binary generators give
composition via $\alpha(\mu^{p_1,p_2}_{abc})$. The defining relations
give $(d^{\mathrm h})^2=0$, associativity, and the horizontal Leibniz
rule; since $\alpha$ is a dg-operad map, $d^{\mathrm h}$ commutes with
$d^{\mathrm v}$ and composition satisfies the vertical Leibniz rule.
We call these algebras \textit{non-unital bidg-categories}. For one
object this recovers a \textit{bidga}.

\subsubsection{Derived categories and their Koszul resolution}
A \textit{non-unital derived $A_\infty$ category} consists of vertical
bicomplexes $\mathcal A(a,b)$ and operations
\[
 \mu_{u,n}^{a_0\cdots a_n}:
 \bigotimes_{j=1}^{n}\mathcal A(a_{j-1},a_j)
 \longrightarrow\mathcal A(a_0,a_n),
 \qquad u\geq0,\ n\geq1,
\]
of bidegree $(u,2-u-n)$, with $\mu_{0,1}=d^{\mathrm v}$, satisfying
\begin{equation}\label{eq_derived_Ainfinity_category_relations}
 \sum_{\substack{u=i+v,\ n=r+q+t\\i,v,r,t\geq0,\ q\geq1}}
 (-1)^{r+qt+v(r+1+t)}
 \mu_{i,r+1+t}
 (\id^{\otimes r}\otimes\mu_{v,q}\otimes\id^{\otimes t})=0.
\end{equation}
For one object, this recovers Sagave's derived $A_\infty$ algebras
\cite{sagave2007dg} with our convention. The sign conversion
is $\mu_{u,n}=(-1)^{\binom n2+u+1}\overline\mu_{u,n}$, where
$\overline\mu$ uses the convention of
\cite[\S3.1]{Livernet_Roitzheim_Whitehouse_derived_Ainfinity}.
Besides, concentration in horizontal degree zero forces $\mu_{u,n}=0$ for
$u>0$, since the output degree would be negative. The remaining
$\mu_{0,n}$ have degree $2-n$ and satisfy the ordinary $A_\infty$
category identities \eqref{eq_Ainfty_category_Stasheff_relations}.

\begin{thm}\label{thm_bidgCat_resolution_derived_Ainfinity_categories}
The operad $\mathsf{bidgCat}_V$ is Koszul. Algebras over its Koszul
resolution $(\mathsf{bidgCat}_V)_\infty
=\Omega_{Q}(\mathsf{bidgCat}_V^{\antishriek})$ are precisely
non-unital derived $A_\infty$ categories with object set $V$.
\end{thm}
\begin{proof}[Sketch of proof]
The one-object operad $\mathsf{bidgCat}_{\{*\}}$ is essentially the bidga operad $\mathsf{dAs}$ and is Koszul by
\cite[Proposition~3.3]{Livernet_Roitzheim_Whitehouse_derived_Ainfinity}.
Let $q$ be the map of Example~\ref{ex_change_quiver_terminal_folding}(iii)
for the complete quiver on $V$. The presentations give
$\mathsf{bidgCat}_V\cong q^*\mathsf{bidgCat}_{\{*\}}$, so
Theorem~\ref{thm_change_quiver_koszul}(ii) gives its Koszul resolution.
The formulas below identify the resolution algebras with
\eqref{eq_derived_Ainfinity_category_relations}.
\end{proof}

For later use, we record the dual data in the normalization of
\cite[Theorem~2.8 and Corollary~2.9]
{Livernet_Roitzheim_Whitehouse_derived_Ainfinity}.
Write $\nu_{u,n}=\nu_{u,n;a_0\cdots a_n}^{p_1,\ldots,p_n}$ and
$p_0=\sum_jp_j-u\geq0$. The dual has zero differential and components
\[
 \mathsf{bidgCat}_V^{\antishriek}
 (e^{p_1}_{a_0a_1},\ldots,e^{p_n}_{a_{n-1}a_n};e^{p_0}_{a_0a_n})
 =\Bbbk\nu_{u,n},\qquad |\nu_{u,n}|=1-u-n\quad(u\geq0,n\geq1),
\]
all other components being zero. Its underlying collection is
$D^{\antishriek}\circ_{Q}P^{\antishriek}$, where $P$ is the $Q$-operad
generated by the multiplications $\mu^{p_1,p_2}_{abc}$ subject to
associativity, and $D$ is generated by the unary operations
$\lambda^p_{ab}$ subject to the square-zero relation.
The basis element is
represented by an associative dual operation followed by $u$ unary
dual operations. In particular, $\nu_{0,1}=\id$,
$\nu_{1,1}=s^{-1}\lambda$, and $\nu_{0,2}=s^{-1}\mu$.
Suppressing the uniquely determined object and horizontal labels,
\begin{equation}\label{eq_bidgCat_dual_infinitesimal_comultiplication}
 \Delta_{(1)}\nu_{u,n}
 =\sum_{\substack{u=i+v,\ n=r+q+t\\i,v,r,t\geq0,\ q\geq1}}
 (-1)^{r(1-v-q)+vt+v\sum_{\ell=1}^{r}p_\ell}
 \nu_{i,r+1+t}\otimes_{r+1}\nu_{v,q}.
\end{equation}
The inner output has colour $\sum_{\ell=r+1}^{r+q}p_\ell-v$ and
endpoints $(a_r,a_{r+q})$. The extra horizontal sign is exactly that
of the tensor convention in \S~\ref{sss_bidg_categories}.
The cobar generators are
$\mathbf m_{u,n}=(-1)^{\binom n2+n-1}s\nu_{u,n}$ of degree $2-u-n$,
excluding $(u,n)=(0,1)$. Their differential follows from
\eqref{eq_bidgCat_dual_infinitesimal_comultiplication} and
\eqref{eq_cobar_d2_explicit}; adjoining $\mu_{0,1}=d^{\mathrm v}$ gives
\eqref{eq_derived_Ainfinity_category_relations}, whose differential
terms combine to $(-1)^u\partial\mu_{u,n}$. The canonical resolution map sends
$\mathbf m_{1,1}$ to $\lambda$ and $\mathbf m_{0,2}$ to $\mu$, and
all other generators to zero.

\subsubsection{Bimodules}
\label{sss_derived_Ainfinity_category_bimodules}
Let $\mathcal A,\mathcal B$ be derived $A_\infty$ categories on disjoint
object sets $V,W$. A \textit{derived $A_\infty$
$(\mathcal A,\mathcal B)$-bimodule} consists of vertical bicomplexes
$\mathcal M(a,b)$ and operations
\[
\begin{aligned}
 \mu^{\mathcal M}_{u,k,\ell}:{}
 &\bigotimes_{i=1}^{k}\mathcal A(a_{i-1},a_i)
 \otimes\mathcal M(a_k,b_0)
 \otimes\bigotimes_{j=1}^{\ell}\mathcal B(b_{j-1},b_j)
 \longrightarrow\mathcal M(a_0,b_\ell),\\
 &u,k,\ell\geq0,\qquad \mu^{\mathcal M}_{0,0,0}=d^{\mathrm v}_{\mathcal M},
\end{aligned}
\]
of bidegree $(u,1-u-k-\ell)$. They satisfy
\eqref{eq_derived_Ainfinity_category_relations} on strings
with one distinguished module input, taking $n=k+\ell+1$.
An inner block before that input uses an $\mathcal A$-operation, one
after it uses a $\mathcal B$-operation, and one containing it uses an
$\mathcal M$-operation; the outer operation is the corresponding module
operation. The signs and the negative-colour convention are unchanged.

Let $Q^{\mathrm{bi}}_{V,W}$ have arrows $e^{A,p}_{aa'}$, $e^{M,p}_{ab}$,
and $e^{B,p}_{bb'}$ for $p\geq0$. Define
$\mathsf{bidgBMod}_{V,W}$ by refining the generators and five relation
families of \eqref{eq_dg_category_bimodule_operad} by additive horizontal
colours, adjoining degree-zero unary generators lowering each arrow's
colour by one, and imposing their square-zero relations and the horizontal Leibniz
relation above for each product and action. Its algebras are pairs of
non-unital bidg-categories with a bidg-bimodule.

\begin{thm}\label{thm_derived_Ainfinity_category_bimodules}
The operad $\mathsf{bidgBMod}_{V,W}$ is Koszul. Algebras over its
Koszul resolution
$\Omega_{Q^{\mathrm{bi}}_{V,W}}
(\mathsf{bidgBMod}_{V,W}^{\antishriek})$
are precisely pairs of non-unital derived $A_\infty$ categories on
$V,W$ together with a derived $A_\infty$ bimodule.
\end{thm}
\begin{proof}
Let $q$ be the map of Example~\ref{ex_change_quiver_terminal_folding}(iii)
for the triangular quiver on $V\sqcup W$. The presentations give
$\mathsf{bidgBMod}_{V,W}\cong q^*\mathsf{bidgCat}_{\{*\}}$.
Theorems~\ref{thm_bidgCat_resolution_derived_Ainfinity_categories}
and~\ref{thm_change_quiver_koszul}(ii) give its Koszul resolution.
The pulled-back cobar formulas on mixed profiles, with $n=k+\ell+1$,
give generators of degree $1-u-k-\ell$ and the stated module identities.
\end{proof}
Singleton object sets recover bimodules over pairs of derived
$A_\infty$ algebras; when both algebra structures coincide, this is
the two-sided representation of
\cite[Theorem~4.20 and Remark~4.21]
{Aponte_Livernet_Robertson_Whitehouse_Ziegenhagen_representations}.
Concentration in horizontal degree zero recovers ordinary $A_\infty$
bimodules.

\section{Change of quiver: diagonal bimodules and tensor products}
\label{sss_diagonal_bimodule_pullback}
We specialize the restriction and extension functors of
Propositions~\ref{prop_change_quiver_algebras}--\ref{prop_change_quiver_algebra_adjunction}
to diagonal bimodules and tensor products, using the compatibility
of Koszul resolutions with pullback in
Theorem~\ref{thm_change_quiver_koszul}.
All Hom sets below consist of strict operadic algebra morphisms
as in Definition~\ref{defn_quiver_operad_algebra_morphisms}.

\subsection{Diagonal $A_\infty$ bimodules by quiver pullback}
\label{sss_diagonal_bimodule_construction}
\begin{figure}[htbp]
\centering
\resizebox{\linewidth}{!}{\begingroup
\definecolor{leftcat}{HTML}{285E9A}
\definecolor{rightcat}{HTML}{A66A19}
\definecolor{mixed}{HTML}{85529B}
\begin{minipage}{17.5cm}
\centering
\begin{tikzpicture}[
  x=1cm,y=1cm,font=\small,line width=.65pt,
  arr/.style={postaction={decorate},decoration={markings,
    mark=at position .77 with {\arrow{Stealth[length=4.2pt,width=3.4pt]}}}},
  A/.style={arr,draw=leftcat},
  B/.style={arr,draw=rightcat},
  M/.style={arr,draw=mixed},
  vertex/.style={circle,fill=black,inner sep=1.8pt},
  lab/.style={fill=white,inner sep=1.8pt,font=\scriptsize}
]

% Dashed frames identify the two blocks of internal arrows.
\draw[draw=leftcat!70,dashed,rounded corners=5pt,line width=.55pt]
  (-3.65,.85) rectangle (3.65,2.7);
\node[fill=white,inner xsep=5pt,inner ysep=2pt,text=leftcat]
  at (-2.55,2.7) {$V_A\times V_A$};
\draw[draw=rightcat!70,dashed,rounded corners=5pt,line width=.55pt]
  (-3.65,-2.7) rectangle (3.65,-.85);
\node[fill=white,inner xsep=5pt,inner ysep=2pt,text=rightcat]
  at (-2.55,-2.7) {$V_B\times V_B$};

% The two complete blocks, including their loops.
\node[vertex,label={[label distance=3pt]90:$a_A$}] (aA) at (-2,1.35) {};
\node[vertex,label={[label distance=3pt]90:$b_A$}] (bA) at (2,1.35) {};
\node[vertex,label={[label distance=3pt]270:$a_B$}] (aB) at (-2,-1.35) {};
\node[vertex,label={[label distance=3pt]270:$b_B$}] (bB) at (2,-1.35) {};

\draw[A] (aA) .. controls (-.9,2.75) and (.9,2.75) ..
  node[lab,pos=.5,text=leftcat] {$e^A_{ab}$} (bA);
\draw[A] (bA) .. controls (.9,1.95) and (-.9,1.95) ..
  node[lab,pos=.5,text=leftcat] {$e^A_{ba}$} (aA);
\draw[A] (aA) .. controls (-3.5,2.35) and (-3.5,.35) ..
  node[lab,pos=.5,text=leftcat] {$e^A_{aa}$} (aA);
\draw[A] (bA) .. controls (3.5,.35) and (3.5,2.35) ..
  node[lab,pos=.5,text=leftcat] {$e^A_{bb}$} (bA);

\draw[B] (aB) .. controls (-.9,-1.95) and (.9,-1.95) ..
  node[lab,pos=.5,text=rightcat] {$e^B_{ab}$} (bB);
\draw[B] (bB) .. controls (.9,-2.75) and (-.9,-2.75) ..
  node[lab,pos=.5,text=rightcat] {$e^B_{ba}$} (aB);
\draw[B] (aB) .. controls (-3.5,-.35) and (-3.5,-2.35) ..
  node[lab,pos=.5,text=rightcat] {$e^B_{aa}$} (aB);
\draw[B] (bB) .. controls (3.5,-2.35) and (3.5,-.35) ..
  node[lab,pos=.5,text=rightcat] {$e^B_{bb}$} (bB);

% All mixed arrows go from the upper (left-copy) block to the lower block.
\draw[M] (aA) -- node[lab,pos=.5,left=2pt,text=mixed] {$e^M_{aa}$} (aB);
\draw[M] (bA) -- node[lab,pos=.5,right=2pt,text=mixed] {$e^M_{bb}$} (bB);
\draw[M] (aA) -- node[lab,pos=.27,text=mixed] {$e^M_{ab}$} (bB);
\draw[M] (bA) -- node[lab,pos=.27,text=mixed] {$e^M_{ba}$} (aB);

% Quiver morphism. The large arrow is not an edge of either quiver.
\draw[-{Stealth[length=7pt,width=5pt]},line width=.8pt]
  (4.15,0) -- node[above=5pt,font=\large] {$f$} (5.95,0);

% Complete target quiver.
\node[vertex,label={[label distance=5pt]90:$a$}] (a) at (7.6,0) {};
\node[vertex,label={[label distance=5pt]90:$b$}] (b) at (10.8,0) {};
\draw[arr] (a) .. controls (8.4,1.05) and (10,1.05) ..
  node[lab,pos=.5] {$e_{ab}$} (b);
\draw[arr] (b) .. controls (10,-1.05) and (8.4,-1.05) ..
  node[lab,pos=.5] {$e_{ba}$} (a);
\draw[arr] (a) .. controls (6.3,.95) and (6.3,-.95) ..
  node[lab,pos=.5] {$e_{aa}$} (a);
\draw[arr] (b) .. controls (12.1,-.95) and (12.1,.95) ..
  node[lab,pos=.5] {$e_{bb}$} (b);

\node[font=\small] at (9.2,-1.7) {$a_A,a_B\longmapsto a$};
\node[font=\small] at (9.2,-2.15) {$b_A,b_B\longmapsto b$};

\path[use as bounding box] (-3.85,-3.05) rectangle (12.25,3.05);
\end{tikzpicture}

\smallskip
{\small
\textcolor{leftcat}{$A$: left category}\qquad
\textcolor{mixed}{$M$: bimodule}\qquad
\textcolor{rightcat}{$B$: right category}
\par}
\end{minipage}
\endgroup}
\caption{The folding map $f:Q^{\mathrm{bi}}\to Q_V$ for $V=\{a,b\}$.
The dashed frames enclose the two blocks of internal arrows.
Each target arrow has one preimage of each type $A$, $M$, and $B$.
Compare the bimodule-quiver diagram in
\cite[\S~6.6, p.~45]{Yuan_higher_operad_Fuk}.}
\label{fig:diagonal-quiver-folding}
\end{figure}
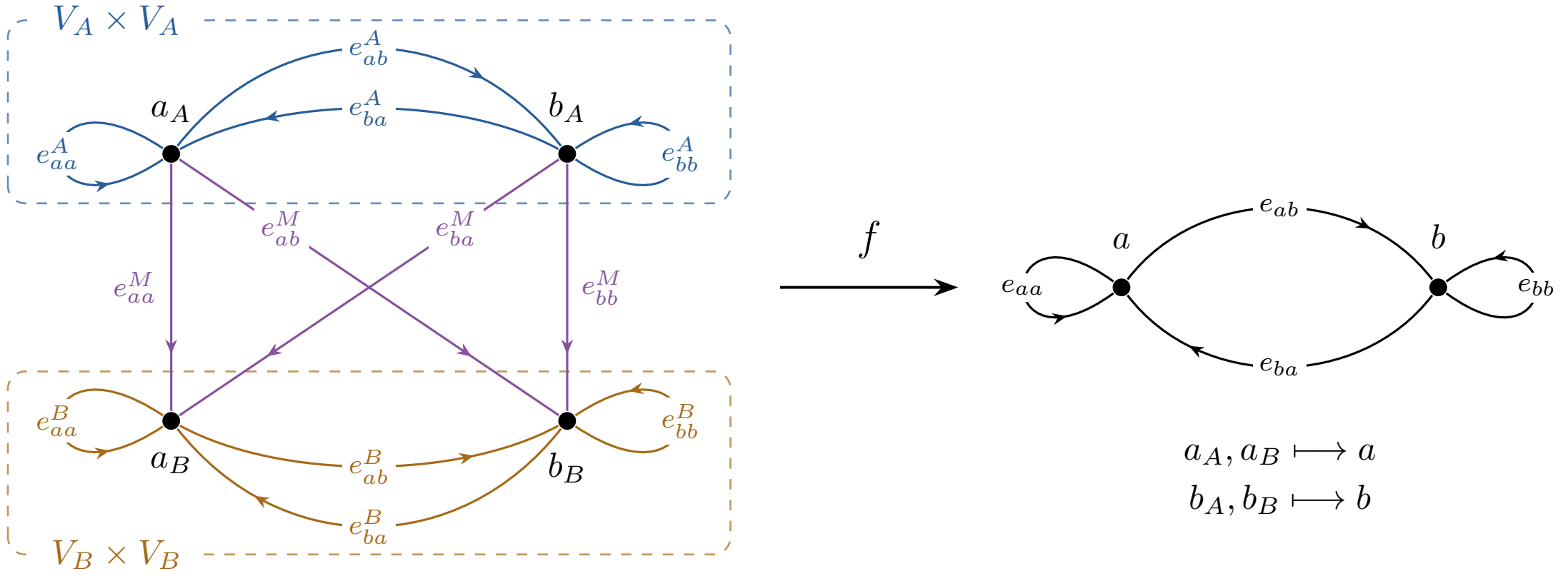

Write $Q_V$ for the complete quiver on a set $V$ and
$\mathsf{dgCat}_V$ for its category operad. Let $V_A,V_B$ be
disjoint copies of $V$, with $a_A\in V_A$ and $a_B\in V_B$
denoting the two copies of each $a\in V$.
Let $Q^{\mathrm{bi}}=Q^{\mathrm{bi}}_{V_A,V_B}$ be the quiver
defined as in \S~\ref{sss_dg_category_bimodules}: its vertex
set is $V_A\sqcup V_B$ and its arrow set is
$(V_A\times V_A)\sqcup(V_A\times V_B)\sqcup(V_B\times V_B)$.
For $a,b\in V$, write
$ e^A_{ab}:a_A\to b_A$, $
 e^M_{ab}:a_A\to b_B$, and $
 e^B_{ab}:a_B\to b_B$.
Define the folding map (cf. Figure~\ref{fig:diagonal-quiver-folding})
\begin{equation}\label{eq_diagonal_quiver_folding}
\begin{gathered}
 f:Q^{\mathrm{bi}}\longrightarrow Q_V,
 \qquad f(a_A)=f(a_B)=a, \qquad 
 f(e^A_{ab})=f(e^M_{ab})=f(e^B_{ab})=e_{ab}.
\end{gathered}
\end{equation}

\begin{prop}\label{prop_diagonal_bimodule_pullback}
The map $f$ has unique diagram lifting. There are canonical
isomorphisms of dg $Q^{\mathrm{bi}}$-operads
\begin{equation}\label{eq_diagonal_operad_pullback}
\begin{aligned}
 \mathscr B\mathsf{Mod}_{V_A,V_B}
 &\cong f^*\mathsf{dgCat}_V,\\
 (\mathscr B\mathsf{Mod}_{V_A,V_B})_\infty
 &\cong f^*((\mathsf{dgCat}_V)_\infty),
\end{aligned}
\end{equation}
and the second isomorphism identifies their canonical resolution maps
to the operads in the first line.
\end{prop}
\begin{proof}
By Example~\ref{ex_change_quiver_terminal_folding}(ii), $f$ has unique
diagram lifting. Write $\mathscr P=\mathsf{dgCat}_V
=\mathscr M_{Q_V}(S,R)$. Since $f$ forgets only copy labels,
the generator \eqref{eq_mu_abc} and relation \eqref{eq_relations_dgCat}
pull back to the following families of
\S~\ref{sss_dg_category_bimodules}, with object indices suppressed:
\[
\begin{array}{c|c|c}
 & \text{input types} & \text{corresponding families} \\ \hline
 f^*S & AA,\ BB,\ AM,\ MB
      & \mu^A,\ \mu^B,\ \lambda,\ \rho \\
 f^*R & AAA,\ BBB,\ AAM,\ MBB,\ AMB
      & r^A,\ r^B,\ r^l,\ r^r,\ r^{lr}
\end{array}
\]
These are precisely the generator and relation collections
$S^{\mathrm{bi}}_{V_A,V_B}$ and $R^{\mathrm{bi}}_{V_A,V_B}$ from
\S~\ref{sss_dg_category_bimodules}, with $(V,W)=(V_A,V_B)$.
Hence, $
 f^*S=S^{\mathrm{bi}}_{V_A,V_B}$ and $
 f^*R=R^{\mathrm{bi}}_{V_A,V_B}$ under these identifications.
 Therefore,
Proposition~\ref{prop_change_quiver_quadratic} gives
$\mathscr B\mathsf{Mod}_{V_A,V_B}\cong f^*\mathscr P$.
Since $\mathscr P$ is Koszul, Theorem~\ref{thm_change_quiver_koszul}
yields
$ (f^*\mathscr P)_\infty\cong f^*(\mathscr P_\infty)$ and $
 \rho_{\kappa_{f^*\mathscr P}}
 =f^*(\rho_{\kappa_{\mathscr P}})$.
The first shows that forming the Koszul resolution commutes with
pullback, giving the second isomorphism in
\eqref{eq_diagonal_operad_pullback}. The equality of the $\rho$-maps
expresses its compatibility with the canonical resolution maps.
\end{proof}

For the resolution generators, suppress the copy labels in object
indices. Write $\mathbf n_{a_0\cdots a_p;b_0\cdots b_q}$ for the
mixed generator whose inputs consist of $p$ left-category arrows,
the arrow from $a_p$ to $b_0$, and $q$ right-category arrows.
The second isomorphism in \eqref{eq_diagonal_operad_pullback} is
given by
\begin{equation}\label{eq_diagonal_generator_folding}
\begin{aligned}
 \mathbf m^A_{a_0\cdots a_n}&\longmapsto
 \mathbf m_{a_0\cdots a_n},
 &\mathbf m^B_{b_0\cdots b_n}&\longmapsto
 \mathbf m_{b_0\cdots b_n},\\
 \mathbf n_{a_0\cdots a_p;b_0\cdots b_q}
 &\longmapsto \mathbf m_{a_0\cdots a_p b_0\cdots b_q},
\end{aligned}
\end{equation}
where $n\geq2$ in the pure cases and $p+q\geq1$ in the mixed case.
These assignments retain the ordered vertex labels in the
decorated-tree construction, so they preserve the differentials
with the unshifted signs already specified above.

\begin{cor}\label{cor_diagonal_Ainfinity_bimodule}
Let $\mathcal A$ be an $A_\infty$ category on $V$, with operations
$m^{\mathcal A}_{a_0\cdots a_n}$ in the convention of
\eqref{eq_Ainfty_category_Stasheff_relations}, where
$m^{\mathcal A}_{ab}=d_{a,b}$ is the differential on
$\mathcal A(a,b)$. Pullback along $f$ gives its two copies
$\mathcal A_A,\mathcal A_B$ and their diagonal $A_\infty$ bimodule
$\mathcal A_\Delta$, with $\mathcal A_\Delta(a_A,b_B)=\mathcal A(a,b)$.
For $p,q\geq0$, its operations, with copy labels suppressed, are
\begin{equation}\label{eq_diagonal_bimodule_operations}
 n^\Delta_{a_0\cdots a_p;b_0\cdots b_q}
 =m^{\mathcal A}_{a_0\cdots a_p b_0\cdots b_q}: \quad 
 \bigotimes_{i=1}^{p}\mathcal A(a_{i-1},a_i)
 \otimes\mathcal A_\Delta(a_p,b_0)
 \otimes\bigotimes_{j=1}^{q}\mathcal A(b_{j-1},b_j)
 \longrightarrow\mathcal A_\Delta(a_0,b_q), \quad 
\end{equation}
of degree $1-p-q$. In particular, $n^\Delta_{a;b}=d_{a,b}$.
\end{cor}
\begin{proof}
Let $\alpha$ be the structure map on $X(e_{ab})=\mathcal A(a,b)$.
By \eqref{eq_diagonal_operad_pullback} and
Proposition~\ref{prop_change_quiver_algebras}, it pulls back to
\[
 (\mathscr B\mathsf{Mod}_{V_A,V_B})_\infty
 \xrightarrow{\,f^*\alpha\,} f^*\mathbf{End}_{Q_V}(X)
 \cong\mathbf{End}_{Q^{\mathrm{bi}}}(f^*X).
\]
The generator assignments \eqref{eq_diagonal_generator_folding}
give the two copied categories and the stated bimodule operations;
the case $p=q=0$ is the differential of $f^*X$.
\end{proof}

The second isomorphism in \eqref{eq_diagonal_operad_pullback}
defines a morphism of dg quiver-operads
\[
 \Phi_f:
 \bigl(Q^{\mathrm{bi}},(\mathscr B\mathsf{Mod}_{V_A,V_B})_\infty\bigr)
 \longrightarrow
 \bigl(Q_V,(\mathsf{dgCat}_V)_\infty\bigr)
\]
over the folding map $f$.
Let $X=(\mathcal A,\mathcal B,M)$ consist of $A_\infty$ categories
on $V_A,V_B$ and an $A_\infty$ $(\mathcal A,\mathcal B)$-bimodule. Let
$\mathcal C$ be an $A_\infty$ category on $V$.
By Proposition~\ref{prop_change_quiver_algebras} and
Corollary~\ref{cor_diagonal_Ainfinity_bimodule},
$\Phi_f^*\mathcal C=(\mathcal C_A,\mathcal C_B,\mathcal C_\Delta)$,
where $\mathcal C_A,\mathcal C_B$ are the labelled copies of
$\mathcal C$ and $\mathcal C_\Delta$ is their diagonal $A_\infty$
bimodule, with $\mathcal C_\Delta(a_A,b_B)=\mathcal C(a,b)$.
By Proposition~\ref{prop_change_quiver_algebra_adjunction},
$(\Phi_f)_!X$ is the quotient of the free
$(\mathsf{dgCat}_V)_\infty$-algebra of
\eqref{eq_free_algebra_universal_property} on
\[
 G(a,b)=\mathcal A(a_A,b_A)\oplus M(a_A,b_B)\oplus\mathcal B(a_B,b_B)
\]
by the relations imposing the given category and bimodule operations.
Folding permits additional composites among these generators.

In summary, we obtain the following corollary.

\begin{cor}[Adjunction for diagonal bimodules]
\label{cor_diagonal_algebra_adjunction}
There are bijections, natural in $X$ and $\mathcal C$,
\begin{equation*}\operatorname{Hom}_{\operatorname{Alg}((\mathsf{dgCat}_V)_\infty)}
 \bigl((\Phi_f)_!X,\mathcal C\bigr) \cong
 \operatorname{Hom}_{\operatorname{Alg}((\mathscr B\mathsf{Mod}_{V_A,V_B})_\infty)}
 \bigl(X,(\mathcal C_A,\mathcal C_B,\mathcal C_\Delta)\bigr).
\end{equation*}
\end{cor}

The left-hand Hom set consists of strict $A_\infty$ functors
$(\Phi_f)_!X\to\mathcal C$ fixing the objects $V$.
The right-hand Hom set consists of strict $A_\infty$ functors
$F_A:\mathcal A\to\mathcal C_A$ and $F_B:\mathcal B\to\mathcal C_B$
fixing the labelled objects, together with a strict morphism of
$A_\infty$ $(\mathcal A,\mathcal B)$-bimodules
$M\to(F_A,F_B)^*\mathcal C_\Delta$.
Here the target actions are restricted along $F_A,F_B$, and
strictness means that the bimodule morphism consists of cochain
maps preserving all bimodule operations.

\subsection{Tensor products by extension of algebras}
\label{sss_derived_tensor_construction}

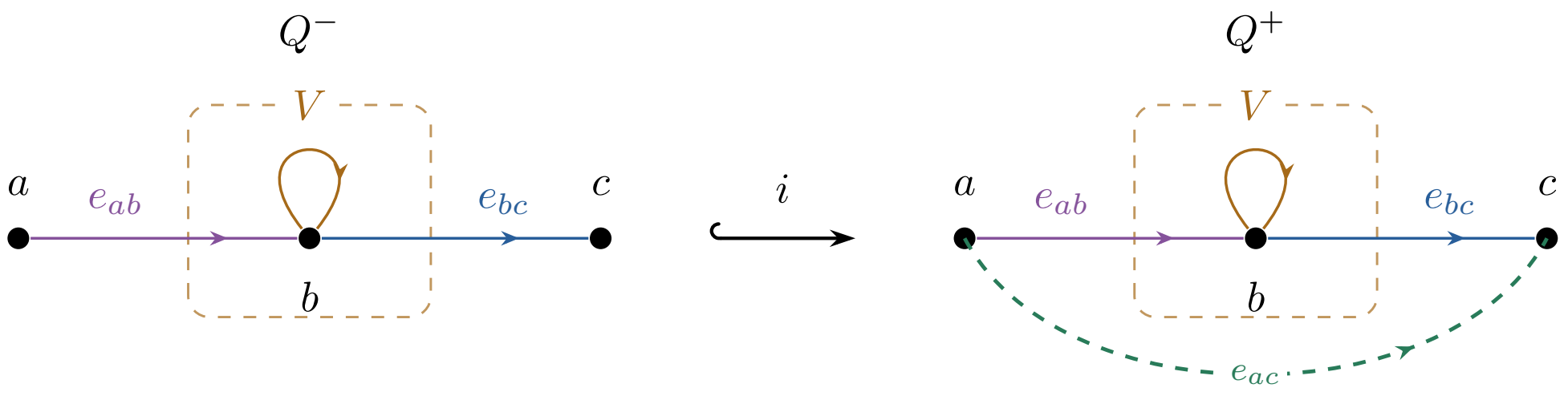
\begin{figure}[htbp]
\centering
\begingroup
\definecolor{barcategory}{HTML}{A66A19}
\definecolor{barleft}{HTML}{85529B}
\definecolor{barright}{HTML}{285E9A}
\definecolor{baroutput}{HTML}{287B59}
\begin{tikzpicture}[
  x=1cm,y=1cm,font=\small,line width=.65pt,
  arr/.style={postaction={decorate},decoration={markings,
    mark=at position .74 with {\arrow{Stealth[length=4.2pt,width=3.4pt]}}}},
  B/.style={arr,draw=barcategory},
  M/.style={arr,draw=barleft},
  N/.style={arr,draw=barright},
  vertex/.style={circle,fill=black,inner sep=1.8pt},
  lab/.style={fill=white,inner sep=1.5pt,font=\scriptsize}
]

% The two panels have the same vertices and all the same old arrows.
\foreach \panelshift in {0,7.8} {
  \begin{scope}[xshift=\panelshift cm]
    \draw[draw=barcategory!70,dashed,rounded corners=5pt,line width=.55pt]
      (1.4,-.65) rectangle (3.4,1.1);
    \node[fill=white,inner xsep=3pt,inner ysep=2pt,text=barcategory]
      at (2.4,1.1) {$V$};

    \node[vertex,label={[label distance=4pt]90:$a$}] (a) at (0,0) {};
    \node[vertex,label={[label distance=4pt]270:$b$}] (b) at (2.4,0) {};
    \node[vertex,label={[label distance=4pt]90:$c$}] (c) at (4.8,0) {};

    % One middle object and its endomorphism loop are shown.
    \draw[B] (b) .. controls (1.65,.95) and (3.15,.95) .. (b);

    \draw[M] (a) -- (b);
    \draw[N] (b) -- (c);
    \node[text=barleft] at (.8,.3) {$e_{ab}$};
    \node[text=barright] at (4,.3) {$e_{bc}$};
  \end{scope}
}

\node at (2.4,1.7) {$Q^-$};
\node at (10.2,1.7) {$Q^+$};
\draw[right hook-{Stealth[length=6pt,width=4pt]},line width=.8pt]
  (5.7,0) -- node[above=5pt] {$i$} (6.9,0);

% Only the completed quiver has the composite output arrow.
\draw[arr,draw=baroutput,dashed,line width=.9pt]
  (7.8,0) .. controls (8.6,-1.5) and (11.8,-1.5) ..
  node[lab,pos=.5,text=baroutput] {$e_{ac}$} (12.6,0);
\path[use as bounding box] (-.7,-1.35) rectangle (13.3,1.9);
\end{tikzpicture}
\endgroup
\caption{The inclusion $i:Q^-\hookrightarrow Q^+$ adds the green
output edge $e_{ac}$. One middle vertex $b\in V$ and its loop
are shown; neither endpoint has a loop.}
\label{fig:bar-tensor-quiver}
\end{figure}

Let $V$ be a set and let $a,c$ be two auxiliary vertices.
Let $Q^-$ have vertex set $\{a\}\sqcup V\sqcup\{c\}$ and arrows
$ e_{ab}:a\to b$, $e_{bb'}:b\to b'$, $e_{bc}:b\to c$ for
 $b,b'\in V$.
Thus $Q^-$ joins the right- and left-module quivers of
\S~\ref{sss_dg_category_bimodules} along their common complete
quiver on $V$. Let $Q^+$ be obtained by adjoining $e_{ac}:a\to c$,
and define the inclusion and the terminal map
(cf. Figure~\ref{fig:bar-tensor-quiver})
\[
 i:Q^-\hookrightarrow Q^+,
 \qquad q:Q^+\longrightarrow *.
\]

At the operad level, extension along $i$ supplies the missing
composites from $a$ to $c$. Set
\[
 \mathscr P=q^*\mathsf{As},\qquad
 \mathscr O=i^*\mathscr P,\qquad
 \Phi=(i,\mathrm{id}):(Q^-,\mathscr O)\longrightarrow(Q^+,\mathscr P).
\]
The quiver $Q^+$ is transitive and thin, so $q$ has unique diagram
lifting for diagrams without nullary faces by
Example~\ref{ex_change_quiver_terminal_folding}(i).
The same holds for $i$: a path with output in $Q^-$ cannot cross
from $a$ to $c$, so none of its diagram edges can be $e_{ac}$.
Since $\mathsf{As}$ is Koszul,
using Theorem~\ref{thm_change_quiver_koszul} shows that
$\mathscr P$ and $\mathscr O$ are Koszul and gives canonical
isomorphisms
\begin{equation}\label{eq_tensor_cooperad_pullback}
 \mathscr P^{\antishriek}
 \cong q^*(\mathsf{As}^{\antishriek}),\qquad
 \mathscr O_\infty\cong i^*\mathscr P_\infty.
\end{equation}
The second isomorphism defines a morphism
$\Phi_\infty=(i,\mathrm{id})$ between the Koszul resolutions.

The algebras over these four operads have the following descriptions.
\begin{itemize}
\item An $\mathscr O$-algebra consists of a non-unital dg category
$\mathcal B$ on $V$, a right dg $\mathcal B$-module $M$, and a left
dg $\mathcal B$-module $N$.
\item A $\mathscr P$-algebra consists of such a triple
$(\mathcal B,M,N)$, a cochain complex $T$ on $e_{ac}$, and a cochain
map $M\otimes_{\mathcal B}N\to T$, where
$\otimes_{\mathcal B}$ denotes the ordinary balanced tensor product.
\item An $\mathscr O_\infty$-algebra consists of a non-unital
$A_\infty$ category $\mathcal B$ on $V$, a right $A_\infty$
$\mathcal B$-module $M$, and a left $A_\infty$ $\mathcal B$-module $N$.
\item A $\mathscr P_\infty$-algebra consists of such an
$A_\infty$ triple, a cochain complex $T$ on $e_{ac}$, and higher
pairings
\[
 \beta_r^{b_0\cdots b_r}:
 M(b_0)\otimes\mathcal B(b_0,b_1)\otimes\cdots
 \otimes\mathcal B(b_{r-1},b_r)\otimes N(b_r)
 \longrightarrow T,\qquad r\geq0,
\]
of degree $-r$, satisfying the $A_\infty$ relations on strings
$a\to b_0\to\cdots\to b_r\to c$.
\end{itemize}
In each case, $\mathcal B(b,b')$, $M(b)$, and $N(b)$ are the
complexes on $e_{bb'}$, $e_{ab}$, and $e_{bc}$, respectively.
The restriction functors $\Phi^*$ and $\Phi_\infty^*$ forget $T$
and its pairing operations.

For the tensor-product construction below, fix
$X=(\mathcal B,M,N)$ and assume that the category and both modules
are strictly unital. Thus $X$ is an $\mathscr O_\infty$-algebra,
or an $\mathscr O$-algebra in the dg case.

Throughout, $M\otimes_{\mathcal B}^{\mathbb L}N$ denotes the
two-sided bar model of the derived tensor product,
with underlying graded space
\begin{equation}\label{eq_tensor_bar_model}
\begin{aligned}
 M\otimes_{\mathcal B}^{\mathbb L}N
 =\bigoplus_{\substack{r\ge0\\b_0,\ldots,b_r\in\Ob\mathcal B}}
M(b_0)\otimes\mathcal B(b_0,b_1)[1]\otimes\cdots \otimes\mathcal B(b_{r-1},b_r)[1]\otimes N(b_r).
\end{aligned}
\end{equation}
Here $K[j]^i=K^{i+j}$ with differential $(-1)^j d_K$, so
$K[1]=s^{-1}K$.

\Needspace{14\baselineskip}
\begin{prop}[Tensor products by extension]
\label{prop_tensor_algebra_extension}
For $X$ as above, extension of
algebras has the following properties.
\begin{enumerate}[(i)]
\item For dg data, the adjunction unit $X\to\Phi^*\Phi_!X$ is an
isomorphism, and $(\Phi_!X)(e_{ac})$ is canonically the ordinary
balanced tensor product $M\otimes_{\mathcal B}N$.
\item For $A_\infty$ data, the adjunction unit
$X\to\Phi_\infty^*(\Phi_\infty)_!X$ is an isomorphism, and there
is a canonical isomorphism of complexes
\begin{equation}\label{eq_tensor_extension_bar}
 \bigl((\Phi_\infty)_!X\bigr)(e_{ac})
 \cong M\otimes_{\mathcal B}^{\mathbb L}N.
\end{equation}
\end{enumerate}
These are the initial extensions of $X$, with morphisms fixing $X$,
in the corresponding categories of operadic algebras.
\end{prop}
\begin{proof}
No arrow enters $a$ or leaves $c$, so no operation of arity at
least two accepts $e_{ac}$. The free-algebra quotient of
Proposition~\ref{prop_change_quiver_algebra_adjunction} therefore
recovers $X$ on $Q^-$, proving the unit assertions.
Let $R_X$ and $R_X^\infty$ be the $e_{ac}$-components of the dg
relation ideals \eqref{eq_change_quiver_algebra_relations}
for $\Phi$ and $\Phi_\infty$, respectively.

For (i), the relations evaluate consecutive input segments with
output in $Q^-$. Associativity reduces each word to $M(b)\otimes N(b)$;
the two evaluations of $m\otimes u\otimes n$ give precisely the
balancing relation. Hence
\begin{align*}
 (\Phi_!X)(e_{ac})
 &\cong
 \left(\bigoplus_{\substack{r\ge0\\b_0,\ldots,b_r\in V}}
 M(b_0)\otimes\bigotimes_{j=1}^r\mathcal B(b_{j-1},b_j)
 \otimes N(b_r)\right)\bigg/R_X\\
 &\cong
 \frac{\displaystyle\bigoplus_{b\in V}M(b)\otimes N(b)}
 {\bigl\langle (m u)\otimes n-m\otimes(u n)
 \mid m\in M(b),\ u\in\mathcal B(b,b'),\ n\in N(b')\bigr\rangle} =M\otimes_{\mathcal B}N.
\end{align*}
These identifications commute with the differentials because
the compositions and module actions are cochain maps.

For (ii), write $\theta_{\mathbf b}=s\nu_{a b_0\cdots b_r c}$,
of degree $-r$, labelling the elementary face with input path
$a\to b_0\to\cdots\to b_r\to c$ and output $e_{ac}$.
Since $(i_\#X)(e_{ac})=0$, every contributing planar diagram
has a root face with output $e_{ac}$ and all other face profiles in $Q^-$.
Modulo $R_X^\infty$, evaluating the subdiagrams attached to the
root face leaves that face with inputs in $X$.

Freeness of the underlying graded operad therefore gives
\begin{align*}
 \bigl((\Phi_\infty)_!X\bigr)(e_{ac})
 &\cong
 \frac{\mathbb T_{\mathscr P_\infty}(i_\#X)(e_{ac})}{R_X^\infty}\\
 &\cong
 \bigoplus_{r,\mathbf b}\Bbbk\theta_{\mathbf b}\otimes M(b_0)
 \otimes\bigotimes_{j=1}^r\mathcal B(b_{j-1},b_j)\otimes N(b_r)\\
 &\xrightarrow[\;\sigma\;]{\cong}
 \bigoplus_{r,\mathbf b}M(b_0)
 \otimes\bigotimes_{j=1}^r\mathcal B(b_{j-1},b_j)[1]\otimes N(b_r) =M\otimes_{\mathcal B}^{\mathbb L}N,
\end{align*}
as graded spaces, where $r\ge0$ and $\mathbf b=(b_0,\ldots,b_r)\in V^{r+1}$; suppressing object labels, the suspension isomorphism is
\begin{equation}\label{eq_tensor_suspension_sign}
\begin{gathered}
 \sigma:\quad
 \theta_r\otimes m\otimes\beta_1\otimes\cdots\otimes\beta_r\otimes n
 \longmapsto
 (-1)^{\epsilon_r}m\otimes s^{-1}\beta_1\otimes\cdots
 \otimes s^{-1}\beta_r\otimes n
\end{gathered}
\end{equation}
where the sign is $\epsilon_r=r|m|+\sum_{j=1}^r(r-j)|\beta_j|$.
The cobar differential inserts an internal edge, splitting off
a face on an interior, initial, or final input segment.
Evaluating that face gives a category, right-module, or left-module
operation. Together with the internal differentials, these are
the bar terms, with signs transported by
\eqref{eq_cobar_d2_explicit} and \eqref{eq_tensor_suspension_sign}.
Thus $\sigma$ is an isomorphism of complexes.
Initiality follows from the adjunction and its invertible unit.
\end{proof}

\begin{rmk}
	For strictly unital
$(\mathcal A,\mathcal B)$- and $(\mathcal B,\mathcal D)$-bimodules
$M,N$ over strictly unital $A_\infty$ categories, the objectwise bar
construction gives
$ T(a,c)=M(a,-)\otimes_{\mathcal B}^{\mathbb L}N(-,c)$.
The left actions with $p>0$ sum $n^M_{p,j}$ on initial interior
segments, and the right actions with $q>0$ sum $n^N_{j,q}$ on final
segments; operations with $p,q>0$ vanish. With the suspension signs
of \eqref{eq_bimodule_suspension_convention}, the identities for
$M,N$ make $T$ a strictly unital $(\mathcal A,\mathcal D)$-bimodule,
denoted $M\otimes_{\mathcal B}^{\mathbb L}N$; see
\cite[\S~3.1.2]{Ganatra_Pardon_Shende_sectors}.
The diagonal bimodules give canonical quasi-isomorphisms of
$A_\infty$ bimodules
$M\otimes_{\mathcal B}^{\mathbb L}\mathcal B_\Delta\xrightarrow{\sim}M$
and
$\mathcal B_\Delta\otimes_{\mathcal B}^{\mathbb L}N\xrightarrow{\sim}N$;
see \cite[Lemma~3.7]
{Ganatra_Pardon_Shende_sectors}.
\end{rmk}

For $X=(\mathcal B,M,N)$ as above, write
\[
\begin{aligned}
 E(X)&=(\mathcal B,M,N,M\otimes_{\mathcal B}N),\\
 E_\infty(X)&=(\mathcal B,M,N,M\otimes_{\mathcal B}^{\mathbb L}N),
\end{aligned}
\]
with the canonical $\mathscr P$- and $\mathscr P_\infty$-algebra
structures of Proposition~\ref{prop_tensor_algebra_extension}
in the dg and $A_\infty$ cases, respectively.
Each extension retains the category $\mathcal B$, its right module
$M$, and its left module $N$, and places the indicated tensor-product
complex on the new arrow $e_{ac}$, with the canonical balanced
pairing or higher pairings, respectively.

\begin{cor}[Tensor-product adjunctions]
\label{cor_tensor_algebra_adjunction}
There are natural bijections
\begin{enumerate}[(i)]
\item for dg data $X$ and $Y\in\operatorname{Alg}(\mathscr P)$,
\begin{equation}\label{eq_tensor_dg_algebra_adjunction}
 \operatorname{Hom}_{\operatorname{Alg}(\mathscr P)}\bigl(E(X),Y\bigr)
 \cong
 \operatorname{Hom}_{\operatorname{Alg}(\mathscr O)}\bigl(X,\Phi^*Y\bigr);
\end{equation}
\item for $A_\infty$ data $X$ and
$Y\in\operatorname{Alg}(\mathscr P_\infty)$,
\begin{equation}\label{eq_tensor_Ainfinity_algebra_adjunction}
 \operatorname{Hom}_{\operatorname{Alg}(\mathscr P_\infty)}
 \bigl(E_\infty(X),Y\bigr)
 \cong
 \operatorname{Hom}_{\operatorname{Alg}(\mathscr O_\infty)}
 \bigl(X,\Phi_\infty^*Y\bigr).
\end{equation}
\end{enumerate}
\end{cor}
\begin{proof}
Proposition~\ref{prop_tensor_algebra_extension} identifies
$E(X)\cong\Phi_!X$ and
$E_\infty(X)\cong(\Phi_\infty)_!X$ as operadic algebras, with
their restrictions identified with $X$ by the adjunction units.
Apply Proposition~\ref{prop_change_quiver_algebra_adjunction}
to $\Phi$ and $\Phi_\infty$, respectively.
\end{proof}

\bibliographystyle{abbrv}

	\bibliography{mybib_haho}

\end{document}